\documentclass[ 10pt,reqno]{amsart}
\usepackage{amsaddr}
\usepackage[a4paper, total={6.5in, 8.4in}]{geometry}

\usepackage[symbol]{footmisc}

\usepackage{amsmath}
\usepackage{amssymb}
\usepackage{amsthm}
\usepackage[latin1]{inputenc}
\usepackage{eurosym}
\usepackage{graphicx}
\usepackage{epsfig}
\usepackage{hyperref}
\usepackage{dsfont}
\usepackage[space]{cite}
\usepackage{subcaption}
\usepackage{caption}
\usepackage{placeins}
\usepackage[ruled,vlined]{algorithm2e}
\usepackage[title]{appendix}

\usepackage[displaymath,mathlines]{lineno}
\modulolinenumbers[1]

\allowdisplaybreaks

\usepackage{hyperref}
\usepackage{ifthen}
\usepackage{comment}

\makeindex
\usepackage{color}
\usepackage{float}
\usepackage{commath}

\renewcommand{\div}{\operatorname{div}}

\newcommand{\epsi}{\varepsilon}

\def\leq{\leqslant}
\def\geq{\geqslant}

\numberwithin{equation}{section}

\newtheoremstyle{thmlemcorr}{10pt}{10pt}{\itshape}{}{\bfseries}{.}{10pt}{{\thmname{#1}\thmnumber{
#2}\thmnote{ (#3)}}}
\newtheoremstyle{thmlemcorr*}{10pt}{10pt}{\itshape}{}{\bfseries}{.}\newline{{\thmname{#1}\thmnumber{
#2}\thmnote{ (#3)}}}
\newtheoremstyle{defi}{10pt}{10pt}{\itshape}{}{\bfseries}{.}{10pt}{{\thmname{#1}\thmnumber{
#2}\thmnote{ (#3)}}}
\newtheoremstyle{remexample}{10pt}{10pt}{}{}{\bfseries}{.}{10pt}{{\thmname{#1}\thmnumber{
#2}\thmnote{ (#3)}}}
\newtheoremstyle{ass}{10pt}{10pt}{}{}{\bfseries}{.}{10pt}{{\thmname{#1}\thmnumber{
A#2}\thmnote{ (#3)}}}

\theoremstyle{thmlemcorr}
\newtheorem{theorem}{Theorem}
\numberwithin{theorem}{section}

\newtheorem{corollary}[theorem]{Corollary}

\theoremstyle{thmlemcorr*}
\newtheorem{theorem*}{Theorem}
\newtheorem{lemma*}[theorem]{Lemma}
\newtheorem{corollary*}[theorem]{Corollary}
\newtheorem{proposition*}[theorem]{Proposition}
\newtheorem{problem*}[theorem]{Problem}
\newtheorem{conjecture*}[theorem]{Conjecture}

\theoremstyle{defi}

\newtheorem{hyp}{Assumption}
\newtheorem{claim}{Claim}

\theoremstyle{remexample}
\newtheorem{remark}[theorem]{Remark}

\theoremstyle{ass}

\newtheorem*{notations*}{Notations}

\allowdisplaybreaks

\begin{document}

\title[]{Numerical methods for  Mean field Games based on Gaussian Processes and Fourier Features}

\author{Chenchen Mou$^1$, Xianjin Yang$^{2,3,*}$, Chao Zhou$^4$}

\address[C. MOU]{
	$^1$Department of Mathematics, City University of Hong Kong, Hong Kong SAR, China.}
\email{chencmou@cityu.edu.hk}

\thanks{$^*$Corresponding author.}
\address[X. Yang]{$^2$Yau Mathematical Sciences Center, Tsinghua  University, Haidian  District, Beijing, 100084, China.}
\address{$^3$Beijing Institute of Mathematical Sciences and Applications, Huairou District, Beijing, 101408, China.}
\email{yxjmath@gmail.com}

\address[C. Zhou]{$^4$
	Department of Mathematics and Risk Management Institute, National University of Singapore, Singapore.}
\email{matzc@nus.edu.sg}

 %Existence problems: global existence, local existence, non-existence

\thanks{
}
\date{\today}

\begin{abstract}
In this article, we propose two numerical methods, the Gaussian Process (GP) method and the Fourier Features (FF) algorithm, to solve mean field games (MFGs). The GP  algorithm  approximates the solution of a MFG with maximum \textit{a posteriori} probability estimators of GPs conditioned on the partial differential equation (PDE) system of the MFG at a finite number of sample points.  The main bottleneck of the GP method is to compute the inverse of a square gram matrix, whose size is proportional to the number of sample points. To improve the performance, we introduce the FF method, whose insight comes from the recent trend of approximating positive definite kernels with random Fourier features. The FF algorithm seeks approximated solutions in the space  generated by sampled Fourier features. In the FF method, the size of the matrix to be inverted depends only on the number of Fourier features selected, which is much less than the size of sample points. Hence, the FF method reduces the precomputation time, saves the memory, and achieves comparable accuracy to the GP method. We give the existence and the convergence proofs for both algorithms. The convergence argument of the GP method does not depend on any monotonicity condition, which suggests the potential applications of the GP method to solve MFGs with non-monotone couplings in future work. We show the efficacy of our algorithms through experiments on a stationary MFG with a non-local coupling and on a time-dependent planning problem. We believe that the FF method can also serve as an alternative algorithm to solve general PDEs. 
\end{abstract}

\maketitle 
\section{Introduction}
Mean field games (MFGs)  \cite{lasry2006jeux1, lasry2006jeux, lasry2007mean, huang2006large, huang2007large, huang2007nash, huang2007invariance} study the behavior of a large population of rational and indistinguishable agents as the number of agents goes to infinity. The Nash equilibrium of a typical MFG
is formulated by a coupled system of two partial differential equations (PDEs), a Hamilton--Jacobi--Bellman
(HJB) equation and a Fokker--Plank (FP) equation. The HJB equation gives the value function of agents and the FP equation determines agents' distribution. Recently, MFG models have found widespread applications, see \cite{gueant2011mean, gomes2015economic, lee2020mean, gao2021belief, gomes2021mean, evangelista2018first, lee2021mean}.
\begin{comment}
\cite{gueant2011mean, achdou2014partial, gomes2015economic, achdou2017income, lachapelle2016efficiency, cardaliaguet2018mean, casgrain2019algorithmic, firoozi2017optimal, lachapelle2011mean, burger2013mean, aurell2018mean, achdou2019mean, kizilkale2019integral, de2019mean, gomes2021mean, weinan2019mean, welch2019describing, evangelista2017radially, evangelista2018first}.
\end{comment}
Up to now, the well-posedness of MFGs is well understood in various settings and
the first results date back to the original works of Lasry and Lions and have been given
in the course of Lions at Coll{\`e}ge de France (see \cite{lionslec}). 
However, very few MFG models admit explicit solutions. Hence, numerical computations of MFGs play an essential role in obtaining quantitative descriptions of underlying models. 

Here, we propose two new algorithms to solve MFGs. The first one, called the Gaussian Process (GP) method, applies the algorithm in \cite{chen2021solving} to solve MFGs. The authors in \cite{chen2021solving} give a GP regression framework to solve nonlinear PDEs. The solution to a PDE is found by solving an optimal recovery problem, whose minimizer is viewed as a maximum \textit{a posteriori} probability (MAP) estimator of a GP conditioned on the PDE evaluated at sample points. The main bottleneck of the GP method lies in the computation of the inverse of a square gram matrix, whose size is proportional to the product of the size of sample points times the number of  linear operators in the PDE. To improve the performance, we propose the Fourier Features (FF) method, where the optimal recovery problem seeks minimizers in the space generated by sampled trigonometric functions. The FF method is based on the recent trend of approximating positive definite kernels with randomized trigonometric functions in Gaussian regressions, see \cite{rahimi2007random, yu2016orthogonal, tompkins2020periodic}. Using the technique we propose in Remark \ref{AAqr}, the dimension of the matrix needed to be inverted in the FF algorithm depends only on the number of sampled Fourier features, which is less than the size of sample points. The numerical experiments show that the FF algorithm reduces the precomputation time and the amount of storage with comparable accuracy compared to the GP method. 

Meanwhile, we also prove the convergence of our algorithms. The proofs for the GP method are based on the compactness arguments in \cite{chen2021solving} and do not depend on any  monotonicity condition that guarantees stability and uniqueness of the solution to a typical MFG. This feature implies the potential applications of the GP method to MFGs with non-monotone couplings. On the other hand, since the Fourier features space the FF method uses lacks compactness, the same  arguments of the GP method cannot be adapted to the setting of the FF method. Instead, the Lasry--Lions monotonicity arguments provide the uniform bounds for the errors of numerical solutions and lead to the convergence of the FF algorithm. In future work, we plan to investigate the convergence of the FF method under the setting of MFGs without monotone couplings. 
\subsection{Related Works}
\label{relatedworks}
By now, there have been various numerical methods for MFGs. Here, we briefly summarize the numerical algorithms solving MFGs that are closely related to the methods proposed in this paper. For clarification, we group them  into different categories. The first group consists of mesh-based algorithms:
\begin{itemize}
	\item Finite difference methods. As far as we know, the first finite difference method for MFGs is introduced in  \cite{achdou2010mean} for both stationary and time-dependent MFGs. The Lax-Friedrichs or Godunov type schemes are introduced to approximate Hamiltonians. The Fokker--Plank equation is discretized in such a way that it preserves the adjoint structure in the MFG. The authors proved the existence and uniqueness of the discretized algebraic systems. One can solve the discretized equations using iterative methods \cite{achdou2012iterative}. To know more about details and applications of the finite difference method, we refer readers to \cite{ lauriere2021numerical, achdou2020mean, achdou2012mean, achdou2019mean};
	\item Optimization algorithms. In the seminal paper \cite{lasry2007mean}, Lasry and Lions explain  the concept of a variational MFG, which can be interpreted as the optimality condition of  a PDE-driven optimization problem. From this point of view, several techniques from Optimization have been introduced to solve the PDE constrained minimization problem. In \cite{benamou2015augmented}, the authors propose augmented Lagrangian methods. Later, the performances of different algorithms are compared in  \cite{briceno2018proximal, briceno2019implementation}, and the Chambolle-Pock algorithm outperforms other methods. 
	\item Monotone flows. In \cite{almulla2017two}, the authors have proposed the monotone flow method to solve stationary MFGs. Later, in \cite{gomes2021numerical}, the authors apply the monotone flow to solve  stationary and time-dependent MFGs with finite states. The insight is that solving a MFG with monotonicity is equivalent to computing a zero of a monotone operator, which is the stationary point of the monotone flow. However, the monotone flow may not guarantee the non-negativity of the probability density in the flow. In \cite{gomes2020hessian}, the authors design the Hessian Riemannian flow method to preserve positivity. 
\end{itemize}
Despite the rigorous theories behind the algorithms mentioned above, they are meshed-based methods and prone to the curse of dimensionality. The following methods are mesh-free methods developed recently to keep pace with growing problem sizes: 
\begin{itemize} 
	\item Lagrangian Methods. As far as we know, the first Lagrangian method for MFGs appears in \cite{nurbekyan2019fourier}. The authors use Fourier expansions of the non-local coupling term in the HJB equation and parameterize the MFG. The parameterized MFG is reformulated as the optimality condition of a convex optimization problem  over a finite-dimensional subspace of continuous curves, which is independent of the dimension. Based on this work, the authors in \cite{liu2020splitting, liu2021computational} further develop this idea to solve non-potential MFGs with mixed non-local couplings. Even though the methods proposed in our paper use kernels and Fourier features, we work in a different way and approximate the solution of a MFG in reproducing kernel Hilbert spaces (RKHSs) or Fourier feature spaces. 
	\item Neural Networks. Since neural networks are compositions of nonlinear maps and have the power to represent sufficient complicated functions, the authors in \cite{carmona2019convergence, carmona2021convergence} use neural networks to approximate unknowns in PDEs and solve  ergodic MFGs and time-dependent MFGs. At the same time, in \cite{ruthotto2020machine}, the authors use the Lagrangian Method in \cite{nurbekyan2019fourier} to reformulate variational MFGs and solve the resulting system using methods from neural networks. For more applications of neural networks in MFGs, see \cite{lin2020apac, carmona2021deep}.  
\end{itemize}  
We refer readers to the recent comprehensive surveys \cite{lauriere2021numerical, achdou2020mean} for  more details and applications of numerical methods for MFGs. 

Compared to the above-mentioned algorithms, our methods for MFGs admit the following features:
\begin{enumerate}
	\item The GP and the FF algorithms are meshfree and flexible to the shape of  domains, compared to algorithms in the first group mentioned above. Especially, when we solve MFGs in the whole Euclidean space, to approximate derivatives,  mesh-based algorithms have to impose artificial conditions on the boundary of the domain we choose to work on. Our methods parameterize the values of derivatives evaluated at sample points and solve reformulated finite-dimensional minimization problems as in \eqref{1ssmfgoptfnpk} and \eqref{1doptffrcfd}. Hence, we do not impose extra conditions to deal with derivatives near the boundary;
	\item Compared to the neural network methods, our algorithms base on theories of RKHSs and Fourier series, for which the math backgrounds are well understood. On the other hand, our algorithms are equivalent to the neural network methods in the following perspective: a neural network with a single inner layer, with the activation functions being feature maps parameterized at sample points or being random sampled trigonometric functions, and with a linear  output layer can be viewed as a function in RKHSs or the Fourier features space, and vice versa; 
	\item The choice of the kernel for the GP method and the selection of Fourier features have a profound impact on the convergence and the accuracy of our approximations. We leave the study of hyperparameter learning (see \cite{owhadi2019kernel, chen2021consistency}) to future work.
	\item The convergence of the GP method does not depend on the Lasry--Lions monotonicity condition of the coupling terms in MFGs, which is required by the monotone flow 	
	\cite{almulla2017two, gomes2021numerical, gomes2020hessian}
	and the Lagrangian methods \cite{nurbekyan2019fourier, liu2020splitting, liu2021computational}. Hence, we plan to study the application of the GP method to solve MFGs with displacement monotonicity \cite{meszaros2021mean} or non-monotone couplings in future work. 
	\item In general, it is less costly to compute a linear map of a Fourier feature than to calculate the same linear transformation of a kernel function. For instance, in Subsection \ref{secErgodicMFG}, we need to parametrize the linear operator $L=(1-\Delta)^{-1}(1-\Delta)^{-1}$ at sample points, where $\Delta$ is a Laplacian operator. In the GP method, the representer theorem gives expressions involving the computation of $LK_2(x,\cdot)$ for $x\in \mathbb{T}^2$, where $K_2$ is the kernel of the RKHS we choose to find an approximation for the probability density. Since $LK_2(x,\cdot)$ does not admit an explicit formula, one can use the fast Fourier transform to compute it numerically. On the other hand, for any $\omega\in \mathbb{N}^2$, we observe that $L\sin(\omega^Tx)=\sin(w^Tx)/(1+|\omega|^2)^2$ and $L\cos(\omega^Tx)=\cos(w^Tx)/(1+|\omega|^2)^2$ for $x\in \mathbb{T}^2$. Hence, if we choose trigonometric functions as features in the FF method, the action of $L$ on a Fourier feature admits an explicit formula and is easier to compute.  
	\item The main bottleneck of the GP method is to compute the Cholesky decomposition of the gram matrix, whose size increases as the number of sample points grows. The FF method  relieves the pressure by approximating functions in the FF space and by using the technique in Remark \ref{AAqr}. However, as the dimension of the problem at hand increases, we should also enlarge the FF space in the absence of information about the solution. Hence, a clever selection of features is required to apply the FF method in large dimensions. In this paper, we choose the Fourier series in the periodic settings and use random Fourier features for the non-periodic cases. We leave the study of other selections to future work. 
	
\end{enumerate}

\subsection{Outline}
\label{outline}
This article is organized as follows. In Section \ref{secApp}, we present our algorithms by solving a one-dimensional MFG which admits an explicit smooth solution. We also give the existence and convergence analysis of our methods. A simple numerical experiment follows and verifies the efficacy of our methods. Section \ref{secGF} presents the general frameworks of the GP method and the FF algorithm. We also show the existence and convergence for the GP and the FF methods in general settings. The proofs related to the GP method are shown in Section \ref{secGF} and those of the FF method are presented in Appendix \ref{appendix}. Numerical experiments follow in Section \ref{secNumerical}. Conclusions and future work appear in Section \ref{secConft}.

\begin{notations*}For a vector $v$ with real values, we denote by $|v|$ the Euclidean norm of $v$ and by $v^T$ the transpose of $v$. We represent the inner product of two real valued vectors $u$ and $v$ by $\langle u, v\rangle$ or $u^Tv$. Let $\mathbb{Z}^d$ be the space of $d$-dimensional vectors with integer elements. Given $N\in\mathbb{N}$, we define $\mathbb{Z}^d_N$ as the set $\{i|i\in \mathbb{Z}^d, 0<|i|\leq N\}$. Let $\Omega$ be a subset of $\mathbb{R}^d$, let $\operatorname{int}\Omega$ be the interior of $\Omega$, and let $\partial\Omega$ be the boundary of $\Omega$.  We denote by $L^2(\Omega)$ the space of square-integrable functions on $\Omega\subset\mathbb{R}^d$. Meanwhile, let $H^1(\Omega)$ be the space such that $\forall f\in H^1(\Omega)$, $f\in L^2(\Omega)$ and $\nabla f\in (L^2(\Omega))^d$. Furthermore, we represent the space of functions with finite sup norm by $L^\infty(\Omega)$.   For a normed vector space $V$, we denote by $\|\cdot\|_V$ the norm of $V$. 
Let $\mathcal{U}$ be a Banach Space endowed with a quadratic norm $\|\cdot\|_{\mathcal{U}}$. We denote by $\mathcal{U}^*$ the dual of $\mathcal{U}$ and by $[\cdot, \cdot]$ the duality pairing. We assume that there exists a covariance operator $\mathcal{K}_{\mathcal{U}}: \mathcal{U}^*\mapsto \mathcal{U}$, which is linear, bijective, symmetric ($[\mathcal{K}_{\mathcal{U}}\phi, \psi]=[\mathcal{K}_{\mathcal{U}}\psi, \phi]$), and positive ($[\mathcal{K}_{\mathcal{U}}\phi, \phi]>0$ for $\phi\not=0$), such that 
\begin{align*}
	\|u\|_{\mathcal{U}}^2 = [\mathcal{K}^{-1}u, u], \forall u\in \mathcal{U}. 
\end{align*}
Let $\phi_1, \dots, \phi_P$ be $P\in\mathbb{N}$ elements of $\mathcal{U}^*$ and let  $\boldsymbol{\phi}:=(\phi_1, \dots, \phi_P)$ be an element of the product space ${(\mathcal{U}^*)}^{\bigotimes P}$. Then, for $u\in \mathcal{U}$, the pairing $[\boldsymbol{\phi}, u]$ is denoted by
\begin{align*}
[\boldsymbol{\phi}, u]:=([{\phi}_1, u], \dots, [{\phi}_P, u]).
\end{align*}
Furthermore, for $\boldsymbol{u}:=(u_1,\dots,u_S)\in \mathcal{U}^{\bigotimes S}$, $S\in \mathbb{N}$, we represent by $[\boldsymbol{\phi}, \boldsymbol{u}]\in\mathbb{R}^{P\times S}$ the matrix with entries $[\phi_i, u_j]$.
Finally, we denote by $C$ a positive real number whose value  may change line by line. 
\end{notations*}

\section{An appetizer: Solving a stationary MFG with an explicit solution}\label{summaryofmethods}
\label{secApp}
To show our algorithms, we apply our methods to solve the one-dimensional stationary MFG in  \cite{almulla2017two}, which admits an explicit solution. 
Let $\mathbb{T}$ be the one-dimensional torus and be characterized by $[0, 1)$.  
Given smooth functions $V:\mathbb{T}\mapsto \mathbb{R}$ and $b:\mathbb{T}\mapsto\mathbb{R}$, we find $(u, m, \overline{H})\in C^\infty(\mathbb{T})\times C^\infty(\mathbb{T})\times \mathbb{R}$ solving 
\begin{align}
	\label{1dsmfg}
	\begin{cases}
		\frac{u_x^2}{2}+V(x)+b(x)u_x = \ln m + \overline{H}, \ \text{on}\ \mathbb{T},\\
		-(m(u_x+b(x)))_x = 0, \ \text{on}\ \mathbb{T},\\
		\int_{\mathbb{T}}m \dif x = 1, \int_{\mathbb{T}}u \dif x = 0,
	\end{cases}
\end{align}
where $u$ is the value function, $m$ presents the probability density of the population, and $\overline{H}$ is a real number. 
According to \cite{almulla2017two}, when $\int_{\mathbb{T}}b(x)\dif x=0$, \eqref{1dsmfg} admits a unique smooth solution, which has  the following explicit formulas
\begin{align}
	\label{exp1du}
	\begin{cases}
		u(x)=-\int_0^yb(y)\dif y+\int_{\mathbb{T}}\int_0^xb(y)\dif y\dif x,\\
		m(x)=\frac{e^{V(x)-\frac{b^2(x)}{2}}}{\int_{\mathbb{T}}e^{V(y)-\frac{b^2(y)}{2}}\dif y}, \\
		\overline{H}=\ln \bigg(\int_{\mathbb{T}}e^{V(y)-\frac{b^2(y)}{2}}\dif y\bigg).
	\end{cases} 
\end{align}
For ease of presentation, we denote 
\begin{align}
	\label{defH}
	H(x,p) := V(x) + \frac{p^2}{2} + b(x)p, \forall x\in \mathbb{T}, p\in \mathbb{R}. 
\end{align}
\subsection{The Gaussian Process Method}
\label{secgp1d}
Here, we use the GP method in \cite{chen2021solving} to solve \eqref{1dsmfg}.  First, we sample a collection of $M$ points $\{x_i\}_{i=1}^M$ in $[0, 1)$. Let $\mathcal{U}$ and $\mathcal{V}$ be the RKHSs associated with positive definite kernels $K_1$ and $K_2$, respectively. We denote by $\|\cdot\|_{\mathcal{U}}$ the norm of $\mathcal{U}$ and represent the norm of $\mathcal{V}$ by $\|\cdot\|_{\mathcal{V}}$. 
Following \cite{chen2021solving}, we approximate the solution $(u^*, m^*, \overline{H}^*)$  of \eqref{1dsmfg} by a minimizer of the following optimal recovery problem
\begin{align}
	\label{1ssmfgoptpk}
	\begin{cases}
		\min\limits_{(u, m ,\overline{H})\in \mathcal{U}\times\mathcal{V}\times\mathbb{R}} \|u\|_{\mathcal{U}}^2+\|m\|_{\mathcal{V}}^2+|\overline{H}|^2 + \beta\bigg|\frac{1}{M}\sum\limits_{i=1}^{M}m(x_i)-1\bigg|^2\\
		\quad\quad\quad\quad\quad\quad+\beta\bigg|\frac{1}{M}\sum\limits_{i=1}^{M}u(x_i)\bigg|^2\\
		\text{s.t.}\  H(x_i, u_x(x_i))=\ln m(x_i) + \overline{H}, \forall i = 1,\dots,M,\\
		\quad \ \  m_{x}(x_i)(u_x(x_i) + b(x_i)) + m(x_i)(u_{xx}(x_i)+b_x(x_i))  = 0, \forall i = 1,\dots,M,
	\end{cases}
\end{align}
where $\beta>0$ is a penalization parameter. 
We assume that the kernels $K_1$ and $K_2$  are properly chosen such that $\mathcal{U}\subset C^2(\mathbb{T})$ and $\mathcal{V}\subset C^1(\mathbb{T})$. Hence, the constraints of \eqref{1ssmfgoptpk} are well-defined. 
Let $(u^\dagger, m^\dagger, \overline{H}^\dagger)$ be a minimizer to \eqref{1ssmfgoptpk}, whose existence is given by Theorem \ref{exi11} later. Then, $u^\dagger$ and $m^\dagger$ can be viewed as MAP estimators of two  GPs  conditioned on the MFG at the sample points $\{x_i\}_{i=1}^M$ \cite{chen2021solving}. We note that the conditioned GPs are not Gaussian since the constraints in \eqref{1ssmfgoptpk} are nonlinear.

The key idea of the GP method is to characterize the minimizer of \eqref{1ssmfgoptpk} via a finite-dimensional representation formula using the representer theorem \cite[Sec. 17.8]{owhadi2019operator}. Following \cite{chen2021solving}, we rewrite \eqref{1ssmfgoptpk} as a two-level optimization problem
\begin{align}
	\label{1ssmfgopt2lvpk}
	\begin{cases}
		\min\limits_{\boldsymbol{z}\in\mathbb{R}^{3M},\boldsymbol{\rho}\in\mathbb{R}^{2M},\lambda\in\mathbb{R}}\begin{cases}
			\min\limits_{(u, m ,\overline{H})\in \mathcal{U}\times\mathcal{V}\times\mathbb{R}} \|u\|_{\mathcal{U}}^2+\|m\|_{\mathcal{V}}^2+|\overline{H}|^2\\
			\text{s.t.}\  u(x_i) = z_i^{(1)}, u_x(x_i) = z_i^{(2)}, u_{xx}(x_i) = z_i^{(3)},\forall i=1,\dots,M,\\
			\quad\ \  m(x_i) = \rho_i^{(1)},  m_x(x_i) = \rho_i^{(2)}, \overline{H} = \lambda, \forall i=1,\dots,M,
		\end{cases}\\
		\quad\quad\quad\quad\quad\quad+ \beta\bigg|\frac{1}{M}\sum\limits_{i=1}^{M}\rho_i^{(1)}- 1\bigg|^2+\beta\bigg|\frac{1}{M}\sum\limits_{i=1}^{M}z_i^{(1)}\bigg|^2\\
		\text{s.t.}\ H(x, z^{(2)}_i)=\ln \rho_i^{(1)} + \lambda,\forall i=1,\dots,M,\\
		\quad\ \   \rho_i^{(2)}(z_i^{(2)} + b(x_i)) + \rho_i^{(1)}(z_i^{(3)}+b_x(x_i))  = 0, \forall i=1,\dots,M,
	\end{cases}
\end{align}
where 
\begin{align}
	\label{zform}
	\boldsymbol{z}=(z^{(1)}_1,\dots, z^{(1)}_M, z^{(2)}_1,\dots, z^{(2)}_M, z^{(3)}_1,\dots, z^{(3)}_M)
\end{align}
and 
\begin{align}
	\label{rhoform}
	\boldsymbol{\rho}=(\rho^{(1)}_1,\dots, \rho^{(1)}_M, \rho^{(2)}_1,\dots, \rho^{(2)}_M).
\end{align}
Let $\delta_x$ be the Dirac delta function concentrated at $x$. We define $\phi_i^{(1)}=\delta_{x_i}$,  $\phi^{(2)}_i=\delta_{x_i}\circ \partial_x$, and $\phi^{(3)}_i=\delta_{x_i}\circ \partial_{xx}$ for $i=1,\dots,M$. Let $\boldsymbol{\phi}^{(1)}$,  $\boldsymbol{\phi}^{(2)}$, and $\boldsymbol{\phi}^{(3)}$ be the $M$-dimensional vectors with entries ${\phi}^{(1)}_i$, ${\phi}^{(2)}_i$, and ${\phi}^{(3)}_i$, separately,  and let  $\boldsymbol{\phi}$ be  the $(3M)$-dimensional vector obtained by concatenating $\boldsymbol{\phi}^{(1)}$, $\boldsymbol{\phi}^{(2)}$, and $\boldsymbol{\phi}^{(3)}$. Denote by ${\boldsymbol{\psi}}$ the  $(2M)$-dimensional vector obtained by concatenating $\boldsymbol{\phi}^{(1)}$ and $\boldsymbol{\phi}^{(2)}$. For ease of presentation, we also write  $\phi_i$ and $\psi_i$ as the $i$th components of $\boldsymbol{\phi}$ and $\boldsymbol{\psi}$, separately. 

For $k=1,2$, let $K_k(x, \boldsymbol{\phi})$ be the vector with entries $\int_{\mathbb{T}}K_k(x,x' )\phi_i(x')\dif x'$  and let $K_k(\boldsymbol{\phi}, \boldsymbol{\phi})$, called the gram matrix, be with entries $\int_{\mathbb{T}}\int_{\mathbb{T}}K_k(x, x)\phi_i(x')\phi_j(x')\dif x\dif x'$.  Similarly, we define $K_k(x, \boldsymbol{\psi})$ and $K_k(\boldsymbol{\psi}, \boldsymbol{\psi})$. By the representer theorem (see \cite[Sec. 17.8]{owhadi2019operator}), the first level optimization problem in \eqref{1ssmfgopt2lvpk} yields
\begin{align}
	\label{1drepre}
	\begin{cases}
		u(x)=\langle K_1(x,\boldsymbol{\phi}), K_1(\boldsymbol{\phi}, \boldsymbol{\phi})^{-1}\boldsymbol{z}\rangle,\\
		m(x)=\langle K_2(x,{\boldsymbol{\psi}}), K_2({\boldsymbol{\psi}}, {\boldsymbol{\psi}})^{-1}\boldsymbol{\rho}\rangle,\\
		\overline{H} =\lambda. 
	\end{cases}
\end{align}
Thus, we have
\begin{align*}
	\begin{cases}
		\|u\|_{\mathcal{U}}=\boldsymbol{z}^TK_1(\boldsymbol{\phi}, \boldsymbol{\phi})^{-1}\boldsymbol{z},\\
		\|m\|_{\mathcal{V}}=\boldsymbol{\rho}^TK_2({\boldsymbol{\psi}}, {\boldsymbol{\psi}})^{-1}\boldsymbol{\rho}. 
	\end{cases}
\end{align*}
Hence, we can formulate \eqref{1ssmfgopt2lvpk} as a finite-dimensional optimization problem
\begin{align}
	\label{1ssmfgoptfnpk}
	\begin{cases}
		\min\limits_{\boldsymbol{z}\in\mathbb{R}^{3M},\boldsymbol{\rho}\in\mathbb{R}^{2M},\lambda\in\mathbb{R}} \boldsymbol{z}^TK_1(\boldsymbol{\phi}, \boldsymbol{\phi})^{-1}\boldsymbol{z} + \boldsymbol{\rho}^TK_2({\boldsymbol{\psi}}, {\boldsymbol{\psi}})^{-1}\boldsymbol{\rho} + |\lambda|^2\\
		\quad\quad\quad\quad\quad\quad\quad + \beta\bigg|\frac{1}{M}\sum\limits_{i=1}^{M}\rho_i^{(1)}- 1\bigg|^2+\beta\bigg|\frac{1}{M}\sum\limits_{i=1}^{M}z_i^{(1)}\bigg|^2\\
		\text{s.t.}\ H(x, z^{(2)}_i)=\ln \rho_i^{(1)} + \lambda,\forall i = 1,\dots, M\\
		\quad\ \  \rho_i^{(2)}(z_i^{(2)} + b(x_i)) + \rho_i^{(1)}(z_i^{(3)}+b_x(x_i))  = 0,\forall i = 1,\dots, M.
	\end{cases}
\end{align}
\begin{comment}
Furthermore, following \cite{chen2021solving}, we can eliminate the last two linear constraints in \eqref{1ssmfgoptfnpk}. Let
\begin{align*}
\bar{\boldsymbol{z}}:=&(z_1^{(1)}, \dots, z_{M-1}^{(1)}, z_1^{(2)}, \dots, z_{M}^{(2)}, z_{1}^{(3)}, \dots, z_{M}^{(3)})
\end{align*} 
and 
\begin{align*}
\bar{\boldsymbol{\rho}}:=(\rho_1^{(1)},\dots,\rho_{M-1}^{(1)}, \rho_{1}^{(2)},\dots,\rho_{M}^{(2)}).
\end{align*}
For a little abuse of notations, we denote by
\begin{align*}
{\boldsymbol{z}}:=(z_1^{(1)}, \dots, z_{M-1}^{(1)}, -\sum_{i=1}^{M_{}-1}z_i^{(1)}, z_1^{(2)}, \dots, z_{M_{}}^{(2)}, z_{1}^{(3)}, \dots, z_{M_{}}^{(3)})
\end{align*}
and 
\begin{align*}
{\boldsymbol{\rho}}:=&(\rho_1^{(1)},\dots,\rho_{M-1}^{(1)}, M - \sum_{i=1}^{M-1}\rho_i^{(1)}, \rho_{1}^{(2)},\dots,\rho_{M}^{(2)}). 
\end{align*}
Then, we reformulate \eqref{1ssmfgoptfnpk} as
\begin{align}
\label{1ssmfgoptfnelpk}
\begin{cases}
\min\limits_{\bar{\boldsymbol{z}}\in\mathbb{R}^{3M-1}, \bar{\boldsymbol{\rho}}\in\mathbb{R}^{2M-1},\lambda\in\mathbb{R}} \boldsymbol{z}^TK_1(\boldsymbol{\phi}, \boldsymbol{\phi})^{-1}\boldsymbol{z} + \boldsymbol{\rho}^TK_2({\boldsymbol{\psi}}, {\boldsymbol{\psi}})^{-1}\boldsymbol{\rho} + |\lambda|^2\\
\text{s.t.}\ H(x, z^{(2)}_i)=\ln \rho_i^{(1)} + \lambda,\forall i=1,\dots,M, \\
\quad \ \   \rho_i^{(2)}(z_i^{(2)} + b(x_i)) + \rho_i^{(1)}(z_i^{(3)}+b_x(x_i))  = 0, \forall i=1,\dots,M. 
\end{cases}
\end{align}
\end{comment}
To deal with the nonlinear constraints in \eqref{1ssmfgoptfnpk}, we introduce a prescribed penalization parameter $\gamma>0$ and consider the following relaxation 
\begin{align}
	\label{relax1D}
	\begin{split}
		\min\limits_{{\boldsymbol{z}}\in\mathbb{R}^{3M}, {\boldsymbol{\rho}}\in\mathbb{R}^{3M},\lambda\in\mathbb{R}} \boldsymbol{z}^TK_1(\boldsymbol{\phi}, \boldsymbol{\phi})^{-1}\boldsymbol{z} + \boldsymbol{\rho}^TK_2({\boldsymbol{\psi}}, {\boldsymbol{\psi}})^{-1}\boldsymbol{\rho} + |\lambda|^2
		+ \beta\bigg|\frac{1}{M}\sum\limits_{i=1}^{M}\rho_i^{(1)}- 1\bigg|^2+\beta\bigg|\frac{1}{M}\sum\limits_{i=1}^{M}z_i^{(1)}\bigg|^2\\+\gamma \sum_{i=1}^{M}|e^{H(x,z_i^{(2)})-\lambda}-\rho_i^{(1)}|^2+\gamma \sum_{i=1}^{M}|\rho_i^{(2)}(z_i^{(2)} + b(x_i)) + \rho_i^{(1)}(z_i^{(3)}+b_x(x_i))|^2.\quad\quad\quad\quad\,
	\end{split}
\end{align} 
The problem \eqref{relax1D} is the foundation of the GP method for solving \eqref{1dsmfg}. 
We use the Gauss--Newton method to solve \eqref{relax1D},  which is detailed in Section 3 of  \cite{chen2021solving}.
\begin{remark}
\label{rmk1dpositi}
	In \eqref{relax1D}, we use the exponential form for the  constraints from the HJB equation to avoid possible numerical issues when evaluating the logarithm function. 
\end{remark}
\begin{remark}
	\label{rmk1dnugget}
	The matrices $K_1(\boldsymbol{\phi}, \boldsymbol{\phi})$ and $K_2(\boldsymbol{\psi}, \boldsymbol{\psi})$ are ill-conditioned in general. To compute $K_1(\boldsymbol{\phi}, \boldsymbol{\phi})^{-1}$ and $K_2(\boldsymbol{\psi}, \boldsymbol{\psi})^{-1}$, we perform the  Cholesky decomposition on $K_1(\boldsymbol{\phi}, \boldsymbol{\phi})+\eta_1R_1$ and $K_2(\boldsymbol{\psi}, \boldsymbol{\psi})+\eta_2R_2$, where $\eta_1>0$, $\eta_2>0$ are  chosen regularization  constants, and $R_1$, $R_2$ are block diagonal nuggets constructed using the approach introduced in \cite{chen2021solving}. In the numerical experiments, we precompute the Cholesky decomposition of  $K_1(\boldsymbol{\phi}, \boldsymbol{\phi})+\eta_1R_1$ and $K_2(\boldsymbol{\psi}, \boldsymbol{\psi})+\eta_2R_2$, and store Cholesky factors for further uses. 
\end{remark} 
Following the arguments in \cite{chen2021solving}, the next theorem shows that  there exists a solution to \eqref{1ssmfgoptpk}.
\begin{theorem}
	\label{exi11}
	The minimization problem in \eqref{1ssmfgoptpk} admits a minimizer $(u^\dagger, m^\dagger, \overline{H}^\dagger)$ such that 
	\begin{align}
		\label{exp1df}
		\begin{cases}
			u^\dagger(x)=\langle K_1(x,\boldsymbol{\phi}), K_1(\boldsymbol{\phi}, \boldsymbol{\phi})^{-1}\boldsymbol{z}^\dagger\rangle,\\
			m^\dagger(x)=\langle  K_2(x, {\boldsymbol{\psi}}), K_2({\boldsymbol{\psi}}, {\boldsymbol{\psi}})^{-1}\boldsymbol{\rho}^\dagger\rangle,\\
			\overline{H}^\dagger = \lambda^\dagger,
		\end{cases}
	\end{align}
	where $(\boldsymbol{z}^\dagger, \boldsymbol{\rho}^\dagger, \lambda^\dagger)$ is a minimizer of  \eqref{1ssmfgoptfnpk}.   
\end{theorem}
\begin{proof}
	By the above arguments, \eqref{1ssmfgoptpk} is equivalent to \eqref{1ssmfgoptfnpk}. Hence, the key is to prove the existence of a minimizer to \eqref{1ssmfgoptfnpk}. 
	The argument here is similar to the proof of Theorem 1.1 in \cite{chen2021solving}, which
	proves the existence of a minimizer for the problem of the following form
	\begin{align}
		\label{chenform1}
		\begin{cases}
			\min\limits_{{\boldsymbol{z}}} \boldsymbol{z}^T\Theta^{-1}\boldsymbol{z}\\
			\text{s.t.}\ G(\boldsymbol{z})=0.
		\end{cases}
	\end{align}
	Here, $\Theta$ is assumed to be an invertible gram  matrix and $G$ is continuous. The arguments of \cite{chen2021solving} proceed by constructing a vector  $\boldsymbol{z}_*$ from the solution to the corresponding PDE such that $G(\boldsymbol{z}_*)=0$. Then, \eqref{chenform1} is equivalent to 
	\begin{align*}
		\begin{cases}
			\min\limits_{{\boldsymbol{z}}} \boldsymbol{z}^T\Theta^{-1}\boldsymbol{z}\\
			\text{s.t.}\ \mathcal{C}=\{\boldsymbol{z}|G(\boldsymbol{z})=0\}
			\cap \{\boldsymbol{z}| \boldsymbol{z}^T\Theta^{-1}\boldsymbol{z}\leq \boldsymbol{z}_*^T\Theta^{-1}\boldsymbol{z}_*\}.
		\end{cases}
	\end{align*}
	Since $\mathcal{C}$ is compact (by the continuity of $G$) and nonempty ($\mathcal{C}$ contains $\boldsymbol{z}_*$), the objective function $\boldsymbol{z}^T\Theta^{-1}\boldsymbol{z}$ admits a minimum in $\mathcal{C}$. 
	
	The above arguments can be extended  to our MFG setting. Let $(u^*, m^*, \overline{H}^*)$ be the solution to \eqref{1dsmfg}. We define the tuple $(\boldsymbol{z}_*, \boldsymbol{\rho}_*, \lambda_*)$ such that $\boldsymbol{z}_*$ is the vector with entries $z_{*,i}^{(1)}= u^*(x_i)$, $z_{*,i}^{(2)}=u^*_x(x_i)$, and $z_{*,i}^{(3)}=u^*_{xx}(x_i)$, $\boldsymbol{\rho}_*$ consists of elements $\rho^{(1)}_{*,i} = m^*(x_i)$ and $\rho^{(2)}_{*,i} = m^*_x(x_i)$, and $\lambda_*=\overline{H}^*$ for $i=1,\dots, M$. Then, $(\boldsymbol{z}_*, \boldsymbol{\rho}_*, \lambda_*)$ satisfies the constraints in \eqref{1ssmfgoptfnpk}. For a little abuse of notations, we denote by $G(\boldsymbol{z}, \boldsymbol{\rho},\lambda)=0$ the constraints in \eqref{1ssmfgoptfnpk} and define
	\begin{align*}
		\begin{split}
			\mathcal{C}_1=
			&\bigg\{(\boldsymbol{z}, \boldsymbol{\rho}, \lambda)| \boldsymbol{z}^TK_1(\boldsymbol{\phi}, \boldsymbol{\phi})^{-1}\boldsymbol{z} + \boldsymbol{\rho}^TK_2({\boldsymbol{\psi}}, {\boldsymbol{\psi}})^{-1}\boldsymbol{\rho} + |\lambda|^2 + \beta\bigg|\frac{1}{M}\sum\limits_{i=1}^{M}\rho_i^{(1)}- 1\bigg|^2+\beta\bigg|\frac{1}{M}\sum\limits_{i=1}^{M}z_i^{(1)}\bigg|^2\\
			&\quad\quad\quad\quad\leq \boldsymbol{z}_*^TK_1(\boldsymbol{\phi}, \boldsymbol{\phi})^{-1}\boldsymbol{z}_* + \boldsymbol{\rho}_*^TK_2({\boldsymbol{\psi}}, {\boldsymbol{\psi}})^{-1}\boldsymbol{\rho}_* + |\lambda_*|^2 + \beta\bigg|\frac{1}{M}\sum\limits_{i=1}^{M}\rho_{*,i}^{(1)}- 1\bigg|^2+\beta\bigg|\frac{1}{M}\sum\limits_{i=1}^{M}z_{*,i}^{(1)}\bigg|^2\bigg\}.
		\end{split}
	\end{align*} Then, \eqref{1ssmfgoptfnpk} is equivalent to 
	\begin{align}
		\label{mfg1dform2}
		\begin{cases}
			\min\limits_{\boldsymbol{z},\boldsymbol{\rho},\lambda} \boldsymbol{z}^TK_1(\boldsymbol{\phi}, \boldsymbol{\phi})^{-1}\boldsymbol{z} + \boldsymbol{\rho}^TK_2({\boldsymbol{\psi}}, {\boldsymbol{\psi}})^{-1}\boldsymbol{\rho} + |\lambda|^2\\
			\quad\quad\quad\quad + \beta\bigg|\frac{1}{M}\sum\limits_{i=1}^{M}\rho_i^{(1)}- 1\bigg|^2+\beta\bigg|\frac{1}{M}\sum\limits_{i=1}^{M}z_i^{(1)}\bigg|^2\\
			\text{s.t.}\  (\boldsymbol{z}, \boldsymbol{\rho}, \lambda)\in \mathcal{C}:=\bigg\{(\boldsymbol{z}, \boldsymbol{\rho},\lambda)|G(\boldsymbol{z}, \boldsymbol{\rho}, \lambda)=0\bigg\}\cap \mathcal{C}_1.
		\end{cases}
	\end{align}
	By the continuity of $G$ in $(\boldsymbol{z}, \boldsymbol{\rho}, \lambda)$ and the fact that $(\boldsymbol{z}_*, \boldsymbol{\rho}_*, \lambda_*)\in \mathcal{C}$, $\mathcal{C}$ is compact and non-empty. Hence, the objective function \eqref{mfg1dform2} achieves a minimum on $\mathcal{C}$. Thus, \eqref{1ssmfgoptfnpk} admits a minimizer. We conclude \eqref{exp1df} by \eqref{1drepre}. 
\end{proof}

Using a similar argument as in \cite{chen2021solving}, we obtain the following convergence theory. 
\begin{theorem}
	\label{cvtmfg1d}
	Assume that the kernels $K_1$ and $K_2$ are chosen such that $\mathcal{U}\subset\subset H^{s_1}(\mathbb{T})$ and $\mathcal{V}\subset\subset H^{s_2}(\mathbb{T})$ for some $s_1 > 3$ and $s_2>2$. Let $(u^*, m^*, \overline{H}^*)$ be the solution to \eqref{1dsmfg}. Denote by $(u^\dagger_{M,\beta}, m^\dagger_{M,\beta}, \overline{H}^\dagger_{M,\beta})$ a minimizer of \eqref{1ssmfgoptpk} with $M$ different sample points $\{x_i\}_{i=1}^M$ and the penalization parameter $\beta$. Suppose further that as $M\rightarrow\infty$, 
	\begin{align}
		\label{assXpoints}
		\sup_{x\in\mathbb{T}}\min_{1\leq i\leq M}|x-x_i|\rightarrow 0. 
	\end{align}
	Then, as $\beta$ and $M$ go to infinity, up to a subsequence, $(u^\dagger_{M,\beta}, m^\dagger_{M,\beta}, \overline{H}^\dagger_{M,\beta})$ converges to $(u^*, m^*, \overline{H}^*)$ pointwisely in $\mathbb{T}$ and in $H^{t_1}(\mathbb{T})\times H^{t_2}(\mathbb{T})\times\mathbb{R}$ for any $t_1\in (3, s_1)$ and any $t_2\in (2, s_2)$. 
\end{theorem}
\begin{proof}
	The argument is a direct adaptation of the proof of Theorem 1.2 in \cite{chen2021solving}. Let $(u^*, m^*, \overline{H}^*)$ be the solution of \eqref{1dsmfg}. Clearly, $(u^*, m^*, \overline{H}^*)$ satisfies the constraints in \eqref{1ssmfgoptpk}.
	According to \eqref{exp1du},  $u^*$ and $m^*$ are  smooth. Thus, a minimizer $(u^\dagger_{M,\beta}, m^\dagger_{M,\beta}, \overline{H}^\dagger_{M,\beta})$ to \eqref{1ssmfgoptpk} with $M$ different sample points and the penalization constant $\beta$ satisfies
	\begin{align}
		\label{bdummmhm}
		\begin{split}
			&\|u^\dagger_{M,\beta}\|_{\mathcal{U}}^2+\|m^\dagger_{M,\beta}\|_{\mathcal{V}}^2+|\overline{H}^\dagger_{M,\beta}|^2 + \beta\bigg|\frac{1}{M}\sum\limits_{i=1}^{M}m^\dagger_{M,\beta}(x_i)- 1\bigg|^2+\beta\bigg|\frac{1}{M}\sum\limits_{i=1}^{M}u^\dagger_{M,\beta}(x_i)\bigg|^2\\
			\leq& 	\|u^*\|_{\mathcal{U}}^2 + \|m^*\|_{\mathcal{V}}^2 +  |\overline{H}^*|^2+ \beta\bigg|\frac{1}{M}\sum\limits_{i=1}^{M}m^*(x_i)- 1\bigg|^2+\beta\bigg|\frac{1}{M}\sum\limits_{i=1}^{M}u^*(x_i)\bigg|^2.
		\end{split}
	\end{align}
	By $\int_{\mathbb{T}}u^*\dif x=0
	$, $\int_{\mathbb{T}}m^*\dif x = 1$ and \eqref{assXpoints}, there exist a sequence $\{(\beta_{p},M_p)\}_{p=1}^\infty$ and a constant $C$ such that
	\begin{align}
		\label{bdctgp}
		\beta_p\bigg|\frac{1}{M_p}\sum\limits_{i=1}^{M_p}m^*(x_i)- 1\bigg|^2+\beta_p\bigg|\frac{1}{M_p}\sum\limits_{i=1}^{M_p}u^*(x_i)\bigg|^2\leq C, \text{ for } p\geq 1. 
	\end{align} 
	Thus, using \eqref{bdummmhm} and \eqref{bdctgp}, we get 
	\begin{align}
		\label{bdumdag}
		\|u^\dagger_{M,\beta}\|_{\mathcal{U}}^2+\|m^\dagger_{M,\beta}\|_{\mathcal{V}}^2+|\overline{H}^\dagger_{M,\beta}|^2\leq C. 
	\end{align}
	Meanwhile, by $\mathcal{U}\subset\subset H^{s_1}$ and $\mathcal{V}\subset\subset H^{s_2}$, $s_1>3$ and $s_2>2$,  there exists a constant $C$ such that
	\begin{align}
		\label{bdummmhm2}
		\|u^\dagger_{M,\beta}\|_{H^{s_1}(\mathbb{T})}\leq 	C\|u^\dagger_{M,\beta}\|_{\mathcal{U}} \text{ and } \|m^\dagger_{M,\beta}\|_{H^{s_2}(\mathbb{T})}\leq 	C\|m^\dagger_{M,\beta}\|_{\mathcal{V}}.
	\end{align}
	For any $t_1\in (3, s_1)$ and any $t_2\in (2, s_2)$, we have $H^{s_1}(\mathbb{T})\subset\subset H^{t_1}(\mathbb{T})$ and $H^{s_2}(\mathbb{T})\subset\subset H^{t_2}(\mathbb{T})$. Thus, according to  \eqref{bdumdag} and \eqref{bdummmhm2},  there exits a limit $(u^\dagger_\infty, m^\dagger_\infty, \overline{H}^\dagger_\infty)\in H^{t_1}(\mathbb{T})\times H^{t_2}(\mathbb{T})\times\mathbb{R}$ such that, up to a subsequence, 
	\begin{align*}
		u^\dagger_{M_p,\beta_p} \rightarrow u^\dagger_\infty\  \text{in}\ H^{t_1}(\mathbb{T}),\\
		m^\dagger_{M_p,\beta_p} \rightarrow m^\dagger_\infty\  \text{in}\ H^{t_2}(\mathbb{T}), 
	\end{align*}
	and
	\begin{align*}
		\overline{H}_{M_p,\beta_p}^\dagger \rightarrow \overline{H}_\infty^\dagger\ \text{in}\ \mathbb{R},
	\end{align*}
	as $p\rightarrow\infty$.
	Since $H^{t_1}(\mathbb{T})\subset \subset C^2(\mathbb{T})$ and $H^{t_2}(\mathbb{T})\subset \subset C^1(\mathbb{T})$, 
	the limit $(u^\dagger_\infty, m^\dagger_\infty, \overline{H}^\dagger_\infty)$ satisfies the constraints in \eqref{1ssmfgoptpk}. By \eqref{assXpoints}, the collection of points $\{x_i\}_{i=1}^{M_p}$ forms a dense subset of $\mathbb{T}$ as $M_p\rightarrow\infty$. Thus, dividing both sides of \eqref{bdummmhm} by $\beta$ and passing $M$ and $\beta$ to infinity, we get $\int_{\mathbb{T}}m^\dagger_\infty\dif x=1$ and $\int_{\mathbb{T}}u^\dagger_\infty\dif x=0$. Hence, by the regularity of $u^\dagger_\infty$ and $m^\dagger_\infty$,  $(u^\dagger_\infty, m^\dagger_\infty, \overline{H}^\dagger_\infty)$ is a classical solution to \eqref{1ssmfgoptfnpk} and by the uniqueness of the solution to \eqref{1dsmfg}, $(u^\dagger_\infty, m^\dagger_\infty, \overline{H}^\dagger_\infty)=(u^*, m^*, \overline{H}^*)$. The limit $(u^\dagger_\infty, m^\dagger_\infty, \overline{H}^\dagger_\infty)$ is independent of the choice of the samples and the convergent subsequence. Therefore, we conclude that $(u^\dagger_{M,\beta}, m^\dagger_{M, \beta}, \overline{H}^\dagger_{M,\beta})$ converges to $(u^*, m^*, \overline{H}^*)$ pointwisely in $\mathbb{T}$ and in $H^{t_1}(\mathbb{T})\times H^{t_2}(\mathbb{T})\times\mathbb{R}$. 
\end{proof}

\subsection{The Fourier Features Method}
The main bottleneck of solving \eqref{1ssmfgoptfnpk} is to compute the inverses of $K_1(\boldsymbol{\phi}, \boldsymbol{\phi})$ and $K_2({\boldsymbol{\psi}},{\boldsymbol{\psi}})$, whose dimensions increase with the product of the size of samples and the number of linear operators in the MFG. Hence, the computational cost of the GP method grows dramatically as the number of samples increases. We propose the FF method based on the idea of approximating kernels with randomized trigonometric functions (see \cite{rahimi2007random, yu2016orthogonal, tompkins2020periodic}), which reduces the precomputation time and the amount of storage significantly.   
\begin{comment}
The random  Fourier Features method \cite{rahimi2007random} is based on 
Drawing $N/2$ samples $\{\omega_i\}_{i=1}^{N/2}$ from $p(\omega)$  and using the law of large numbers, we get
\begin{align}
\label{kernelap}
\kappa(x-y)=E_{\omega}[\cos(\omega^T(x-y))]\approx\frac{2}{N}\sum_{i=1}^{N/2}\bigg(\begin{bmatrix}
\cos(\omega_i^Tx) \\ \sin(\omega_i^Tx)
\end{bmatrix}^T\begin{bmatrix}
\cos(\omega_i^Ty) \\ \sin(\omega_i^Ty)
\end{bmatrix}\bigg).
\end{align}
Hence, $\kappa(x-y)$ can be approximated by the inner product of two vectors, i.e., 
\begin{align}
\label{kapp}
\kappa(x-y)\approx\zeta(x)^T\zeta(y),
\end{align}
where
\begin{align*}
\zeta(x)=\sqrt{\frac{2}{N}}[\sin(\omega_1^Tx),\dots,\sin(\omega_{N/2}^Tx), \cos(\omega_1^Tx),\dots,\cos(\omega_{N/2}^Tx)]. 
\end{align*}
\end{comment}

The idea of our Fourier Features method comes from the following observations. According to \eqref{exp1df}, the function $u^\dagger$ in the minimizer of \eqref{1ssmfgoptpk} has the following form
\begin{align}
	\label{udaggerform}
	u^\dagger(x) = \sum_{i=1}^{M}a_iK_1(x, x_i) + \sum_{i=1}^{M}b_i\partial_yK_1(x, x_i) + \sum_{i=1}^{M}c_i\partial_{y}^2K_1(x, x_i),
\end{align}
where $a_i, b_i, c_i\in \mathbb{R}$ are coefficients determined by $K_1(\boldsymbol{\phi}, \boldsymbol{\phi})^{-1}\boldsymbol{z}^\dagger$. 
We assume that $K_1$ is shift-invariant and properly scaled. According to the Bochner theorem \cite{rudin2017fourier, rahimi2007random}, the Fourier transform $p$ of $K_1$ is a probability distribution and
\begin{align}
	\label{Bochner}
	K_1(x-y)=E_{w}[\cos(\omega(x-y))],  \forall x, y\in \mathbb{R},
\end{align}
where $w$ is a random variable following $p$. 
Combining \eqref{udaggerform} and \eqref{Bochner}, we get
\begin{align*}
	\begin{split}
		u^\dagger(x)=&\sum_{i=1}^{M}a_iE_{w}[\cos(w(x-x_i))] + \sum_{i=1}^{M}b_iE_{w}[w\sin(w(x-x_i))] \\
		&-\sum_{i=1}^Mc_iE_w[w^2\cos(w(x-x_i))].
	\end{split}
\end{align*}
Thus, drawing $N/2$ samples $\{w_j\}_{j=1}^{N/2}$ from $p$, using trigonometric identities, and by the law of large numbers, we obtain
\begin{align}
	\label{kobsbigo}
	\begin{split}
		u^\dagger(x)\approx \sum_{j=1}^{N/2}\sum_{i=1}^{M}(a_i\cos(w_jx_i)-b_iw_j\sin(w_jx_i)-c_iw_j^2\cos(w_jx_i))\cos(w_jx)\\
		+\sum_{j=1}^{N/2}\sum_{i=1}^{M}(-a_i\sin(w_jx_i)+b_iw_j\cos(w_jx_i)+c_iw_j^2\sin(w_jx_i))\sin(w_jx). 
	\end{split}
\end{align}
\begin{comment}
\begin{align}
\label{kobsbigo}
\begin{split}
u^\dagger(x)\approx&\sum_{i=1}^{M}\sum_{j=1}^{N/2}a_i(\cos(\omega_j x_i)\cos(\omega_j x)  - \sin(\omega_j x_i)\sin(\omega_j x))\\
& + \sum_{i=1}^{M}\sum_{j=1}^{N/2}b_i\omega_j(\cos(\omega_jx_i)\sin(\omega_j x)-\sin(\omega_jx_i)\cos(\omega_jx))) \\
&-\sum_{i=1}^M\sum_{j=1}^{N/2}c_i\omega_j^2(\cos(\omega_jx_i)\cos(\omega_jx))-\sin(\omega_jx_i)\sin(\omega_jx))\\
=& \sum_{j=1}^{N/2}\sum_{i=1}^{M}(a_i\cos(\omega_jx_i)-b_i\omega_j\sin(\omega_jx_i)-c_i\omega_j^2\cos(\omega_jx_i))\cos(\omega_jx)\\
&+\sum_{j=1}^{N/2}\sum_{i=1}^{M}(-a_i\sin(\omega_jx_i)+b_i\omega_j\cos(\omega_jx_i)+c_i\omega_j^2\sin(\omega_jx_i))\sin(\omega_jx). 
\end{split}
\end{align}
\end{comment}
From \eqref{kobsbigo}, we observe that $u^\dagger$ can be approximated by a linear combination of randomized trigonometric functions. The same conclusion holds for $m^\dagger$ in \eqref{exp1df}. Hence, \eqref{kobsbigo} motivates us to approximate solutions of MFGs in the space generated by sampled trigonometric functions. 

Next, we present the FF method by solving \eqref{1dsmfg}. To deal with the periodic boundary condition in \eqref{1dsmfg}, inspired by \cite{tompkins2020periodic}, we approximate $u$ and $m$ by functions in the space generated by the Fourier series.  
Non-periodic settings are discussed in Section  \ref{secGF}, where we present a general framework of the FF method.

Given $N\in\mathbb{N}$, we define the set
\begin{align}
	\label{defG}
	\mathcal{G}^N=\bigg\{\phi\bigg|\phi(x)=c+\sum_{i=1}^{N}\alpha_i\sin(2\pi i x)+\sum_{i=1}^{N}\beta_i\cos(2\pi ix), c, \alpha_i, \beta_i\in \mathbb{R}\bigg\}. 
\end{align}
We call $\mathcal{G}^N$ the Fourier features space. Meanwhile, we consider functions $\sin(2\pi ix)$ and $\cos(2\pi i x)$ for $i=1,\dots, N$ as features in $\mathcal{G}^N$. 
Given $\gamma>0$, we define the functional 
\begin{align}
	\label{defca}
	J_\gamma(u, m, \overline{H})=\int_{\mathbb{T}}|u|^2\dif x + \int_{\mathbb{T}}|m|^2\dif x+\overline{H}^2+\gamma\mathcal{Q}(u, m, \overline{H}),
\end{align}
where 
\begin{align}
	\label{defmq}
	\begin{split}
		\mathcal{Q}(u, m, \overline{H}) =  
		&\int_{\mathbb{T}}|e^{H(x, u_x)-\overline{H}}-m|^2\dif x+\biggl |\int_{\mathbb{T}}u\dif x\biggl|^2+\biggl|\int_{\mathbb{T}}m\dif x - 1\biggl|^2\\
		&+\int_{\mathbb{T}}|m(u_{xx}+b_x(x))+m_x(u_x+b(x))|^2\dif x. 
	\end{split}
\end{align}
Then, we approximate the solution of \eqref{1dsmfg} by a minimizer of the following problem  
\begin{align}
	\label{1dcts}
	\min_{u^N\in \mathcal{G}^N, m^N\in \mathcal{G}^N, \overline{H}^N\in \mathbb{R}}J_\gamma(u^N, m^N, \overline{H}^N).
\end{align}
\begin{comment}
\begin{remark}
The following recovery problem
\begin{align} 
\label{1doptffrc}
\begin{cases}
\min\limits_{u^N\in \mathcal{G}^N, m^N\in\mathcal{G}^N, \overline{H}^N\in\mathbb{R}} \|u^N\|^2_{\mathcal{G}^N}+\|m^N\|^2_{\mathcal{G}^N}+|\overline{H}^N|^2,\\
\text{s.t.}\  H(x, u^N_x(x)) = \ln m^N+\overline{H}^N,\\
\quad\quad (m^N(u^N_x+b(x)))_x=0,\\
\quad\ \ \int_{\mathbb{T}}u^N\dif x = 0, \ \int_{\mathbb{T}}m^N\dif x=1,
\end{cases}
\end{align}
do not admit a minimizer since for a fixed $N\in\mathbb{T}$, $\mathcal{G}^N$ is not dense in the space of smooth functions. 
\end{remark}
\end{comment}
The following theorem gives the existence of a solution to \eqref{1dcts}.
\begin{theorem}
	\label{ffexi}
	The minimization problem \eqref{1dcts} admits a minimizer. 
\end{theorem}
\begin{proof}
	Let $N\in \mathbb{N}$. We define
	\begin{align}
		\label{defzeta1d}
		\boldsymbol{\zeta}(x)=[1, \sin(2\pi x),\dots,\sin(2\pi Nx), \cos(2\pi x),\dots,\cos(2\pi Nx)]^T.
	\end{align} 
	Then, for any $(u^N, m^N)\in \mathcal{G}^N\times\mathcal{G}^N$, there exits $\boldsymbol{\alpha}\in \mathbb{R}^{2N+1}$ and $\boldsymbol{\beta}\in \mathbb{R}^{2N+1}$ such that 
	\begin{align*}
		u^{N} = \boldsymbol{\alpha}^T\boldsymbol{\zeta} \text{ and } m^{N} = \boldsymbol{\beta}^T\boldsymbol{\zeta}. 
	\end{align*} 
	Therefore, \eqref{1dcts} is equivalent to 
	\begin{align}
		\label{1dctsfnt}
		\min_{\boldsymbol{\alpha}\in \mathbb{R}^{2N+1}, \boldsymbol{\beta}\in \mathbb{R}^{2N+1}, \overline{H}^N\in \mathbb{R}}J_{\gamma}(\boldsymbol{\alpha}^T\boldsymbol{\zeta}, \boldsymbol{\beta}^T\boldsymbol{\zeta}, \overline{H}^N).
	\end{align}
	We observe that the objective function in \eqref{1dctsfnt} is co-coercive and continuous. Thus, \eqref{1dctsfnt} admits a solution. Therefore, there is a minimizer to \eqref{1dcts}. 
\end{proof}
Then, we prove the convergence of the minimizer of \eqref{1dcts} to the solution of \eqref{1dsmfg} as $\gamma$ and $N$ go to infinity. First, we give a bound for the minimum of $J_\gamma$ in the following theorem. 
\begin{theorem}
	\label{convUpBdFF1}
	Let $(u^*, m^*, \overline{H}^*)$ be the solution to \eqref{1dsmfg}, let $\mathcal{G}^N$ be as in \eqref{defG},  and let $J_\gamma$ be as in \eqref{defca} for $\gamma>0$. Then,  for any sufficiently small $\epsi>0$, there exist a constant $C>0$, and functions $(u^N, m^N)\in \mathcal{G}^N\times\mathcal{G}^N$ such that 
	\begin{align}
		\label{bdJ}
		J_\gamma(u^N, m^N, \overline{H}^*) \leq 2\|u^*\|_{L^2(\mathbb{T})}^2+2\|m^*\|_{L^2(\mathbb{T})}^2+|\overline{H}^*|^2+C(1+\gamma)\epsilon^2.
	\end{align}
\end{theorem}
We refer readers to Appendix \ref{appendix} for the proof of Theorem \ref{convUpBdFF1}, which directly yields the following corollary.
\begin{corollary}
	\label{convFF1Q}
	Let $\mathcal{G}^N$ be as in \eqref{defG} and let $\mathcal{Q}$ be given in \eqref{defmq}. For any sufficient small $\epsilon>0$, there exist a constant $C>0$,  a sufficiently large $\gamma>0$, and $N>0$ such that any minimizer $({u}^{N,\gamma}, {m}^{N,\gamma}, {\overline{H}}^{N,\gamma})\in \mathcal{G}^N\times\mathcal{G}^N\times \mathbb{R}$ of \eqref{1dcts} corresponding to $\gamma$ satisfies
	\begin{align}
		\label{bdQ}
		\mathcal{Q}({u}^{N, \gamma}, {m}^{N, \gamma}, {\overline{H}}^{N, \gamma})\leq C\epsilon. 
	\end{align}
\end{corollary}
See Appendix \ref{appendix} for the proof of the above corollary. The following theorem shows that any minimizer of \eqref{1dcts} converges to the solution of \eqref{1dsmfg} as $N$  and $\gamma$ increase. 
\begin{theorem}
	\label{conv1dFF}
	There exists a sequence $\{(N_i,\gamma_i)\}_{i=1}^{\infty}$ such that any minimizer $(u^{i}, m^{i}, \overline{H}^{i})$ of \eqref{1dcts} corresponding to $(N_i, \gamma_i)$ satisfies $u^{i}\rightarrow u^*$ in $H^1(\mathbb{T})$, $m^{i}\rightarrow m^*$ in $L^1(\mathbb{T})$, and $\overline{H}^{i}\rightarrow \overline{H}^*$ in $\mathbb{R}$ as $i$ goes to infinity. 
\end{theorem}
The proof is given in Appendix \ref{appendix}. Next, we propose a numerical method to solve \eqref{1dctsfnt}, since \eqref{1dcts} and \eqref{1dctsfnt} are equivalent. Let $N\in \mathbb{N}$ and take $M$ samples $\{x_i\}_{i=1}^M$. Let  $\boldsymbol{\zeta}$ be given in \eqref{defzeta1d}. Introducing penalization parameters $\gamma>0$ and $\beta>0$, we reformulate \eqref{1dctsfnt} as an equivalent  two-level optimization problem
\begin{align}
	\label{1doptffrc2l}
	\begin{split}
		\min_{\boldsymbol{z}\in\mathbb{R}^{3M},\boldsymbol{\rho}\in\mathbb{R}^{2M},\lambda\in\mathbb{R}}\begin{cases}
			\min\limits_{\boldsymbol{\alpha}\in\mathbb{R}^{2N+1}, \boldsymbol{\beta}\in\mathbb{R}^{2N+1},\overline{H}\in\mathbb{R}} \|\boldsymbol{\alpha}\|^2+\|\boldsymbol{\beta}\|^2+|\overline{H}|^2\\
			\text{s.t.}\  \boldsymbol{\alpha}^T\boldsymbol{\zeta}(x_i) = z_i^{(1)}, \boldsymbol{\alpha}^T\boldsymbol{\zeta}_x(x_i) = z_i^{(2)}, \boldsymbol{\alpha}^T\boldsymbol{\zeta}_{xx}(x_i) = z_i^{(3)},\\ \ \ \ \ \  \boldsymbol{\beta}^T\boldsymbol{\zeta}(x_i) = \rho_i^{(1)}, \boldsymbol{\beta}^T\boldsymbol{\zeta}_x(x_i) = \rho_i^{(2)}, \overline{H} = \lambda, 
		\end{cases}\\
		+\beta\biggl | \frac{1}{M}\sum_{i=1}^{M}\rho_i^{(1)}-1\biggl|^2+\beta\biggl | \frac{1}{M}\sum_{i=1}^{M}z_i^{(1)}\biggl|^2+\gamma\sum_{i=1}^M|e^{H(x_i, z_i^{(2)})-\lambda}- \rho_i^{(1)}|^2\\
		+\gamma\sum_{i=1}^M|\rho_i^{(2)}(z_i^{(2)}+b(x_i))+\rho_i^{(1)}(z_i^{(3)}+b_x(x_i))|^2,\quad\quad\quad\quad\quad\quad\quad\ \ 
	\end{split}
\end{align}
where $\boldsymbol{z}$ and  $\boldsymbol{\rho}$ have forms in \eqref{zform} and \eqref{rhoform}. In \eqref{1doptffrc2l}, we use two penalization parameters to take into account different variability of the constraints.
Let $\boldsymbol{\phi}$, $\boldsymbol{\psi}$ be defined as in Subsection \ref{secgp1d}. Denote by $\mathcal{U}^*$ the dual of the Banach space $\mathcal{U}$ and by $[\cdot, \cdot]$ the duality pairing. Let $[\boldsymbol{\phi}, \boldsymbol{\zeta}]$ and $[\boldsymbol{\psi}, \boldsymbol{\zeta}]$ be the matrices with entries $[\phi_i, \zeta_j]$ and $[\psi_i, \zeta_j]$, separately.  Then, the first level optimization problem of \eqref{1doptffrc2l} yields
\begin{align}
	\label{zrholdrep}
	\boldsymbol{\alpha}=[\boldsymbol{\phi}, \boldsymbol{\zeta}]^T([\boldsymbol{\phi}, \boldsymbol{\zeta}][\boldsymbol{\phi}, \boldsymbol{\zeta}]^T)^{-1}\boldsymbol{z}, 
	\boldsymbol{\beta}=[\boldsymbol{\psi}, \boldsymbol{\zeta}]^T([\boldsymbol{\psi}, \boldsymbol{\zeta}][\boldsymbol{\psi}, \boldsymbol{\zeta}]^T)^{-1}\boldsymbol{\rho}, \overline{H}=\lambda. 
\end{align}
From \eqref{zrholdrep}, we get
\begin{align*}
	\|\boldsymbol{\alpha}\|^2=\boldsymbol{z}^T([\boldsymbol{\phi}, \boldsymbol{\zeta}][\boldsymbol{\phi}, \boldsymbol{\zeta}]^T)^{-1}\boldsymbol{z}, \|\boldsymbol{\beta}\|^2=\boldsymbol{\rho}^T([\boldsymbol{\psi}, \boldsymbol{\zeta}][\boldsymbol{\psi}, \boldsymbol{\zeta}]^T)^{-1}\boldsymbol{\rho}. 
\end{align*}
Hence, \eqref{1doptffrc2l} is equivalent to 
\begin{align}
	\label{1doptffrcfd}
	\begin{split}
		\min\limits_{\boldsymbol{z}\in\mathbb{R}^{3M},\boldsymbol{\rho}\in\mathbb{R}^{2M},\lambda\in\mathbb{R}}&\boldsymbol{z}^T([\boldsymbol{\phi}, \boldsymbol{\zeta}][\boldsymbol{\phi}, \boldsymbol{\zeta}]^T)^{-1}\boldsymbol{z}+ \boldsymbol{\rho}^T([\boldsymbol{\psi}, \boldsymbol{\zeta}][\boldsymbol{\psi}, \boldsymbol{\zeta}]^T)^{-1}\boldsymbol{\rho}+|\lambda|^2 \\
		&+\gamma\sum_{i=1}^M|e^{H(x_i, z_i^{(2)})-\lambda}- \rho_i^{(1)}|^2+\gamma\sum_{i=1}^M|\rho_i^{(2)}(z_i^{(2)}+b(x_i))+\rho_i^{(1)}(z_i^{(3)}+b_x(x_i))|^2,\\
		&+\beta\biggl | \frac{1}{M}\sum_{i=1}^{M}\rho_i^{(1)}-1\biggl|^2+\beta\biggl|\frac{1}{M}\sum_{i=1}^{M}z_i^{(1)}\biggl|^2.
	\end{split}
\end{align}
Then, we apply the Gauss--Newton method to solve \eqref{1doptffrcfd}. 
\begin{remark}	
	\label{regFF}
	In general, $[\boldsymbol{\phi}, \boldsymbol{\zeta}][\boldsymbol{\phi}, \boldsymbol{\zeta}]^T$ and $[\boldsymbol{\psi}, \boldsymbol{\zeta}][\boldsymbol{\psi}, \boldsymbol{\zeta}]^T$ are  ill-conditioned. Hence, we introduce two regularization parameters $\mu_1>0$ and $\mu_2>0$, and consider the inverses of  $[\boldsymbol{\phi}, \boldsymbol{\zeta}][\boldsymbol{\phi}, \boldsymbol{\zeta}]^T+\mu_1 I$ and $[\boldsymbol{\psi}, \boldsymbol{\zeta}][\boldsymbol{\psi}, \boldsymbol{\zeta}]^T+\mu_2 I$.
\end{remark}
\begin{remark}
	\label{AAqr}
	We propose an alternative method other than the Cholesky factorization to compute $[\boldsymbol{\phi}, \boldsymbol{\zeta}^N][\boldsymbol{\phi}, \boldsymbol{\zeta}^N]^T+\mu_1 I$ for $\mu_1>0$, which is the cornerstone for the FF method to save the precomputation time and the memory. For simplicity, denote $A:=[\boldsymbol{\phi}, \boldsymbol{\zeta}^N][\boldsymbol{\phi}, \boldsymbol{\zeta}^N]^T$. To solve $(AA^T+\mu_1 I)^{-1}$, we first perform QR decomposition on $A$. We get 
	\begin{align*}
		A = QR = \begin{bmatrix}
			Q_1 \ Q_2
		\end{bmatrix}\begin{bmatrix}
			R_1 \\ 0
		\end{bmatrix} = Q_1R_1, 
	\end{align*}
	where $Q\in \mathbb{R}^{3M\times 3M}$ is orthonormal, $R\in\mathbb{R}^{3M\times (2N+1)}$, $Q_1\in\mathbb{R}^{3M\times (2N+1)}$, $Q_2\in\mathbb{R}^{ 3M\times (3M-2N-1)}$, and $R_1\in\mathbb{R}^{(2N+1)\times (2N+1)}$. 
	Then, we have
	\begin{align*}
		AA^T+\mu_1 I = QRR^TQ^T + \mu_1 I = QRR^TQ^T + \mu_1 QQ^T.
	\end{align*}
	Thus, we get
	\begin{align*}
		(AA^T+\mu_1 I)^{-1} =& Q(RR^T+\mu_1 I)^{-1}Q^{T}=Q\begin{bmatrix}
			R_1R_1^T+\mu_1 I & 0 \\ 0 & \mu_1 I
		\end{bmatrix}^{-1}Q^T \\
		=& Q_1(R_1R_1^T+\mu_1 I)^{-1}Q_1^T + \frac{1}{\mu_1}Q_2Q_2^T. 
	\end{align*}
	Using 
	\begin{align*}
		I = QQ^T = Q_1Q_1^T + Q_2Q_2^T,
	\end{align*}
	we have 
	\begin{align}
		\label{AATinvf}
		(AA^T+\mu_1 I)^{-1} = Q_1(R_1R_1^T+\mu_1 I)^{-1}Q_1^T + \frac{1}{\mu_1}(I-Q_1Q_1^T). 
	\end{align}
	Since $R_1R_1^T+\mu_1 I \in \mathbb{R}^{(2N+1)\times (2N+1)}$, computing the right-hand side of \eqref{AATinvf} consumes less CPU time than calculating the left-hand side of \eqref{AATinvf} if $N<< M$. In the FF method, we precompute the QR decomposition of $A$ and the Cholesky decomposition of $R_1R_1^T+\mu_1I$, and store the results for further uses. We apply the same technique to compute the inverse of $[\boldsymbol{\psi}, \boldsymbol{\zeta}][\boldsymbol{\psi}, \boldsymbol{\zeta}]^T+\mu_2 I$.
\end{remark}
\begin{remark}
	The Woodbury matrix identity admits the same advantage as what we stated for \eqref{AATinvf} in Remark \ref{AAqr}. However, the Woodbury formula is numerically unstable in our experiments, especially for small values of $\mu_1$ and $\mu_2$.  
\end{remark}

\subsection{Numerical Results}
\begin{figure}[!hbtbp]
	\centering
	\begin{tabular}{c c c}            
		\begin{subfigure}[b]{0.3\textwidth}
			\includegraphics[width=\textwidth]{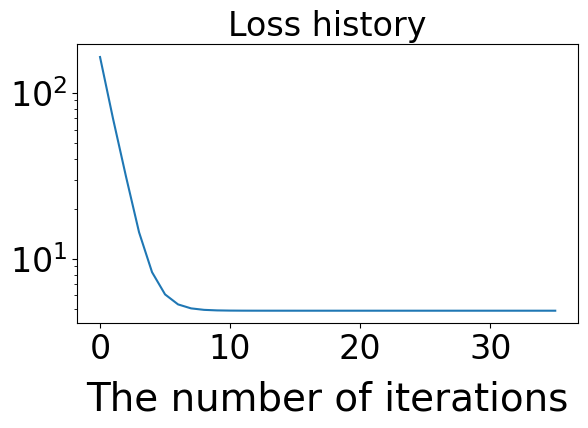}
			\caption{Loss history (GP).}
			\label{fig:GomesMFG:GPloss}
		\end{subfigure} &
		\begin{subfigure}[b]{0.3\textwidth}
			\includegraphics[width=\textwidth]{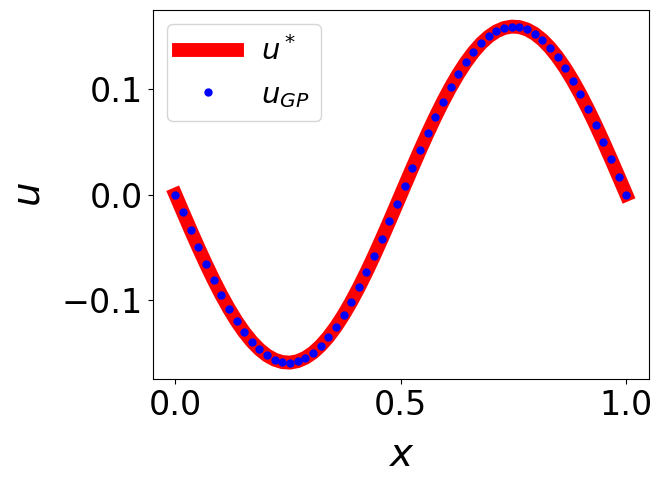}
			\caption{$u_{GP}$ v.s. $u^*$}
			\label{fig:GomesMFG:u}
		\end{subfigure} & 
		\begin{subfigure}[b]{0.3\textwidth}
			\includegraphics[width=\textwidth]{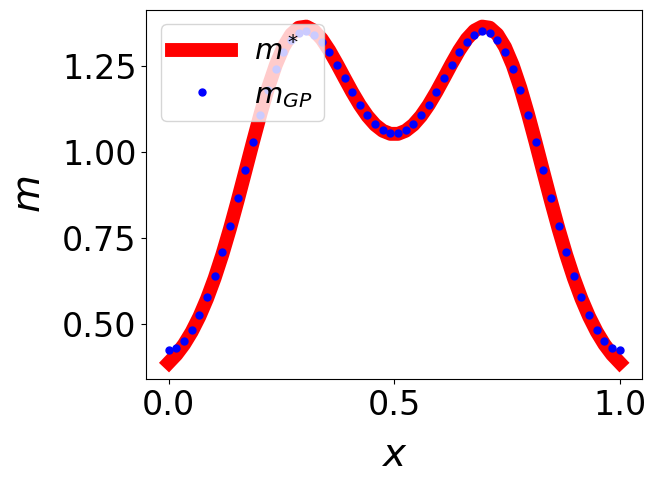}
			\caption{$m_{GP}$ v.s. $m^*$}
			\label{fig:GomesMFG:m}
		\end{subfigure} \\
		\begin{subfigure}[b]{0.3\textwidth}
			\includegraphics[width=\textwidth]{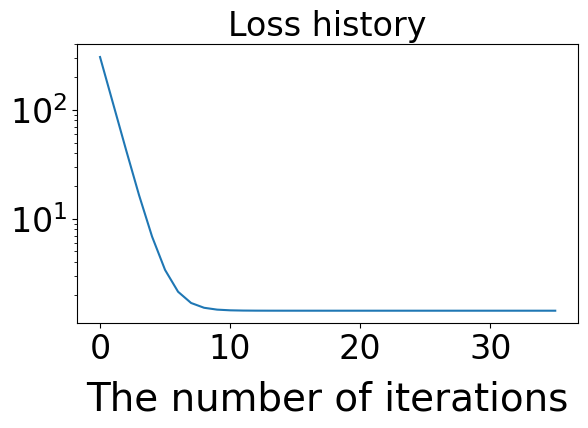}
			\caption{Loss history (FF).}
			\label{fig:GomesMFG:FFloss}
		\end{subfigure} &
		\begin{subfigure}[b]{0.3\textwidth}
			\includegraphics[width=\textwidth]{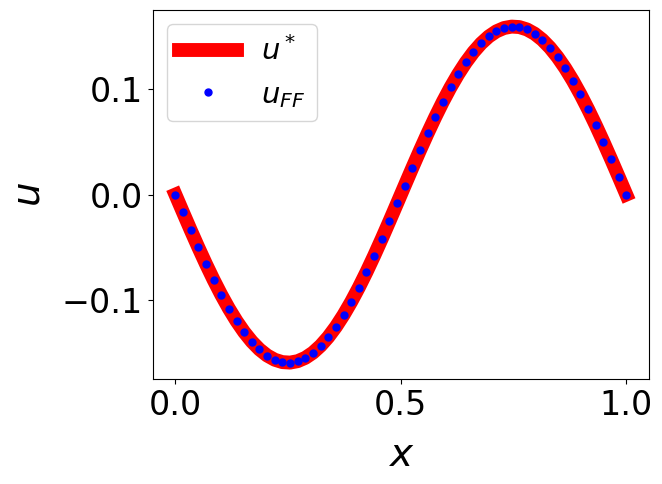}
			\caption{$u_{FF}$ v.s. $u^*$}
			\label{fig:GomesMFG:FFu}
		\end{subfigure} & 
		\begin{subfigure}[b]{0.3\textwidth}
			\includegraphics[width=\textwidth]{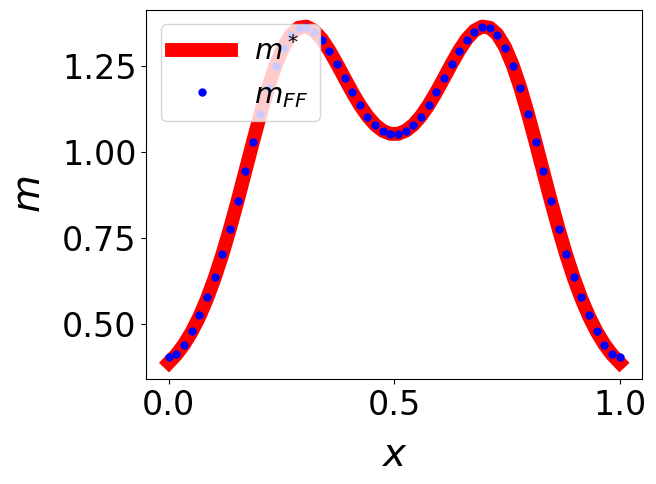}
			\caption{$m_{FF}$ v.s. $m^*$}
			\label{fig:GomesMFG:FFm}
		\end{subfigure} 
	\end{tabular}
	\caption{Numerical results for \eqref{1dsmfg}: (a), (d) the histories of loss functions; (b), (e) the numerical results $u_{GP}$ and $u_{FF}$ v.s. the explicit solution $u^*$; (c), (f)  the numerical approximations $m_{GP}$ and $m_{FF}$ v.s. the explicit solution $m^*$. }
	\label{fig:GomesMFG}
\end{figure}

To demonstrate the efficacy of our algorithms, we show here a simple numerical experiment by solving \eqref{1dsmfg}. We perform the calculation using MacBook Air 2015 (4GB RAM, Intel Core i5 CPU). Let $V(x)=\sin(\pi x)$ and $b(x)=\cos(2\pi x)$ for $x\in\mathbb{T}$. For the GP method, we choose the periodic kernel used in \cite{tompkins2020periodic} for both $K_1$ and $K_2$, i.e., 
\begin{align*}
	K_1(x, y)=K_2(x, y)=e^{\frac{\cos(2\pi (x-y))-1}{\sigma^2}}, \forall x, y \in \mathbb{T},
\end{align*}
with lengthscale  $\sigma>0$. 
We denote by $(u_{GP}, m_{GP}, \overline{H}_{GP})$ the numerical result of the GP method. For the FF algorithm, given $N\in \mathbb{N}$, we approximate the solution of \eqref{1dsmfg} in the space $\mathcal{G}^N$ given by \eqref{defG}. We represent the numerical solution of the FF method by $(u_{FF}, m_{FF}, \overline{H}_{FF})$. Let $(u^*, m^*, \overline{H}^*)$ be the solution of \eqref{1dsmfg} given by \eqref{exp1du}.  We choose $\sigma=0.6$ and $N=10$. For both methods, we use the Gauss--Newton method and take the step size $0.4$. We set the regularization parameters $\eta_1=\eta_2=\mu_1=\mu_2=10^{-6}$ in Remarks \ref{rmk1dnugget} and  \ref{regFF}. Meanwhile, we choose the penalization constants $\beta=10^{-6}$ and $\gamma=1$ in \eqref{relax1D} and \eqref{1doptffrc2l}. Both algorithms start from the same initial point and stop after {35}  iterations. In Figure \ref{fig:GomesMFG}, 
we show the numerical solutions of the GP method and the FF algorithm when we  take uniformly distributed $M=256$ sample points. We present $L^\infty$ errors in Table \ref{gomes1d:Errors}.  We see that the FF algorithm is comparable to the GP method in terms of accuracy. Table \ref{GomesMFG:CPUCP} records the CPU time consumed by the Cholesky and the QR decomposition by both methods, which implies that the FF algorithm costs less precomputation time than the GP method. We see that the CPU time consumed by the Cholesky decomposition in the FF method does not increase with the number of samples since we choose a fixed number of base functions. Surprisingly, the CPU time of the QR decomposition for the FF method even decreases when $M=2048$. We attribute it to the tall-and-skinny property of $[\boldsymbol{\phi}, \boldsymbol{\zeta}][\boldsymbol{\phi}, \boldsymbol{\zeta}]^T$ and $[\boldsymbol{\psi}, \boldsymbol{\zeta}][\boldsymbol{\psi}, \boldsymbol{\zeta}]^T$ when $M>>N$, i.e., $[\boldsymbol{\phi}, \boldsymbol{\zeta}][\boldsymbol{\phi}, \boldsymbol{\zeta}]^T$ and $[\boldsymbol{\psi}, \boldsymbol{\zeta}][\boldsymbol{\psi}, \boldsymbol{\zeta}]^T$ have more rows than columns.   
\begin{table}[ht!]
	\centering
	\begin{tabular}{c c c c}
		\hline
		& &  &  \\
		& \multicolumn{1}{c}{GP} & & \multicolumn{1}{c}{FF}\\
		\cline{2-2} \cline{4-4} 
		&  &  &  \\
		Errors of $u$    & $3.08\times 10^{-6}$ &  & $9.09\times 10^{-5}$ \\
		Errors of $m$ & $3.70\times 10^{-2}$ & & $1.68\times 10^{-2}$  \\
		Errors of $\overline{H}$ & $7.83\times 10^{-4}$ & & $9.87\times 10^{-4}$ \\
		\hline
	\end{tabular}
	\caption{ $L^\infty$ errors of numerical approximations for the solution of \eqref{1dsmfg}, given the number of sample points $M=256$. }
	\label{gomes1d:Errors}
\end{table}
\begin{table}[ht!]
	\centering
	\begin{tabular}{c c c c c}
		\hline
		& &  &  &      \\
		& \multicolumn{1}{c}{GP} & & \multicolumn{2}{c}{FF}  \\
		\cline{2-2} \cline{4-5}
		&  &  &   &  \\
		$M$    & Cholesky & & QR & Cholesky \\
		$256$ & $0.70$ & & $0.30$ & $0.18$ \\
		$512$ & $0.78$ & & $0.32$ & $0.18$\\
		$1024$ & $1.59$ &  & $0.34$  & $0.18$ \\
		$2048$ & $8.08$ & & $0.24$ & $0.18$ \\
		\hline
	\end{tabular}
	\caption{The precomputation time (in seconds) for solving \eqref{1dsmfg} as $M$ increases.}
	\label{GomesMFG:CPUCP}
\end{table}

\section{The General Frameworks}
\label{secGF}
This section presents the general frameworks of the GP method and the FF algorithm for solving MFGs.  We mainly state the settings for stationary MFGs. The arguments can be naturally adapted to the time-dependent cases. Numerical experiments on both stationary and time-dependent MFGs follow in Section \ref{secNumerical}. 
\subsection{The General Forms of MFGs}
Let $\Omega$ be a subset of $\mathbb{R}^d$. Suppose that the stationary MFGs of our interest have the form
\begin{align}
	\label{gsmfgs}
	\begin{cases}
		\mathcal{P}(u^*, m^*, \overline{H}^*)(x)= 0, \forall x \in \operatorname{int}\Omega,\\
		\mathcal{B}(u^*, m^*)(x) = 0, \forall x \in \partial \Omega,\\
		\int_{\Omega}u^*\dif x=0, \int_{\Omega}m^*\dif x = 1. 
	\end{cases}
\end{align}  
Here, $\mathcal{P}$ is a nonlinear differential operator and $\mathcal{B}$ represents a boundary operator.  We assume that \eqref{gsmfgs} admits a unique {classical}  solution $(u^*, m^*, \overline{H}^*)$. { If the solution of \eqref{gsmfgs} is not smooth enough, we suggest using the vanishing viscosity method \cite{cardaliaguet2010notes} or regularizing the MFG with smooth mollifiers  \cite{cesaroni2019stationary} to get a system with a solution of stronger regularity and applying our numerical methods to compute approximated solutions. }
\begin{remark}
In time-dependent settings, let $\Omega$ be a space-time domain. We consider MFGs taking the form
\begin{align}
	\label{gtdmfgs}
	\begin{cases}
		\mathcal{P}(u^*, m^*)(x) = 0, \forall x \in \operatorname{int} \Omega,\\
		\mathcal{B}(u^*, m^*)(x) = 0, \forall x \in \partial \Omega, 
	\end{cases}
\end{align} 
and assume that $(u^*, m^*)$ is a unique {classical} solution to \eqref{gtdmfgs}.  
\end{remark}
\subsection{The Gaussian Process Method}
\label{gpfw}
Using the method in \cite{chen2021solving}, we approximate $(u^*, m^*)$ in the solution of  \eqref{gsmfgs} by two GPs conditioned on PDEs at sampled collocation points in $\Omega$. Then, we compute the solution by calculating the MAP points of such conditioned GPs. More precisely, we take a set of samples $\{x_i\}_{i=1}^M$ in such a way that $x_1, \dots, x_{M_\Omega}\in \operatorname{int}\Omega$ and $x_{M_\Omega+1}, \dots, x_M\in \partial \Omega$ for $1\leq M_\Omega\leq M$. Let $\mathcal{U}$ and $\mathcal{V}$ be Banach Spaces with associated covariance operators $\mathcal{K}_{\mathcal{U}}: \mathcal{U}^*\mapsto \mathcal{U}$ and $\mathcal{K}_{\mathcal{V}}: \mathcal{V}^*\mapsto \mathcal{V}$, respectively. Following \cite{chen2021solving}, we introduce a penalization parameter $\beta>0$ and consider the following problem
\begin{align}
\label{OptGPProb}
\begin{cases}
\min\limits_{(u, m ,\overline{H})\in \mathcal{U}\times\mathcal{V}\times\mathbb{R}} \|u\|_{\mathcal{U}}^2 + \|m\|_{\mathcal{V}}^2 + |\overline{H}|^2\\
\quad\quad\quad\quad\quad\quad\quad+\beta\bigg|\frac{1}{M_\Omega}\sum\limits_{i=1}^{M_\Omega}u(x_i)\bigg|^2+\beta\bigg|\frac{1}{M_\Omega}\sum\limits_{i=1}^{M_\Omega}m(x_i)-1\bigg|^2\\
\text{s.t.}\quad \mathcal{P}(u, m, \overline{H})(x_i) = 0, \quad\text{for}\ i = 1, \dots M_{\Omega},\\
\quad\quad\ \mathcal{B}(u, m)(x_j) = 0, \quad \text{for}\ j = M_{\Omega+1},\dots, M. 
\end{cases}
\end{align}
We further make a similar assumption to Assumption 3.1 in \cite{chen2021solving}  on $\mathcal{P}$ and $\mathcal{B}$. 
\begin{hyp}
	\label{hypPB}
	For $1\leq Q_b\leq  Q$ and $1\leq D_b\leq D$, there exist bounded and linear operators $L_1, \dots, L_{Q_b}\in \mathcal{L}(\mathcal{U};C(\operatorname{\partial}\Omega))$, $L_{Q_b+1}, \dots, L_{Q}\in \mathcal{L}(\mathcal{U};C(\operatorname{int}\Omega))$, $J_1, \dots, J_{D_b}\in \mathcal{L}(\mathcal{V};C(\operatorname{\partial}\Omega))$, $J_{D_b+1}, \dots, J_{D}\in \mathcal{L}(\mathcal{V};C(\operatorname{int}\Omega))$, and  continuous nonlinear maps $P$ and ${B}$ such that 
	\begin{align}
		\label{gsmfgslo}
		\begin{cases}
			\mathcal{P}(u^*, m^*, \overline{H}^*)(x) = \\
			\quad P(L_{Q_{b+1}}(u^*)(x), \dots, L_{Q}(u^*)(x), J_{D_{b+1}}(m^*)(x), \dots, J_D(m^*)(x), \overline{H}^*), \forall x\in \operatorname{int}\Omega,\\
			\mathcal{B}(u^*, m^*)(x) = \\
			\quad B(L_1(u^*)(x), \dots, L_{Q_b}(u^*)(x), J_1(m^*)(x), \dots, J_{D_b}(m^*)(x)), \forall x\in \operatorname{\partial}\Omega.
		\end{cases}
	\end{align}
\end{hyp}
Following \cite{chen2021solving}, under Assumption \ref{hypPB}, we define functionals $\phi_i^{(q)}\in \mathcal{U}^*$ and $\psi_i^{(p)}\in \mathcal{V}^*$ as
\begin{align*}
\phi_{i}^{(q)} = \delta_{x_i}\circ L_q, \text{where}\ \begin{cases}
	M_\Omega+1\leq i\leq M, &\text{if}\ 1\leq q \leq Q_b, \\
	1\leq i\leq M_{\Omega}, &\text{if}\ Q_{b+1}\leq q\leq Q,
\end{cases}
\end{align*} 
and 
\begin{align*}
	\psi_{i}^{(p)} = \delta_{x_i}\circ J_p, \text{where}\ \begin{cases}
		M_\Omega+1\leq i\leq M, &\text{if}\ 1\leq p \leq D_b, \\
		1\leq i\leq M_{\Omega}, &\text{if}\ D_{b+1}\leq p\leq D,
	\end{cases}
\end{align*} 
For ease of presentation, we denote by $\boldsymbol{\phi}^{(q)}$ the vector consisting of $\phi^{(q)}_i$ and define
\begin{align}
\label{vecPhi}
\boldsymbol{\phi} = (\boldsymbol{\phi}^{(1)}, \dots, \boldsymbol{\phi}^{(Q)})\in (\mathcal{U}^*)^{\bigotimes N_{\mathcal{U}}}, \text{where}\ N_\mathcal{U} = (M-M_\Omega)Q_b+M_\Omega(Q-Q_b). 
\end{align}
Similarly, we concatenate $\psi_i^{(p)}$ to get the vector $\boldsymbol{\psi}^{(p)}$ and denote
\begin{align}
\label{vecPsi}
\boldsymbol{\psi} = (\boldsymbol{\psi}^{(1)},\dots, \boldsymbol{\psi}^{(D)})\in (\mathcal{V}^*)^{\bigotimes N_{\mathcal{V}}}, \text{where}\ N_\mathcal{V} = (M-M_\Omega)D_b+M_\Omega(D-D_b). 
\end{align}
According to Assumption \ref{hypPB}, we define the nonlinear map $G$ such that for any $u\in \mathcal{U}$, $m\in \mathcal{V}$, $\overline{H}\in \mathbb{R}$, 
\begin{align}
\label{nonlF}
(G([\boldsymbol{\phi}, u], [\boldsymbol{\psi}, u], \overline{H}))_i=\begin{cases}
P([\phi_i^{(Q_{b+1})}, u], \dots, [\phi^{(Q)}_i, u], [\psi_i^{({D_{b+1}})}, m], \dots, [\psi^{(D)}_i, m], \overline{H})\\
\quad\quad\quad\quad\quad\quad\quad \text{if}\ i \in \{1, \dots, M_{\Omega}\}, \\
B([\phi_i^{(1)}, u], \dots, [\phi^{(Q_b)}_i, u], [\psi_i^{(1)}, m], \dots, [\psi^{(D_b)}_i, m], \overline{H})\\
\quad\quad\quad\quad\quad\quad\quad \text{if}\ i \in \{M_{\Omega}+1, \dots, M\}. \\
\end{cases}
\end{align}
Hence, we can rewire \eqref{OptGPProb} as
\begin{align}
\label{cpctOptGGPProb}
\begin{cases}
	\min\limits_{(u, m ,\overline{H})\in \mathcal{U}\times\mathcal{V}\times\mathbb{R}} \|u\|_{\mathcal{U}}^2 + \|m\|_{\mathcal{V}}^2 + |\overline{H}|^2\\
	\quad\quad\quad\quad\quad\quad\quad+\beta\bigg|\frac{1}{M_\Omega}\sum\limits_{i=1}^{M_\Omega}u(x_i)\bigg|^2+\beta\bigg|\frac{1}{M_\Omega}\sum\limits_{i=1}^{M_\Omega}m(x_i)-1\bigg|^2\\
	\text{s.t. } G([\boldsymbol{\phi}, u], [\boldsymbol{\psi}, m], \overline{H}) = 0. 
\end{cases}
\end{align}
\begin{remark}
Similarly, in time-dependent settings, we consider
\begin{align}
	\label{cpctOptGGPProbtd}
	\begin{cases}
		\min\limits_{(u, m)\in \mathcal{U}\times\mathcal{V}} \|u\|_{\mathcal{U}}^2 + \|m\|_{\mathcal{V}}^2,\\
		\text{s.t. } G([\boldsymbol{\phi}, u], [\boldsymbol{\psi}, m]) = 0, 
	\end{cases}
\end{align}
where $G$ is an analog to the one in \eqref{nonlF}
\end{remark}
The following theorem gives the foundation for the GP method to solve \eqref{cpctOptGGPProb}. 
\begin{theorem}
\label{GPFund}
Suppose that Assumption \ref{hypPB} holds. Let $N_{\mathcal{U}}$, $\boldsymbol{\phi}$,  $N_{\mathcal{V}}$, $\boldsymbol{\psi}$, and $G$ be as in \eqref{vecPhi}, \eqref{vecPsi}, and \eqref{nonlF}. Define matrices $\Theta\in \mathbb{R}^{N_{\mathcal{U}}\times N_{\mathcal{U}}}$, $\Psi\in \mathbb{R}^{N_{\mathcal{V}}\times N_{\mathcal{V}}}$ such that
\begin{align*}
\Theta_{ij} = [\phi_i, \mathcal{K}_{\mathcal{U}}\phi_j], 1\leq i, j \leq N_{\mathcal{U}}, \text{ and }
\Psi_{ks}=[\psi_k, \mathcal{K}_{\mathcal{V}}\psi_s], 1\leq k, s\leq N_{\mathcal{V}}. 
\end{align*}
Assume further that $\Theta$ and $\Psi$ are invertible. Let  $\boldsymbol{\chi}=(\chi_1,\dots,\chi_{N_{\mathcal{U}}})$ and $\boldsymbol{\eta}=(\eta_1,\dots,\eta_{N_{\mathcal{V}}})$ be the vectors with elements
\begin{align*}
	\chi_i = \sum_{n=1}^{N_\mathcal{U}}\Theta_{in}^{-1}\mathcal{K}_{\mathcal{U}}\phi_n, i=1,\dots, N_{\mathcal{U}}, \text{ and }\eta_i = \sum_{n=1}^{N_\mathcal{V}}\Psi_{in}^{-1}\mathcal{K}_{\mathcal{V}}\psi_n, i=1,\dots,N_{\mathcal{V}},  
\end{align*}
where $\Theta_{in}^{-1}$ and $\Psi_{in}^{-1}$ are the elements of $\Theta^{-1}$ and $\Psi^{-1}$ at the $i$th row and the $n$th column. 
Then, $(u^\dagger, m^\dagger, \overline{H}^\dagger)$ is a solution to \eqref{cpctOptGGPProb} if and only if 
\begin{align}
\label{udamdaghdarf}
u^\dagger = \boldsymbol{\chi}^T \boldsymbol{z}^\dagger, m^\dagger = \boldsymbol{\eta}^T \boldsymbol{\rho},\ \text{and}\  \overline{H}^\dagger = \lambda,
\end{align}
where $(\boldsymbol{z}^\dagger$, $\boldsymbol{\rho}^\dagger, \lambda)$ is a minimizer to 
\begin{align}
\label{gpfin}
\begin{cases}
\min\limits_{\boldsymbol{z}\in \mathbb{R}^{N_{\mathcal{U}}}, \boldsymbol{\rho}\in \mathbb{R}^{N_{\mathcal{V}}}, \lambda\in \mathbb{R}}\boldsymbol{z}^T\Theta^{-1}\boldsymbol{z} + \boldsymbol{\rho}^T\Psi^{-1}\boldsymbol{\rho}+|\lambda|^2\\
\quad\quad\quad\quad\quad\quad\quad+\beta\bigg|\frac{1}{M_\Omega}\sum\limits_{i=1}^{M_\Omega}z_i^{(1)}\bigg|^2+\beta\bigg|\frac{1}{M_\Omega}\sum\limits_{i=1}^{M_\Omega}\rho_i^{(1)}-1\bigg|^2\\
\text{s.t. } G(\boldsymbol{z}, \boldsymbol{\rho}, \lambda) = 0. 
\end{cases}
\end{align}
\end{theorem}
\begin{proof}
We conclude using similar arguments as in the proof of Theorem \ref{exi11}. Let $(u^*, m^*, \overline{H}^*)$ be the solution to \eqref{gsmfgs} and define $\boldsymbol{z}_* = [\boldsymbol{\phi}, u^*], \boldsymbol{\rho}_* = [\boldsymbol{\psi}, m^*]$ and $\lambda_* = \overline{H}^*$. Then, $G(\boldsymbol{z}_*, \boldsymbol{\rho}_*, \lambda_*)=0$. Thus, the minimization problem \eqref{gpfin} can be restricted to the form of \eqref{mfg1dform2}. Hence, \eqref{gpfin} admits a minimizer. Following nearly identical steps to the derivation of  \eqref{1ssmfgoptfnpk}, we conclude \eqref{udamdaghdarf}. 
\end{proof} 
Next, following \cite{chen2021solving}, we have the convergence theorem.
\begin{theorem}
	\label{cvtmfgnd}
	Assume that Assumption \ref{hypPB} holds and that the MFG in  \eqref{gsmfgs} has a unique {classical}  solution $(u^*, m^*, \overline{H}^*)$ in the space $\mathcal{U}\times \mathcal{V}\times\mathbb{R}$. Assume further that $\mathcal{U}\subset \subset \mathcal{H}_1\subset\subset C^{t_1}(\Omega)\cap C^{t_1'}(\partial\Omega)$ and that $\mathcal{V}\subset \subset \mathcal{H}_2\subset\subset C^{t_2}(\Omega)\cap C^{t_2'}(\partial\Omega)$, where $\mathcal{H}_1$ and $\mathcal{H}_2$ are Banach spaces and $t_1, t_1', t_2$, and $t_2'$ are sufficiently large. 
	Denote by $\{x_i\}_{i=1}^M$ the collection of samples with $M$ points.  Suppose further that as $M\rightarrow\infty$, 
	\begin{align}
		\label{assnXpoints}
		\sup_{x\in \operatorname{int}\Omega}\min_{1\leq i\leq M_\Omega}|x-x_i|\rightarrow 0\ \text{and}\  \sup_{x\in \partial\Omega}\min_{M_{\Omega}+1\leq i\leq M}|x-x_i|\rightarrow 0.
	\end{align}
Given $M$ and $\beta$, let $(u^\dagger_{M,\beta}, m^\dagger_{M,\beta}, \overline{H}^\dagger_{M,\beta})$ be a minimizer of \eqref{cpctOptGGPProb}.
	Then, as $M$ and $\beta$ go to infinity, up to a sub-sequence,  $(u^\dagger_{M,\beta}, m^\dagger_{M,\beta}, \overline{H}^\dagger_{M,\beta})$ converges to $(u^*, m^*, \overline{H}^*)$ pointwisely in $\Omega$ and in $\mathcal{H}_1\times\mathcal{H}_2\times\mathbb{R}$. 
\end{theorem}
\begin{proof}
The argument is similar to the proof of Theorem \ref{cvtmfg1d}. Given $M$ sample points and the penalization parameter $\beta$, denote by $(u^\dagger_{M,\beta}, m^\dagger_{M,\beta}, \overline{H}^\dagger_{M,\beta})$ a minimizer to \eqref{cpctOptGGPProb}. Using the fact that the {classical}  solution $(u^*, m^*, \overline{H}^*)$ to \eqref{gsmfgs} satisfies the constraints in \eqref{cpctOptGGPProb}, we prove that there exists a sequence $\{(M_p, \beta_p)\}_{p=1}^\infty$ such that  $\|u^\dagger_{M_p,\beta_p}\|_{\mathcal{U}}$, $\|m^\dagger_{M_p,\beta_p}\|_{\mathcal{V}}$, and $|\overline{H}_{M_p,\beta_p}|$ are uniformly bounded for all $p\geq 1$, and that $\int_{\Omega}m^\dagger_{M_p,\beta_p}\dif x\rightarrow 1$ and  $\int_{\Omega}u^\dagger_{M_p,\beta_p}\dif x\rightarrow0$. Since $\mathcal{U}\subset \subset \mathcal{H}_1$, $\mathcal{V}\subset \subset \mathcal{H}_2$, and a closed bounded set in $\mathbb{R}$ is compact, the sequence $(u^\dagger_{M,p}, m^\dagger_{M,p}, \overline{H}_{M,p}^\dagger)$ converge, up to a sub-sequence,  to a limit $(u_\infty^\dagger, m_\infty^\dagger, \overline{H}_\infty^\dagger)$ in $\mathcal{H}_1\times\mathcal{H}_2\times\mathbb{R}$ as $p\rightarrow\infty$. Using  $\mathcal{H}_1\subset\subset C^{t_1}(\Omega)\cap C^{t_1'}(\partial\Omega)$ and that $\mathcal{H}_2\subset\subset C^{t_2}(\Omega)\cap C^{t_2'}(\partial\Omega)$, we conclude that $(u_\infty^\dagger, m_\infty^\dagger, \overline{H}_\infty^\dagger)$ satisfies the constraints in \eqref{OptGPProb} at all points in $\{x_i\}_{i=1}^M$. Due to \eqref{assnXpoints}, $\{x_i\}_{i=1}^M$ is dense in $\Omega$ as $M\rightarrow\infty$. Hence,  $(u_\infty^\dagger, m_\infty^\dagger, \overline{H}_\infty^\dagger)$ solves \eqref{gsmfgs}. Since the solution to \eqref{gsmfgs} is unique, we conclude that $(u_\infty^\dagger, m_\infty^\dagger, \overline{H}_\infty^\dagger)=(u^*, m^*, \overline{H}^*)$. 
\end{proof}
\begin{remark}
To deal with the nonlinear constraints of \eqref{gpfin}, we  introduce two penalization parameters, $\gamma>0$ and $\beta>0$, and consider the relaxation formulation similar to \eqref{relax1D}. 
\end{remark}
\begin{remark}
Under similar Assumptions of Theorems \ref{GPFund} and \ref{cvtmfgnd},  \eqref{cpctOptGGPProbtd} admits a solution $(u^\dagger_M, m^\dagger_M)$ for given sample size $M$, and $(u^\dagger_M, m^\dagger_M)$ converges to the {classical} solution of \eqref{gtdmfgs} as $M$ goes to infinity.  
\end{remark}
\subsection{The Fourier Features Method}
\label{fffw}
Let $\mathcal{O}=\{\zeta(\cdot,\omega)|\omega\in \mathcal{W}\}$ be a family of base functions parametrized over the set $\mathcal{W}$, where  
\begin{align}
\label{ctszeta}
\sup_{x\in\Omega}|\zeta(x, \omega)|\leq 1, \forall w\in \mathcal{W} \text{ and } \int_{\Omega}\zeta^2(x,\omega)\dif x <\infty. 
\end{align}
We propose to approximate $u^*$ and $m^*$ in the solution of \eqref{gsmfgs} by linear combinations of functions sampled from  $\mathcal{O}$.  More precisely, given $N\in\mathbb{N}$, we take samples $\{\omega_i\}_{i=1}^N$ and $\{\tau_i\}_{i=1}^N$ from $\mathcal{W}$ and define vector valued functions $\boldsymbol{\zeta}^N$ and $\boldsymbol{\vartheta}^N$ by
\begin{align}
	\label{defzeta}
	\boldsymbol{\zeta}^N(\cdot)=[\zeta(\cdot;w_1),\dots, \zeta(\cdot;w_{N})]^T \text{ and }\boldsymbol{\vartheta}^N(\cdot)=[\zeta(\cdot;\tau_1),\dots, \zeta(\cdot;\tau_{N_\mathcal{T}})]^T.
\end{align} 
Then, we define spaces
\begin{align}
\label{defSpG}
\mathcal{G}_{\boldsymbol{\zeta}^N}=\{\boldsymbol{\alpha}^T\boldsymbol{\zeta}^N|\boldsymbol{\alpha}\in \mathbb{R}^{N}\}\ \text{and}\ \mathcal{G}_{\boldsymbol{\vartheta}^N}=\{\boldsymbol{\beta}^T\boldsymbol{\vartheta}^N|\boldsymbol{\beta}\in \mathbb{R}^{N}\}.
\end{align}
Meanwhile, for $\zeta$ satisfying \eqref{ctszeta}, we equip the spaces $\mathcal{G}_{\boldsymbol{\zeta}^N}$ and $\mathcal{G}_{\boldsymbol{\vartheta}^N}$ with the norms
\begin{align*}
\bigl\|u_{\boldsymbol{\alpha}}\bigl\|_{\mathcal{G}_{\boldsymbol{\zeta}^N}}=|\boldsymbol{\alpha}|, \forall u_{\boldsymbol{\alpha}}=\boldsymbol{\alpha}^T\boldsymbol{\zeta}^N\in \mathcal{G}_{\boldsymbol{\zeta}^N} \text{ and }\bigl\|m_{\boldsymbol{\beta}}\bigl\|_{\mathcal{G}_{\boldsymbol{\vartheta}^N}}=|\boldsymbol{\beta}|, \forall m_{\boldsymbol{\beta}}=\boldsymbol{\beta}^T\boldsymbol{\vartheta}^N\in \mathcal{G}_{\boldsymbol{\vartheta}^N}.
\end{align*}
By \eqref{ctszeta}, the norms $\|\cdot\|_{\mathcal{G}_{\boldsymbol{\zeta}^N}}$ and $\|\cdot\|_{L^2(\Omega)}$ are equivalent for the space $\mathcal{G}_{\boldsymbol{\zeta}^N}$. Then, we approximate the solution of \eqref{gsmfgs} by the minimizer of the following problem
\begin{align}
	\label{OptFFProb}
	\begin{split}
		\min_{u_{\boldsymbol{\alpha}}\in\mathcal{G}_{\boldsymbol{\zeta}^N}, m_{\boldsymbol{\beta}}\in\mathcal{G}_{\boldsymbol{\vartheta}^N}, \overline{H}\in\mathbb{R}} \|u_{\boldsymbol{\alpha}}\|^2_{L^2(\Omega)}+\|m_{\boldsymbol{\beta}}\|^2_{L^2(\Omega)} +\gamma\|\mathcal{P}(u_{\boldsymbol{\alpha}}, m_{\boldsymbol{\beta}}, \overline{H})\|_{L^2(\Omega)}^2\\
		\quad\quad\quad\quad+\gamma\|\mathcal{B}(u_{\boldsymbol{\alpha}}, m_{\boldsymbol{\beta}})\|_{L^2(\Omega)}^2+\gamma\bigg|\int_{\Omega}m_{\boldsymbol{\beta}}\dif x-1\bigg|^2+\gamma\bigg|\int_{\Omega}u_{\boldsymbol{\alpha}}\dif x\bigg|^2,
	\end{split}
\end{align}
where $\mathcal{P}$ and $\mathcal{B}$ are given in Assumption \ref{hypPB}, and $\gamma>0$.
\begin{remark}
	\label{chbasis}
	When the domain $\Omega=\mathbb{T}^d$, we choose $\mathcal{O}$ such that
	\begin{align*}
		\mathcal{O}=\{\zeta(\cdot, \omega)=\cos(2\pi\langle a, \cdot\rangle + b)|\omega=(a,b)\in \mathbb{Z}^d\times\{0,\pi/2\}\}. 
	\end{align*}
Then, given $N\in \mathbb{N}$, we select a subset $\{\omega_i=(a_i,b_i)\}_{i=1}^N$ in $\mathbb{Z}^d\times\{0,\pi/2\}$. Thus, the vector valued function $\boldsymbol{\zeta}^N$ in \eqref{defzeta} becomes
\begin{align*}
	\boldsymbol{\zeta}^N(\cdot)=[\cos(2\pi\langle a_1,\cdot\rangle + b_1), \dots, \cos(2\pi\langle a_N,\cdot\rangle+b_N)]^T. 
\end{align*}  
Hence, for $\boldsymbol{\alpha}=\{\alpha_{1},\dots, \alpha_{N}\}\in\mathbb{R}^{N}$, the function $u_{\boldsymbol{\alpha}}\in \mathcal{G}_{\boldsymbol{\zeta}^N}$ in \eqref{defSpG} has the following form
\begin{align*}
u_{\boldsymbol{\alpha}}(x)=\sum_{i=1}^{N}\alpha_{i}\cos(2\pi a_i^T x+b_i). 
\end{align*}
The same construction holds for the space $\mathcal{G}_{\boldsymbol{\vartheta}^N}$. 
\end{remark}
\begin{remark}
\label{rmkrff}
	When the domain $\Omega$ is non-periodic, we choose 
	\begin{align*}
		\mathcal{O}=\{\zeta(\cdot,\omega)=\cos(2\pi\langle\omega, \cdot\rangle)|\omega\in \mathbb{R}^d\}\cup\{\zeta(\cdot,\omega)=\sin(2\pi\langle\omega, \cdot\rangle)|\omega\in \mathbb{R}^d\}.
	\end{align*}
Then, we take $N/2$, $N\in\mathbb{N}$ and $N$ is even, samples $\{\omega_i\}_{i=1}^{N/2}$ from $\mathbb{R}^d$ in such a way that the nonlinear map $\boldsymbol{\zeta}^N$ in \eqref{defzeta} is defined as
	\begin{align*}
	\boldsymbol{\zeta}^N(x)=\sqrt{\frac{2}{N}}[\sin(\omega_1^Tx), \dots, \sin(\omega_{N/2}^Tx), \cos(\omega_1^Tx),\dots, \cos(\omega_{N/2}^Tx)]^T, \forall x\in \mathbb{R}^d. 
	\end{align*}
To sample $\{\omega_i\}_{i=1}^{N/2}$, we use the method of orthogonal random features \cite{yu2016orthogonal}. More precisely, the matrix $W=[\omega_1,\dots,\omega_{N/2}]^T$ satisfies
\begin{align}
\label{mtxW}
W = \frac{1}{\varsigma}SQ,
\end{align}
where $\varsigma>0$, $Q$ is a uniformly distributed random orthogonal matrix, and $S$ is a diagonal matrix with entries sampled i.i.d from $\chi$-distribution with $d$ degrees of freedom. We refer readers to \cite{yu2016orthogonal} for more details about the construction of $W$. 
 
Then, given 
	$\boldsymbol{\alpha}=\{\alpha_{1},\dots, \alpha_{N}\}\in\mathbb{R}^{N}$ and $u_{\boldsymbol{\alpha}}\in \mathcal{G}_{\boldsymbol{\zeta}^N}$ in \eqref{defSpG}, we obtain
\begin{align*}
	u_{\boldsymbol{\alpha}}(x)=\sum_{i=1}^{N/2}\alpha_{i}\sin(\omega_i^T x) + \sum_{i=N/2+1}^{N}\alpha_{i}\cos(\omega_i^T x).  
\end{align*}
We use the same method to build the space $\mathcal{G}_{\boldsymbol{\vartheta}^N}$. 
\end{remark}
\begin{remark}
In a time-dependent setting, the space-time domain is non-periodic. Hence, we use features given in Remark \ref{rmkrff} and approximate the solution to \eqref{gtdmfgs} by the minimizer of the following problem
 \begin{align*}
 	\begin{split}
 		\min_{u_{\boldsymbol{\alpha}}\in\mathcal{G}_{\boldsymbol{\zeta}^N}, m_{\boldsymbol{\beta}}\in\mathcal{G}_{\boldsymbol{\vartheta}^N}, \overline{H}\in\mathbb{R}} \|u_{\boldsymbol{\alpha}}\|^2_{L^2(\Omega)}+\|m_{\boldsymbol{\beta}}\|^2_{L^2(\Omega)} +\gamma\|\mathcal{P}(u_{\boldsymbol{\alpha}}, m_{\boldsymbol{\beta}}, \overline{H})\|_{L^2(\Omega)}^2+\gamma\|\mathcal{B}(u_{\boldsymbol{\alpha}}, m_{\boldsymbol{\beta}})\|_{L^2(\Omega)}^2.
 	\end{split}
 \end{align*}
A numerical example for a time-dependent planning MFG is shown in Section \ref{secNumerical}.
\end{remark}
The following theorem gives the existence of a solution to \eqref{OptFFProb}. 
\begin{theorem}
\label{exitnd}
Under Assumption \ref{hypPB}, the minimization problem \eqref{OptFFProb} admits a minimizer for any given $\gamma\geq 0$.
\end{theorem}
\begin{proof}
The arguments are the same as in the proof of Theorem \ref{ffexi}. First, we convert \eqref{OptFFProb} into an equivalent minimization problem similar to \eqref{1dctsfnt}.  Then, from the lower semi-continuity and the coercivity of the objective function, we conclude that a minimizer exits. 
\end{proof}
Next, we study the convergence of \eqref{OptFFProb} as $N$ and $\gamma$ go to infinity. 
We do not try to provide the most general convergence result for all MFGs since the spaces in \eqref{defSpG} lack compactness. Hence, we cannot apply the arguments of Theorem  \ref{cvtmfgnd} here. Instead, we prove the validity of our method in concrete setups of interests. In the rest of this subsection, we build the convergence results of the FF method applied to a stationary MFG with a  Lipschitz coupling and a unique smooth solution, see Subsection \ref{secErgodicMFG} for a numerical experiment. We postpone the study of the convergence of the FF method in settings of time-dependent MFGs to future work. 

More precisely, given a smooth function $V:\mathbb{R}^d\mapsto \mathbb{R}$ and a functional $F:C(\mathbb{T}^d)\mapsto C(\mathbb{T}^d)$, we consider the following MFG
\begin{align}
\label{ConGenMFG}
\begin{cases}
	-\Delta u + \frac{|\nabla u|^2}{2} + V(x)  =  F[m](x) + \overline{H}, \ \text{in}\ \mathbb{T}^d,\\
	-\Delta m - \div(m\nabla u) = 0, \ \text{in}\ \mathbb{T}^d,\\
	\int_{\mathbb{T}^d}m \dif x = 1, \ \int_{\mathbb{T}^d}u\dif x = 0. 
\end{cases}
\end{align}
We assume that \eqref{ConGenMFG} admits a smooth solution. In addition, we suppose further that \eqref{ConGenMFG} satisfies the following assumption, which  guarantees the uniqueness of the solution to  \eqref{ConGenMFG}. 
\begin{hyp}
\label{hypF}
There exists a constant $L_F>0$ such that for any $m, \mu \in C(\mathbb{T}^d)$, 
\begin{align}
	\label{lpF}
	\sup_{x\in\mathbb{T}^d}|F[m](x)-F[\mu](x)|\leq L_F\|m-\mu\|_{L^\infty(\mathbb{T}^d)}. 
\end{align}
Meanwhile, $F$ is monotone, i.e., for any $m, \mu\in C(\mathbb{T}^d)$, $m=\mu$ if and only if 
\begin{align*}
\int_{\mathbb{T}^d}(F[m](x)-F[\mu](x))(m(x)-\mu(x))\dif x \leq  0.  
\end{align*}
\end{hyp}
Next, we study the convergence of our method when we use the Fourier series to approximate $u$ and  $m$. Given $N\in \mathbb{N}$, we define the Fourier features space 
\begin{align}
	\label{defGnd}
	\mathcal{G}^{N}=\bigg\{\phi\bigg|\phi(x)=c + \sum_{i\in \mathbb{Z}_N^d}\alpha_i\sin(2\pi i^T x)+\sum_{i\in \mathbb{Z}_N^d}\beta_{i}\cos(2\pi i^T x), c, \alpha_i, \beta_i\in \mathbb{R}\bigg\}
\end{align}
and equip $\mathcal{G}^N$ with the norm $\|\cdot\|_{\mathcal{G}^N}$ defined as
\begin{align*}
\|\phi\|_{\mathcal{G}^N}=|c|^2 + \sum_{i\in \mathbb{Z}^d_N}|\alpha_i|^2+|\beta_i|^2, \forall \phi\in \mathcal{G}^N. 
\end{align*}
For  ease of presentation, given $\gamma>0$, we define the following functional  
\begin{align}
	\label{defcand}
	J_\gamma(u, m, \overline{H})=\int_{\mathbb{T}^d}|u|^2\dif x + \int_{\mathbb{T}^d}|m|^2\dif x+\overline{H}^2+\gamma\mathcal{Q}(u, m, \overline{H}),
\end{align}
where 
\begin{align}
	\label{defmqnd}
	\begin{split}
		\mathcal{Q}(u, m, \overline{H}) =  
		\int_{\mathbb{T}^d}|-\Delta u + \frac{|\nabla u|^2}{2}+V(x)-F[m](x)-\overline{H}|^2\dif x\quad\quad\quad\quad\quad\\
		+\int_{\mathbb{T}^d}|-\Delta m - \div(m\nabla u)|^2\dif x+\biggl |\int_{\mathbb{T}^d}u\dif x\biggl|^2+\biggl|\int_{\mathbb{T}^d}m\dif x - 1\biggl|^2.
	\end{split}
\end{align}
Then, \eqref{OptFFProb} is equivalent to 
\begin{align}
	\label{ndcts}
	\min_{u^N\in \mathcal{G}^N, m\in \mathcal{G}^N, \overline{H}^N\in \mathbb{R}}J_\gamma(u^N, m^N, \overline{H}^N).
\end{align}
Using the same arguments as in the proof of Theorem \ref{exitnd}, we conclude that $\eqref{ndcts}$ admits a minimizer. Next, we show the convergence of a minimizer of \eqref{ndcts}. First, we give an upper bound for the minimum of \eqref{ndcts} in the following theorem. 
\begin{theorem}
	\label{convUpBdFFnd}
	Let $(u^*, m^*, \overline{H}^*)$ be the solution to \eqref{ConGenMFG} and the functional $J_\gamma$ be as in \eqref{defcand} for $\gamma>0$. Suppose that Assumption \ref{hypF} holds. 
	Then, 
	for any sufficiently small $\epsi>0$, there exist a constant $C>0$, $N\in\mathbb{N}$, and functions $(u^N, m^N)\in \mathcal{G}^N\times\mathcal{G}^N$ such that 
	\begin{align}
		\label{bdJnd}
		J_\gamma(u^N, m^N, \overline{H}^*) \leq 2\|u^*\|_{L^2(\mathbb{T}^d)}^2+2\|m^*\|_{L^2(\mathbb{T}^d)}^2+|\overline{H}^*|^2+C(1+\gamma)\epsilon^2.
	\end{align}
\end{theorem}
See Appendix \ref{appendix} for  the proof of the above theorem. The following corollary follows directly from Theorem \ref{convUpBdFFnd} and proves that there exists  a minimizer $({u}^{N, \gamma}, {m}^{N, \gamma}, {\overline{H}}^{N, \gamma})$ of \eqref{ndcts} such that $\mathcal{Q}({u}^{N, \gamma}, {m}^{N, \gamma}, {\overline{H}}^{N, \gamma})\rightarrow 0$ as $N$ and $\gamma$ go to infinity. 
\begin{corollary}
	\label{convFFQnd} Let $\mathcal{Q}$ be as in \eqref{defmqnd}. Under the same assumptions of Theorem \ref{convUpBdFFnd},  
	for any $\epsilon>0$, there exist a constant $C>0$,  sufficiently large $\gamma>0$ and  $N>0$, and a minimizer of \eqref{ndcts}, which is denoted by $({u}^{N,\gamma}, {m}^{N,\gamma}, {\overline{H}}^{N,\gamma})\in \mathcal{G}^N\times\mathcal{G}^N\times \mathbb{R}$, satisfies
	\begin{align*}
		\mathcal{Q}({u}^{N, \gamma}, {m}^{N, \gamma}, {\overline{H}}^{N, \gamma})\leq C\epsilon. 
	\end{align*}
\end{corollary}
We give the proof of Corollary \ref{convFFQnd} in Appendix \ref{appendix}. The following theorem shows the convergence of minimizers of \eqref{ndcts}  to the solution of \eqref{ConGenMFG} as $N$ and $\gamma$ go to infinity. The proof is presented in Appendix \ref{appendix}.
\begin{theorem}
	\label{convndFF}Let $(u^*, m^*, \overline{H}^*)$ be the solution to \eqref{ConGenMFG}. Under the assumptions of Theorem \ref{convUpBdFFnd},  
	there exists a sequence $\{(N_i,\gamma_i)\}_{i=1}^{\infty}$ such that the sequence $\{(u^{i}, m^{i}, \overline{H}^{i})\}_{i=1}^\infty$, where $(u^{i}, m^{i}, \overline{H}^{i})$ is a minimizer of  \eqref{ndcts} corresponding to $N_i$ and $\gamma_i$, satisfies $u^{i}\rightarrow u^*$ in $H^{1}(\mathbb{T}^d)$, $m^{i}\rightarrow m^*$ in $H^{1}(\mathbb{T}^d)$, and $\overline{H}^{i}\rightarrow \overline{H}^*$ in $\mathbb{R}$ as $i\rightarrow\infty$. 
\end{theorem}

Next, we propose a method to solve \eqref{OptFFProb}. We  take $M$ samples $\{x_i\}_{i=1}^M$ in the domain $\Omega$ such that $\{x_i\}_{i=1}^{M_\Omega}\subset\operatorname{int}\Omega$ and $\{x_i\}_{i=M_\Omega+1}^M\subset\partial\Omega$ for $1\leq M_\Omega\leq M$. To capture different variability of the constraints, we introduce two penalization parameters $\gamma$, $\beta$, and consider 
\begin{align}
	\label{OptFFProbDisc}
	\begin{split}
		\min_{u_{\boldsymbol{\alpha}}\in\mathcal{G}_{\boldsymbol{\zeta}^N}, m_{\boldsymbol{\beta}}\in\mathcal{G}_{\boldsymbol{\vartheta}^N}, \overline{H}\in\mathbb{R}} \|u_{\boldsymbol{\alpha}}\|^2_{\mathcal{G}_{\boldsymbol{\zeta}^N}}+\|m_{\boldsymbol{\beta}}\|^2_{\mathcal{G}_{\boldsymbol{\vartheta}^N}}+|\overline{H}|^2+\beta\bigg|\frac{1}{M_\Omega}\sum_{i=1}^{M_\Omega}m_{\boldsymbol{\beta}}(x_i)-1\bigg|^2\quad\quad\quad \\+\beta\bigg|\frac{1}{M_\Omega}\sum_{i=1}^{M_\Omega}u_{\boldsymbol{\alpha}}(x_i)\bigg|^2+\gamma\sum_{i=1}^{M_\Omega}|\mathcal{P}(u_{\boldsymbol{\alpha}}, m_{\boldsymbol{\beta}}, \overline{H})(x_i)|^2
		+\gamma\sum_{i=M_{\Omega}+1}^{M}|\mathcal{B}(u_{\boldsymbol{\alpha}}, m_{\boldsymbol{\beta}})(x_i)|^2.
	\end{split}
\end{align}
Let $N_\mathcal{U}$ and $N_\mathcal{V}$ be as in  \eqref{vecPhi} and \eqref{vecPsi}. Under Assumption \ref{hypPB}, we reformulate \eqref{OptFFProbDisc} into an equivalent two-level minimization problem
\begin{align}
	\label{OptFFProb2l}
	\begin{split}
	\min\limits_{\boldsymbol{z}\in\mathbb{R}^{N_\mathcal{U}}, \boldsymbol{\rho}\in\mathbb{R}^{N_\mathcal{V}}, \lambda\in\mathbb{R}} \gamma\|G(\boldsymbol{z}, \boldsymbol{\rho}, \lambda)\|^2+\beta\bigg|\frac{1}{M_\Omega}\sum_{i=1}^{M_\Omega}\rho_i^{(1)}-1\bigg|^2+\beta\bigg|\frac{1}{M_\Omega}\sum_{i=1}^{M_\Omega}z_i^{(1)}\bigg|^2\\
	+
	\begin{cases}
		\min\limits_{\boldsymbol{\alpha}\in\mathbb{R}^N, \boldsymbol{\beta}\in\mathbb{R}^N, \overline{H}\in\mathbb{R}} \|\boldsymbol{\alpha}\|^2+\|\boldsymbol{\beta}\|^2+\overline{H}^2,\\
		\text{s.t.}\ {[\boldsymbol{\phi}, \boldsymbol{\zeta}^N]}\boldsymbol{\alpha} = \boldsymbol{z},  {[\boldsymbol{\psi},  \boldsymbol{\vartheta}^N]}\boldsymbol{\beta} = \boldsymbol{\rho}, \overline{H} = \lambda.
	\end{cases}
	\end{split}
\end{align} 
The first level optimization problem gives
\begin{align*}
\begin{cases}
\boldsymbol{\alpha}=[\boldsymbol{\phi}, \boldsymbol{\zeta}^N]^T([\boldsymbol{\phi}, \boldsymbol{\zeta}^N][\boldsymbol{\phi}, \boldsymbol{\zeta}^N]^T)^{-1}\boldsymbol{z}, \\
\boldsymbol{\beta}=[\boldsymbol{\psi}, \boldsymbol{\vartheta}^N]^T([\boldsymbol{\psi}, \boldsymbol{\vartheta}^N][\boldsymbol{\psi}, \boldsymbol{\vartheta}^N]^T)^{-1}\boldsymbol{\rho},\\ \overline{H}=\lambda. 
\end{cases}
\end{align*}
Hence, \eqref{OptFFProb2l} is equivalent to
\begin{align*}
\begin{split}
\min\limits_{\boldsymbol{z}\in \mathbb{R}^{N_\mathcal{U}}, \boldsymbol{\rho}\in \mathbb{R}^{N_\mathcal{V}}, \lambda\in\mathbb{R}} \gamma\|G(\boldsymbol{z}, \boldsymbol{\rho}, \lambda)\|^2+\beta\bigg|\frac{1}{M_\Omega}\sum_{i=1}^{M_\Omega}\rho_i^{(1)}-1\bigg|^2+\beta\bigg|\frac{1}{M_\Omega}\sum_{i=1}^{M_\Omega}z_i^{(1)}\bigg|^2\\
+\boldsymbol{z}^T([\boldsymbol{\phi}, \boldsymbol{\zeta}^N][\boldsymbol{\phi}, \boldsymbol{\zeta}^N]^T)^{-1}\boldsymbol{z}+\boldsymbol{\beta}^T([\boldsymbol{\psi}, \boldsymbol{\vartheta}^N][\boldsymbol{\psi}, \boldsymbol{\vartheta}^N]^T)^{-1}\boldsymbol{\beta}+\lambda^2.
\end{split}
\end{align*}

\begin{remark}
\label{ffreg}
In general, the matrices $[\boldsymbol{\phi}, \boldsymbol{\zeta}^N][\boldsymbol{\phi}, \boldsymbol{\zeta}^N]^T$ and $[\boldsymbol{\psi}, \boldsymbol{\vartheta}^N][\boldsymbol{\psi}, \boldsymbol{\vartheta}^N]^T$ are ill-conditioned. Hence, we choose two regularization parameters, $\mu_1>0$ and $\mu_2>0$, and consider
\begin{align}
	\label{OptFFProbrfreg}
	\begin{split}
	\min\limits_{\boldsymbol{z}\in \mathbb{R}^{N_\mathcal{U}}, \boldsymbol{\rho}\in \mathbb{R}^{N_\mathcal{V}}, \lambda\in\mathbb{R}} \gamma\|G(\boldsymbol{z}, \boldsymbol{\rho}, \lambda)\|^2+\beta\bigg|\frac{1}{M_\Omega}\sum_{i=1}^{M_\Omega}\rho_i^{(1)}-1\bigg|^2+\beta\bigg|\frac{1}{M_\Omega}\sum_{i=1}^{M_\Omega}z_i^{(1)}\bigg|^2\\
	\quad+\boldsymbol{z}^T([\boldsymbol{\phi}, \boldsymbol{\zeta}^N][\boldsymbol{\phi}, \boldsymbol{\zeta}^N]^T+\mu_1 I)^{-1}\boldsymbol{z}\quad\quad\quad\quad\quad\quad\quad\quad\quad\quad\\
	+\boldsymbol{\beta}^T([\boldsymbol{\psi}, \boldsymbol{\vartheta}^N][\boldsymbol{\psi}, \boldsymbol{\vartheta}^N]^T+\mu_2 I)^{-1}\boldsymbol{\beta}+\lambda^2.\quad\quad\quad\quad\quad\quad\quad 
	\end{split}
\end{align}
\end{remark}
\begin{remark}
When necessary, we also eliminate equality constraints in \eqref{OptFFProbrfreg} as discussed in Section 3.3 of \cite{chen2021solving}.
\end{remark}
{
\begin{remark}
According to Theorems \ref{cvtmfgnd} and \ref{convndFF}, using the GP and the FF methods, the non-negativity of the probability measure $m$ is guaranteed at the limit. Thus, unless the coupling term is not well defined when $m$ is non-positive (see Remark \ref{rmk1dpositi} for discussions about a MFG with a log coupling), we do not impose extra non-negativity constraints on the Gauss--Newton iterations. We will also see that the numerical results of the probability measures are non-negative in the next section. 
\end{remark}
}

\section{Numerical Results}
\label{secNumerical}
\begin{comment}
All the graphs and errors in the following experiments are generated on a workstation running Ubuntu Linux. The system has Intel(R)
Xeon(R) Gold 5218 CPUs with a total of 32 cores and 377 GB of
memory. Since the system is shared with other users, we have to use 
\end{comment}
In this section, we implement our methods to solve a non-local stationary MFG in Subsection \ref{secErgodicMFG} and a time-dependent planning problem in Subsection \ref{secPlanning}. The runtimes are measured using MacBook Air 2015 (4GB Ram, Intel Core i5 CPU). Our implementation is based on the code of \cite{chen2021solving}\footnote{https://github.com/yifanc96/NonlinearPDEs-GPsolver.git}, which uses Python with the JAX package for automatic differentiation. Our experiments count only on CPUs. Additional speedups can be achieved by using accelerated hardware such as
Graphics Processing Units (GPU).
\subsection{A Second-Order Non-local Stationary MFG}
\label{secErgodicMFG}
We consider a variant of the non-local stationary MFG in Subsection 6.2.5 of \cite{achdou2010mean}. More precisely, given $\nu>0$, we want to find $(u, m, \overline{H})$ solving 
\begin{align}
\label{ErgodicProb}
\begin{cases}
-\nu\Delta u + H(x, \nabla u)  =  \big((1-\Delta)^{-1}(1-\Delta)^{-1}m\big)(x) + \overline{H}, \ \text{in}\ \mathbb{T}^2,\\
-\nu\Delta m - \div(mD_pH(x,\nabla u)) = 0, \ \text{in}\ \mathbb{T}^2,\\
\int_{\mathbb{T}^2}u\dif x = 0, \int_{\mathbb{T}^2}m\dif x = 1,
\end{cases}
\end{align}
where 
\begin{align*}
	H(x, p)=\sin(2\pi x_1)+\sin(2\pi x_2)+\cos(4\pi x_1)+|p|^{2}, \forall x\in \mathbb{T}^2, p\in \mathbb{R}^2. 
\end{align*}
We use the GP method and the FF algorithm proposed in Section \ref{secGF} to  solve \eqref{ErgodicProb}. In the experiments, we write \eqref{ErgodicProb} in the form of \eqref{gsmfgslo} with $Q_b=D_b=0$, $Q=4$, $D=5$, and with linear operators $L_1(u)=u$, $L_2(u)=\partial_xu$, $L_3(u)=\partial_yu$, $L_4(u)=\Delta u$, $J_1(m)=m$, $J_2(m)=\partial_xm$, $J_3(m)=\partial_ym$, $J_4(m)=\Delta m$, and $J_5(m)=\big((1-\Delta)^{-1}(1-\Delta)^{-1}m\big)$. We compute the action of $J_5$ on a kernel by the Fast Fourier transform. The GP method uses the periodic kernels
\begin{align*}
K_1((x_1, x_2), (y_1, y_2))=K_2(x_1, x_2, y_1, y_2)=e^{\frac{1}{\sigma^2}(\cos(2\pi(x_1-y_1))+\cos(2\pi(x_2-y_2))-2)}
\end{align*}
with lengthscale $\sigma=0.2$. 
For the FF method, we fix $N\in\mathbb{N}$, use the basis 
\begin{align*}
\boldsymbol{\zeta}^N(x_1,x_2)=&\boldsymbol{\vartheta}^N(x_1, x_2)=[1, \sin(2\pi x_1+2\pi x_2),\dots,\sin(2\pi i x_1+2\pi j x_2),\\
&\dots, \sin(2\pi N x_1+2\pi N x_2),\cos(2\pi x_1+2\pi x_2),\\
&\dots,\cos(2\pi i x_1 + 2\pi jx_2),\dots, \cos(2\pi N x_1 + 2\pi N x_2)]^T,  \forall i, j=1,\dots, N,
\end{align*}
and approximate $u$ and $m$ by functions in the space $\mathcal{G}_{\boldsymbol{\zeta}^N}$ as in \eqref{defSpG}. We choose the regularization parameters $\eta_1=\eta_2=\mu_1=\mu_2=10^{-5}$ in Remarks \ref{rmk1dnugget} and \ref{regFF}. To measure the accuracy of our algorithms, we identify $\mathbb{T}^2$ with $[0, 1)\times [0, 1)$, discretize the domain with $100\times 100$ uniformly distributed grid points, and use the FD method in \cite{achdou2010mean} to solve \eqref{ErgodicProb} with high accuracy. The GP method and the FF algorithm use the same sample points and start from the same initial values. We denote by $(u_{GP}, m_{GP}, \overline{H}_{GP})$ and  $(u_{FF}, m_{FF}, \overline{H}_{FF})$ the numerical solutions of the GP method and the FF algorithm, separately. 

\begin{figure}[!hbtbp]
	\centering
	\begin{tabular}{c c c}
		\begin{subfigure}[b]{0.3\textwidth} 
			\centering
			\resizebox{\textwidth}{!}{%
				\begin{tabular}{c c}
					\hline
					$|\overline{H}_{GP}-\overline{H}_{FD}|$ &  	 $|\overline{H}_{FF}-\overline{H}_{FD}|$ \\
					\hline
					$0.08$ & $0.08$  \\
					\hline
				\end{tabular}%
			}
			\caption{Errors of $\overline{H}_{GP}$ and $\overline{H}_{FF}$}
			\label{fig:Ergodicnu01:3}
		\end{subfigure}&              
		\begin{subfigure}[b]{0.3\textwidth}
			\includegraphics[width=\textwidth]{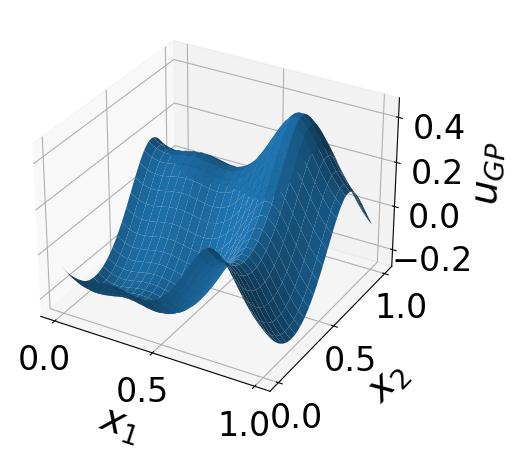}
			\caption{$u_{GP}$.}
			\label{fig:Ergodicnu01:1}
		\end{subfigure} &
		\begin{subfigure}[b]{0.3\textwidth}
			\includegraphics[width=\textwidth]{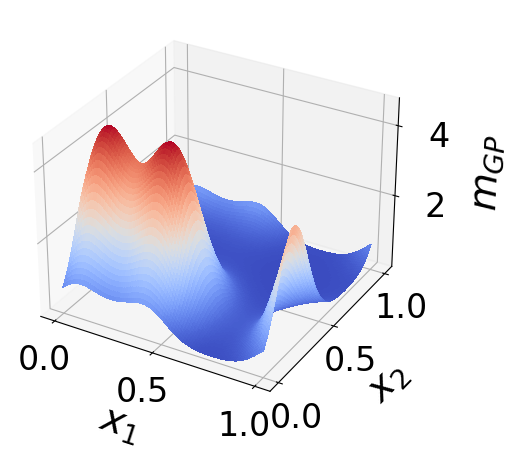}
			\caption{$m_{GP}$}
			\label{fig:Ergodicnu01:2}
		\end{subfigure}\\ 
		\begin{subfigure}[b]{0.3\textwidth}
			\includegraphics[width=\textwidth]{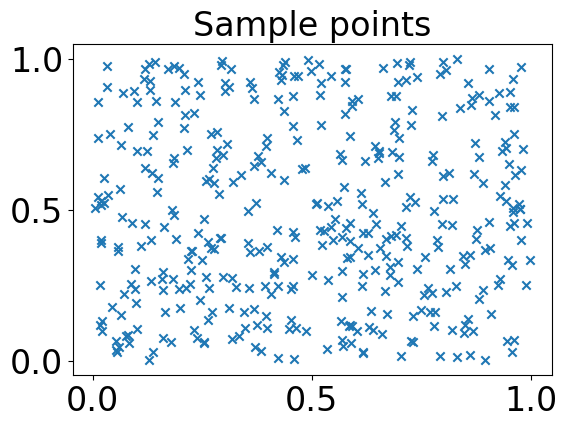}
			\caption{Samples}
			\label{fig:Ergodicnu01:0}
		\end{subfigure} 
		&
		\begin{subfigure}[b]{0.3\textwidth}
			\includegraphics[width=\textwidth]{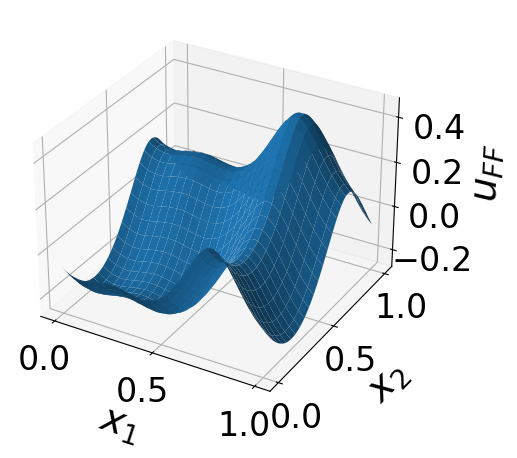}
			\caption{$u_{FF}$}
			\label{fig:Ergodicnu01:4}
		\end{subfigure} &
		\begin{subfigure}[b]{0.3\textwidth}
			\includegraphics[width=\textwidth]{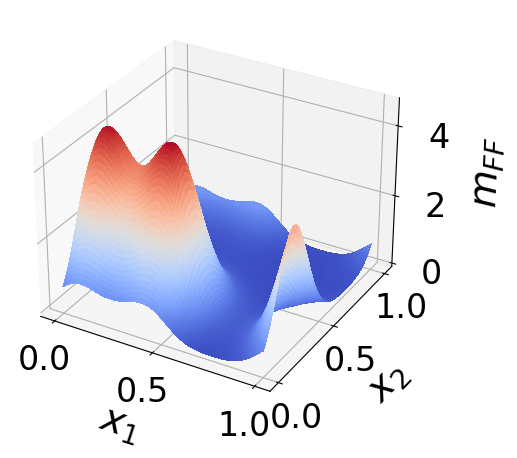}
			\caption{$m_{FF}$}
			\label{fig:Ergodicnu01:5}
		\end{subfigure}\\
		\begin{subfigure}[b]{0.3\textwidth}
			\includegraphics[width=\textwidth]{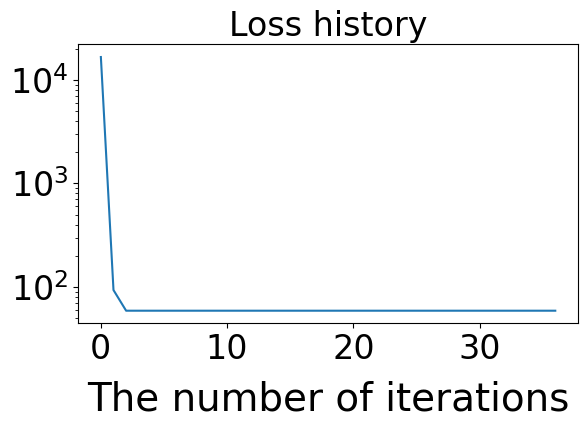}
			\caption{Loss history (GP)}
			\label{fig:Ergodicnu01:6} 
		\end{subfigure}&
		\begin{subfigure}[b]{0.3\textwidth}
			\includegraphics[width=\textwidth]{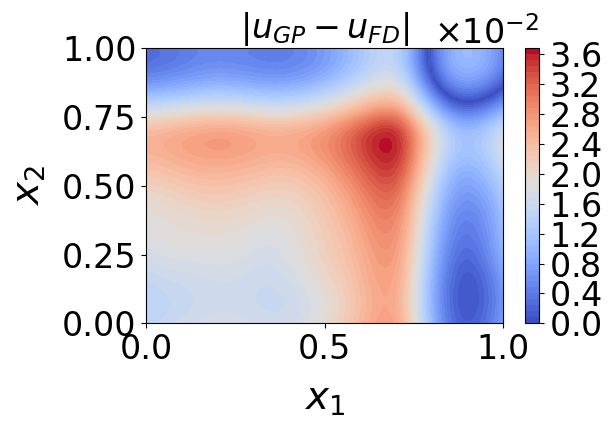}
			\caption{$|u_{GP}-u_{FD}|$}
			\label{fig:Ergodicnu01:7} 
		\end{subfigure} & 
		\begin{subfigure}[b]{0.3\textwidth}
			\includegraphics[width=\textwidth]{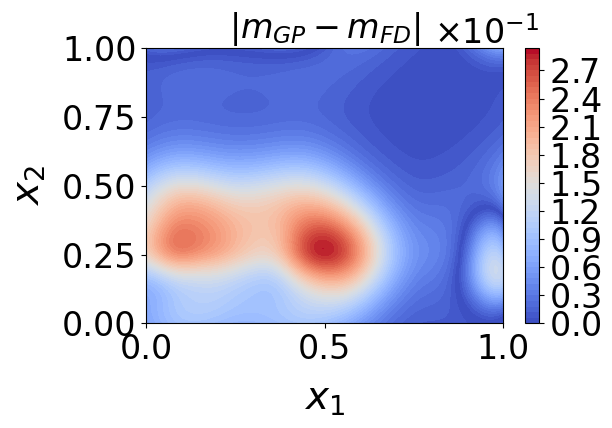}
			\caption{$ |m_{GP}-u_{FD}|$}
			\label{fig:Ergodicnu01:8}
		\end{subfigure}\\
		\begin{subfigure}[b]{0.3\textwidth}
			\includegraphics[width=\textwidth]{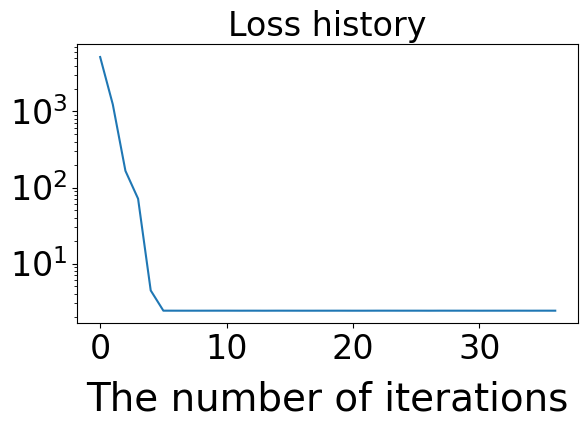}
			\caption{Loss history (FF)}
			\label{fig:Ergodicnu01:9} 
		\end{subfigure}&
		\begin{subfigure}[b]{0.3\textwidth}
			\includegraphics[width=\textwidth]{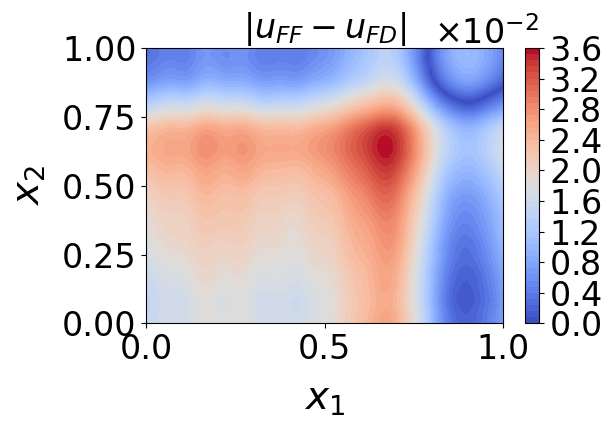}
			\caption{$|u_{FF}-u_{FD}|$}
			\label{fig:Ergodicnu01:10} 
		\end{subfigure} & 
		\begin{subfigure}[b]{0.3\textwidth}
			\includegraphics[width=\textwidth]{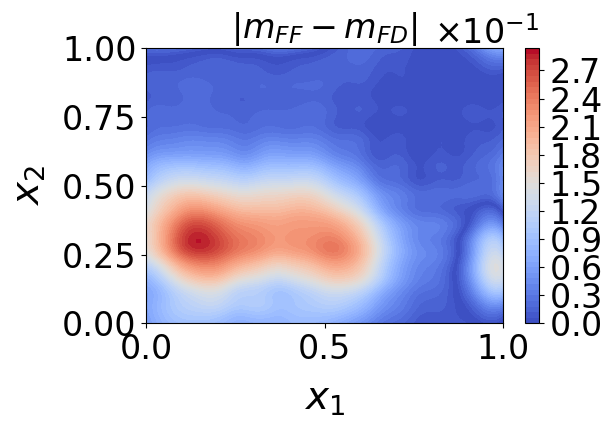}
			\caption{$ |m_{FF}-u_{FD}|$}
			\label{fig:Ergodicnu01:11}
		\end{subfigure}
	\end{tabular}
	\caption{The numerical solutions of the non-local stationary MFG \eqref{ErgodicProb} when $\nu=0.1$. (a) The errors of $\overline{H}_{GP}$ and $\overline{H}_{FF}$ with respect to the reference $\overline{H}_{FD}$; (b), (c), (e), (f) Numerical values of $u_{GP}$, $m_{GP}$, $u_{FF}$, and $m_{FF}$; (d) the samples used by the GP and the FF methods; (g), (j) the evolution of the losses, which are the objective functions in \eqref{gpfin} and \eqref{OptFFProbrfreg}; (h), (i), (k), (l) The contours of the pointwise errors of $u_{GP}$,   $m_{GP}$, $u_{FF}$, and $m_{FF}$. }
	\label{fig:ErgodicN400Nu01}
\end{figure}

In Figure \ref{fig:ErgodicN400Nu01}, we show the numerical results of both algorithms for $\nu=0.1$ after {36} iterations.  We take the same $M=400$ samples for both the GP method and the FF algorithm, which is shown in Figure \ref{fig:Ergodicnu01:0}. We use the Gauss--Newton iteration with step size $1$ to solve the optimization problems, and choose $N=10$ for the FF method. We set the penalization parameters $\beta=1$ and $\gamma=10^{-15}$ for both algorithms.  Figures \ref{fig:Ergodicnu01:1}, \ref{fig:Ergodicnu01:4},  \ref{fig:Ergodicnu01:2}, and \ref{fig:Ergodicnu01:5} plot the graphs of $u_{GP}$, $u_{FF}$, $m_{GP}$, and $m_{FF}$, separately. The convergence histories of the Gauss--Newton iterations are presented in Figures \ref{fig:Ergodicnu01:6} and \ref{fig:Ergodicnu01:9}, which verify the convergence of our algorithms. We attribute the non-monotone decreasing of the loss curves to the non-linearity of the objective functions. The contours of pointwise errors are shown in Figures \ref{fig:Ergodicnu01:7}, \ref{fig:Ergodicnu01:8}, \ref{fig:Ergodicnu01:10}, and \ref{fig:Ergodicnu01:11}. We see that the errors
are smaller in smoother areas. 
The errors of $\overline{H}_{GP}$ and $\overline{H}_{FF}$ are given in Table \ref{fig:Ergodicnu01:3}. 

\begin{figure}[!hbtbp]
	\centering
	\begin{tabular}{c c c}
		\begin{subfigure}[b]{0.3\textwidth} 
			\centering
			\resizebox{\textwidth}{!}{%
				\begin{tabular}{c c}
					\hline
					$|\overline{H}_{GP}-\overline{H}_{FD}|$ &  	 $|\overline{H}_{FF}-\overline{H}_{FD}|$ \\
					\hline
					$1.34\times 10^{-2}$ & $4.96 \times 10^{-3}$  \\
					\hline
				\end{tabular}%
			}
			\caption{Errors of $\overline{H}_{GP}$ and $\overline{H}_{FF}$}
			\label{fig:Ergodicnu1:3}
		\end{subfigure}&              
	\begin{subfigure}[b]{0.3\textwidth}
		\includegraphics[width=\textwidth]{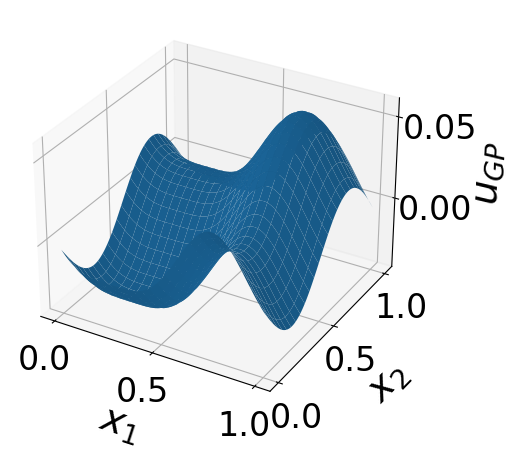}
		\caption{$u_{GP}$.}
		\label{fig:Ergodicnu1:1}
	\end{subfigure} &
	\begin{subfigure}[b]{0.3\textwidth}
		\includegraphics[width=\textwidth]{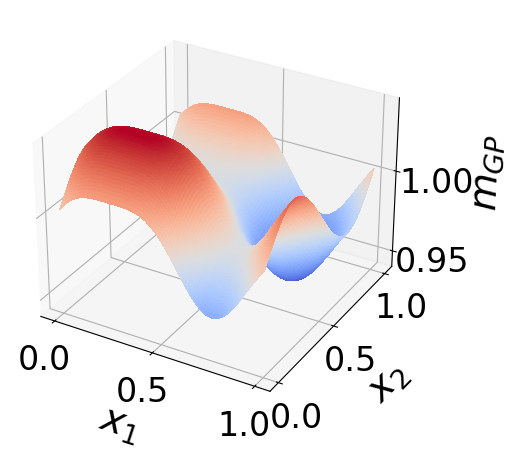}
		\caption{$m_{GP}$}
		\label{fig:Ergodicnu1:2}
	\end{subfigure}\\ 
\begin{subfigure}[b]{0.3\textwidth}
	\includegraphics[width=\textwidth]{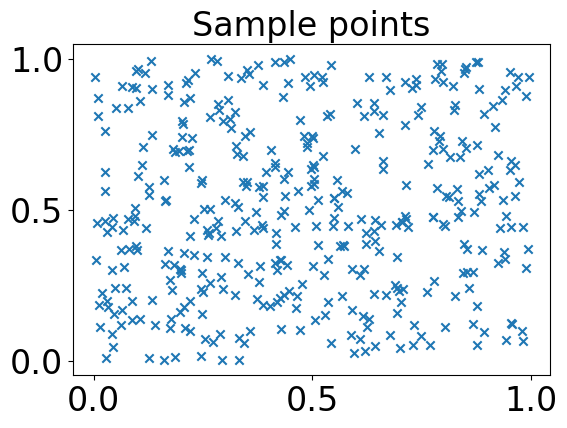}
	\caption{Samples}
	\label{fig:Ergodicnu1:0}
\end{subfigure} 
&
\begin{subfigure}[b]{0.3\textwidth}
	\includegraphics[width=\textwidth]{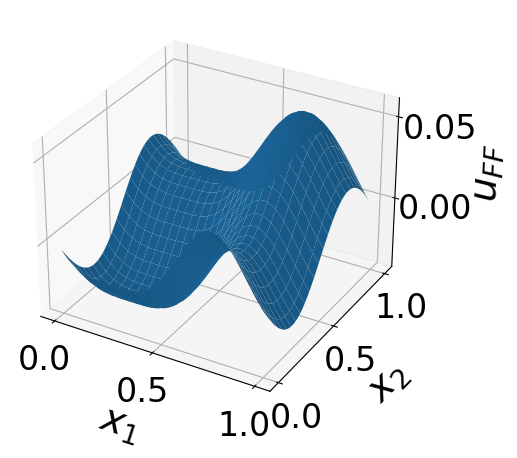}
	\caption{$u_{FF}$}
	\label{fig:Ergodicnu1:4}
\end{subfigure} &
\begin{subfigure}[b]{0.3\textwidth}
	\includegraphics[width=\textwidth]{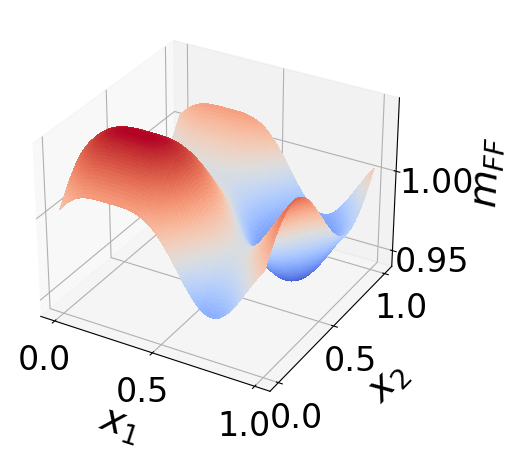}
	\caption{$m_{FF}$}
	\label{fig:Ergodicnu1:5}
\end{subfigure}\\
	\begin{subfigure}[b]{0.3\textwidth}
	\includegraphics[width=\textwidth]{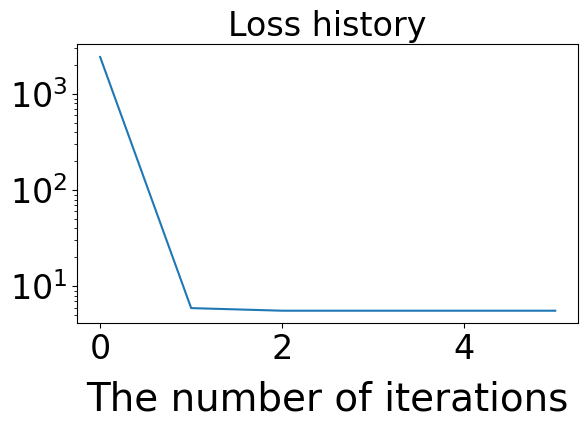}
	\caption{Loss history (GP)}
	\label{fig:Ergodicnu1:6} 
\end{subfigure}&
\begin{subfigure}[b]{0.3\textwidth}
	\includegraphics[width=\textwidth]{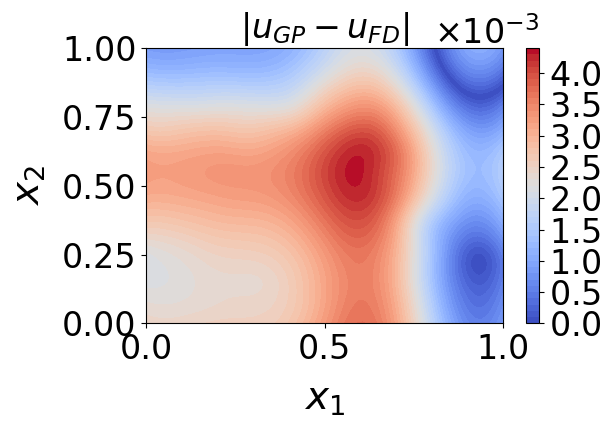}
	\caption{$|u_{GP}-u_{FD}|$}
	\label{fig:Ergodicnu1:7} 
\end{subfigure} & 
\begin{subfigure}[b]{0.3\textwidth}
	\includegraphics[width=\textwidth]{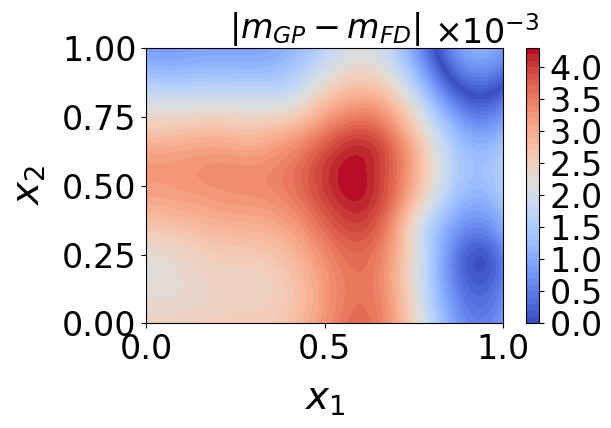}
	\caption{$ |m_{GP}-u_{FD}|$}
	\label{fig:Ergodicnu1:8}
\end{subfigure}\\
	\begin{subfigure}[b]{0.3\textwidth}
		\includegraphics[width=\textwidth]{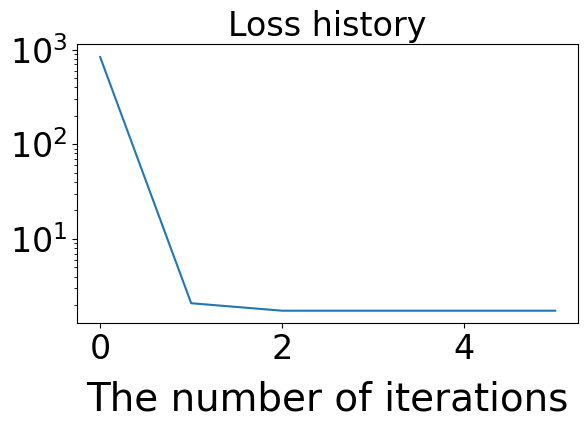}
		\caption{Loss history (FF)}
		\label{fig:Ergodicnu1:9} 
	\end{subfigure}&
		\begin{subfigure}[b]{0.3\textwidth}
			\includegraphics[width=\textwidth]{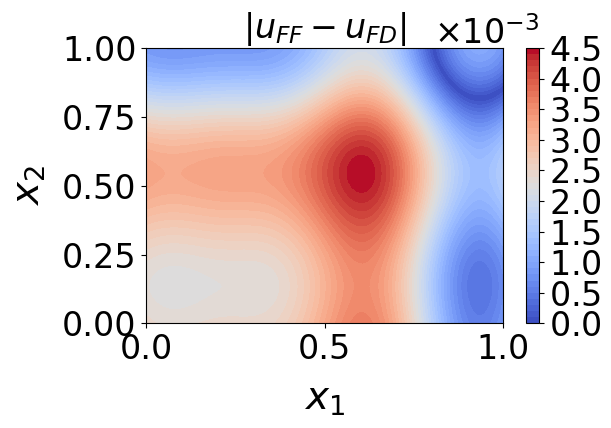}
			\caption{$|u_{FF}-u_{FD}|$}
			\label{fig:Ergodicnu1:10} 
		\end{subfigure} & 
		\begin{subfigure}[b]{0.3\textwidth}
			\includegraphics[width=\textwidth]{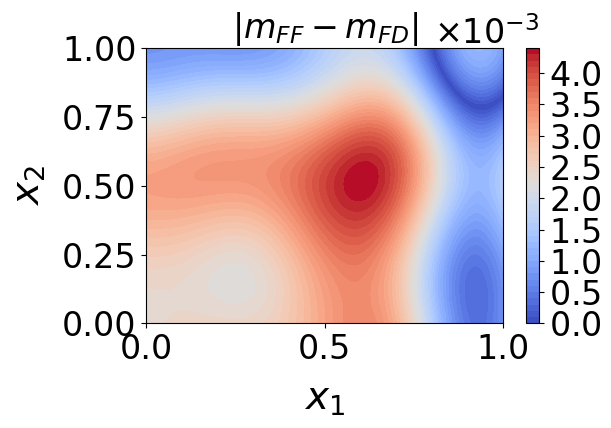}
			\caption{$ |m_{FF}-u_{FD}|$}
			\label{fig:Ergodicnu1:11}
		\end{subfigure}
	\end{tabular}
	\caption{The numerical solutions of the non-local stationary MFG \eqref{ErgodicProb} when $\nu=1$. (a) The errors of $\overline{H}_{GP}$ and $\overline{H}_{FF}$ with respect to the reference $\overline{H}_{FD}$; (b), (c), (e), (f) Numerical values of $u_{GP}$, $m_{GP}$, $u_{FF}$, and $m_{FF}$; (d) the samples used by the GP and the FF methods; (g), (j) the evolution of the losses, which are the objective functions in \eqref{gpfin} and \eqref{OptFFProbrfreg}; (h), (i), (k), (l) The contours of the pointwise errors of $u_{GP}$,   $m_{GP}$, $u_{FF}$, and $m_{FF}$. }
	\label{fig:ErgodicN400Nu1}
\end{figure}

For larger values of $\nu$, the solutions are smoother. Then, fewer bases are enough  for the FF method to achieve higher accuracy, which is shown in Figure \ref{fig:ErgodicN400Nu1}. We take $M=400$, $\nu=1$, and $N=2$. Meanwhile, we set the penalization parameters $\beta=1$, $\gamma=10^{-15}$ for the FF method and $\beta=10^{-2}$, $\gamma=10^{-15}$ for the GP algorithm. The Gauss--Newton method uses step size $1$ for both methods and stops after {5} iterations. Figures \ref{fig:ErgodicN400Nu01} and \ref{fig:ErgodicN400Nu1}  imply that the selection of parameters depends on the data of the model and suggest the need to study hyperparameter learning in future work. 

Table \ref{ergodicPcnu1} records the CPU time of performing the Cholesky and the QR decomposition for both algorithms as $M$ increases. We set  $\nu=0.1$ and $N=10$. We see that the FF algorithm outperforms the GP method in the precomputation stage. 

\begin{table}[ht!]
	\centering
		\begin{tabular}{c c c c c}
			\hline
			& &  &  &      \\
			& \multicolumn{1}{c}{GP} & & \multicolumn{2}{c}{FF}  \\
			\cline{2-2} \cline{4-5}
			&  &  &   &  \\
			$M$    & Cholesky & & QR & Cholesky \\
			$400$ & $1.29$ & & $0.35$ & $0.31$\\
			$800$ & $6.31$ & & $0.50$ & $0.21$\\
			$1200$ & $38.66$ &  & $1.09$ & $0.23$ \\
			$1600$ & $88.79$ & & $1.92$ & $0.30$ \\
			\hline
	\end{tabular}
	\caption{The precomputation time (in seconds) consumed by the GP method and the FF algorithm for solving \eqref{ErgodicProb} when $\nu=0.1$ and $N=10$ as $M$ increases.}
	\label{ergodicPcnu1}
	\label{sndopcp5}
\end{table}

\begin{comment}
{\color{blue}
\eqref{ErgodicProb} is equivalent to 
\begin{align*}
\begin{cases}
-\nu u_{xx} - \nu u_{yy} + V(x_1,x_2) + (u_x^2+u_y^2)^{3/4}=F[m](x)+\overline{H},\\
-\nu m_{xx} - \nu m_{yy} + \frac{3}{4}(u_x^2+u_y^2)^{-\frac{5}{4}}(u_xu_{xx}+u_yu_{xy})u_xm\\
\quad\quad\quad\quad\quad\quad\quad -\frac{3}{2}(u_x^2+u_y^2)^{-\frac{1}{4}}u_{xx}m-\frac{3}{2}(u_x^2+u_y^2)^{-\frac{1}{4}}u_xm_x\\
\quad\quad\quad\quad\quad\quad\quad+\frac{3}{4}(u_x^2+u_y^2)^{-\frac{5}{4}}(u_xu_{xy}+u_yu_{yy})u_ym\\
\quad\quad\quad\quad\quad\quad\quad -\frac{3}{2}(u_x^2+u_y^2)^{-\frac{1}{4}}u_{yy}m-\frac{3}{2}(u_x^2+u_y^2)^{-\frac{1}{4}}u_ym_y=0.
\end{cases}
\end{align*}
}
\end{comment}
\subsection{A Planning Problem}
\label{secPlanning}
\begin{figure}[!hbtbp]
	\centering          
		\begin{subfigure}[b]{0.3\textwidth}
			\includegraphics[width=\textwidth]{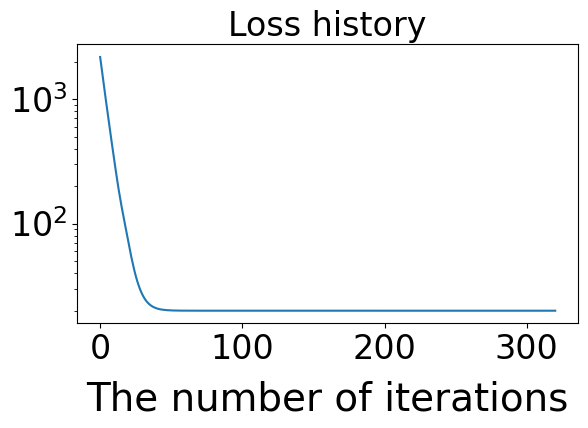}
			\caption{Loss history (GP).}
			\label{fig:Planning:GPloss}
		\end{subfigure} 
		\begin{subfigure}[b]{0.3\textwidth}
			\includegraphics[width=\textwidth]{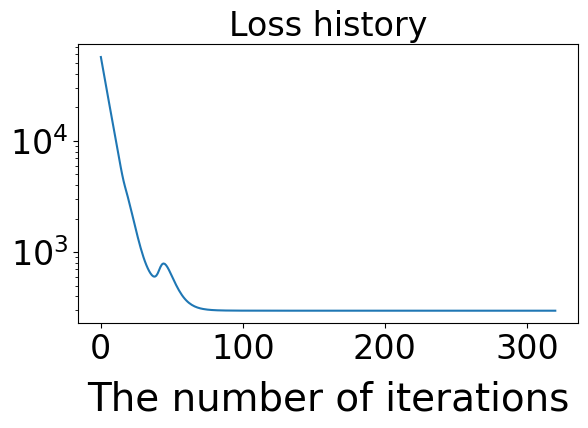}
			\caption{Loss history (FF)}
			\label{fig:Planning:FFloss}
		\end{subfigure}\\ 
		\begin{subfigure}[b]{0.3\textwidth}
			\includegraphics[width=\textwidth]{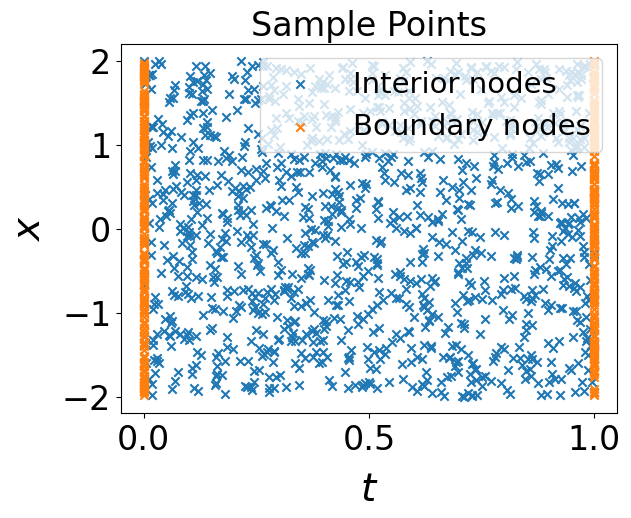}
			\caption{Samples}
			\label{fig:Planning:samples}
		\end{subfigure} 
		\begin{subfigure}[b]{0.3\textwidth}
			\includegraphics[width=\textwidth]{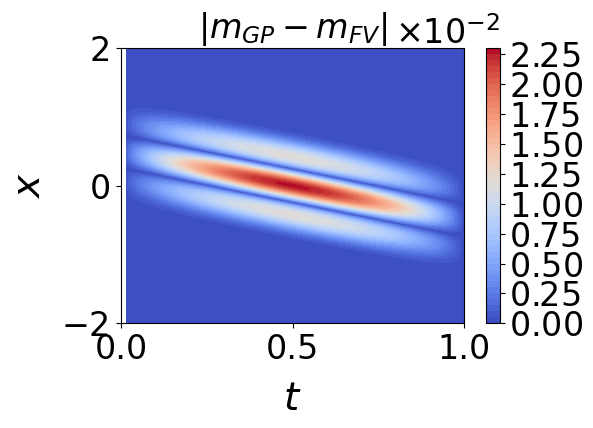}
			\caption{$|m_{GP}-m_{FV}|$}
			\label{fig:Planning:mGPerrors}
		\end{subfigure} 
		\begin{subfigure}[b]{0.3\textwidth}
			\includegraphics[width=\textwidth]{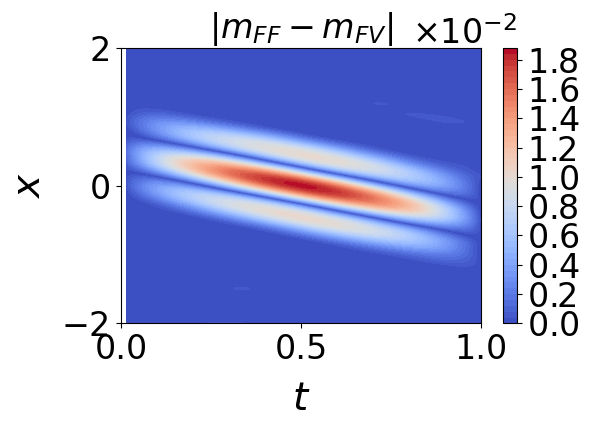}
			\caption{$|m_{FF}-m_{FV}|$}
			\label{fig:Planning:mFFerrors}
		\end{subfigure}\\
		\begin{subfigure}[b]{0.3\textwidth}
			\includegraphics[width=\textwidth]{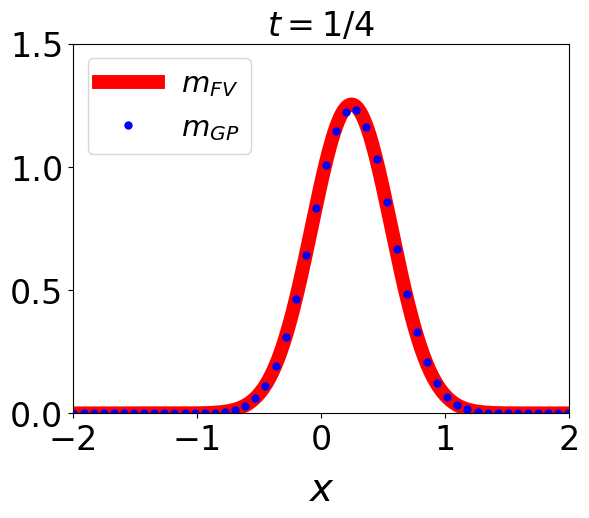}
			\caption{$m_{GP}$ v.s. $m_{FV}$ at $t=1/4$.}
			\label{fig:Planning:m14}
		\end{subfigure}
		\begin{subfigure}[b]{0.3\textwidth}
			\includegraphics[width=\textwidth]{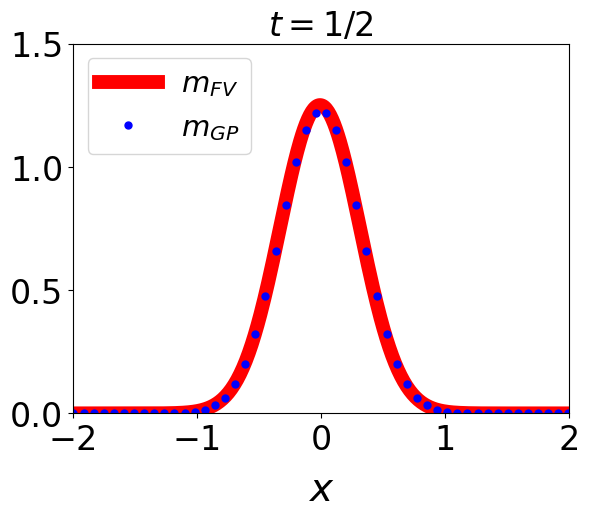}
			\caption{$m_{GP}$ v.s. $m_{FV}$ at $t=1/2$.}
			\label{fig:Planning:m24}
		\end{subfigure} 
		\begin{subfigure}[b]{0.3\textwidth}
			\includegraphics[width=\textwidth]{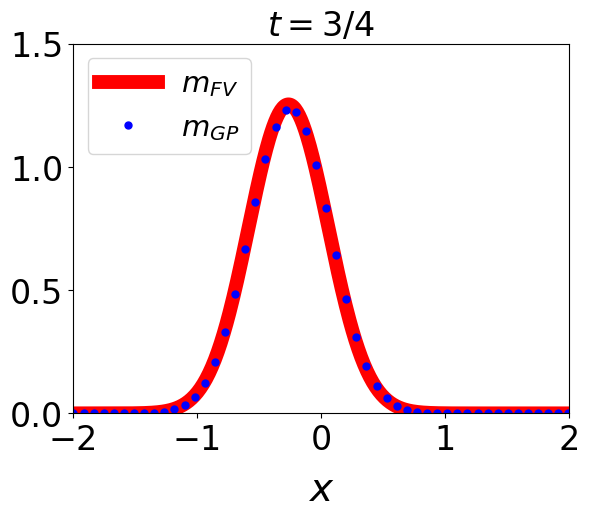}
			\caption{$m_{GP}$ v.s. $m_{FV}$ at $t=3/4$.}
			\label{fig:Planning:m34}
		\end{subfigure}\\
	\begin{subfigure}[b]{0.3\textwidth}
		\includegraphics[width=\textwidth]{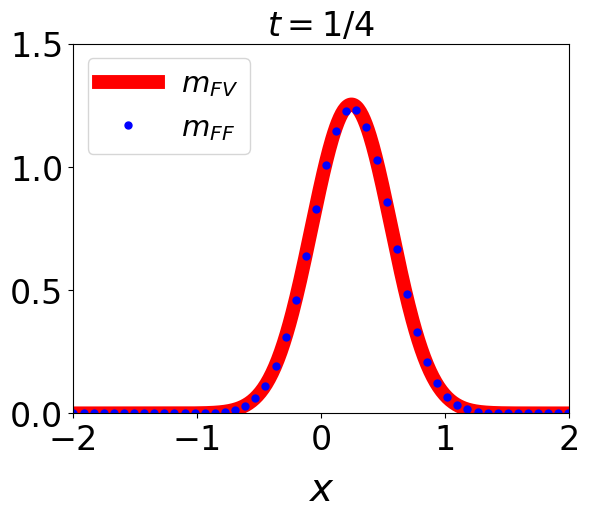}
		\caption{$m_{FF}$ v.s. $m_{FV}$ at $t=1/4$.}
		\label{fig:Planning:ffm14}
	\end{subfigure}
	\begin{subfigure}[b]{0.3\textwidth}
		\includegraphics[width=\textwidth]{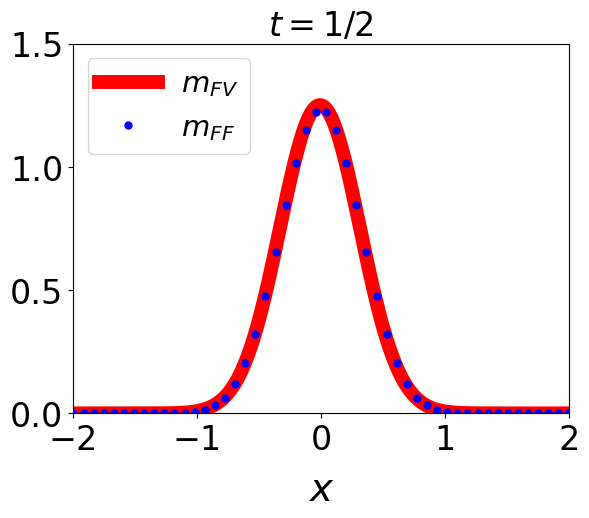}
		\caption{$m_{FF}$ v.s. $m_{FV}$ at $t=1/2$.}
		\label{fig:Planning:ffm24}
	\end{subfigure} 
	\begin{subfigure}[b]{0.3\textwidth}
		\includegraphics[width=\textwidth]{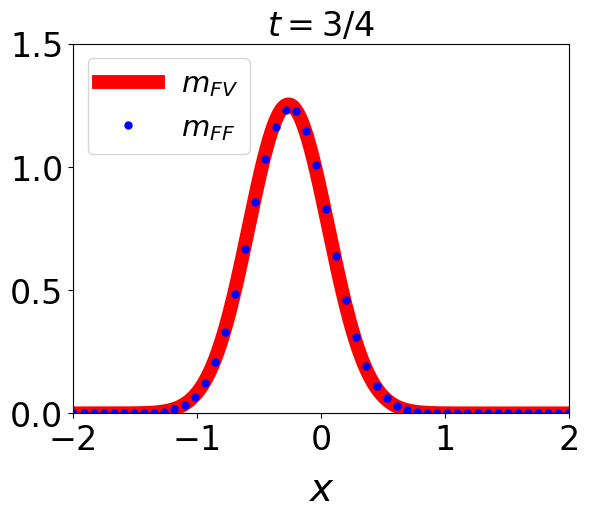}
		\caption{$m_{FF}$ v.s. $m_{FV}$ at $t=3/4$.}
		\label{fig:Planning:ffm34}
	\end{subfigure}
	\caption{Numerical results for the Planing MFG \eqref{CrowdMd3}: (a), (b) The histories of Gauss--Newton iterations; (c) uniformly sampled points in the space-time domain; (d), (e) contours of the pointwise errors of the numerical solutions; (f)-(k) time slices of the numerical solutions and the reference at $t=1/4, 1/2, 3/4$. }
	\label{fig:PlanningMFG}
\end{figure}

We consider a planning MFG, a variant of the crowd motion model given in \cite{ruthotto2020machine}. Let $\rho_G(\cdot, \mu_0, \sigma_0)$ be the probability density function of a one-dimensional Gaussian with mean $\mu_0\in \mathbb{R}$ and standard variance $\sigma_0$. We seek $(u, m)$ solving 
\begin{comment}
\begin{align}
\label{CrowdMd}
\begin{cases}
-\frac{\partial u}{\partial t} - 0.05\Delta u + \frac{8}{(1+m)^{\frac{3}{4}}}|\nabla u|^2 = \frac{1}{3200},\\
\frac{\partial m}{\partial t} - 0.05\Delta m - 16\div\bigg(\frac{m\nabla u}{(1+m)^{\frac{3}{4}}}\bigg)=0,\\
u(T, x) = 0, m(0, x)=m_0(x), \\
u(t, x) = 0, m(t, x)=0\ \text{on}\ [0, T]\times\Gamma_D,\\
\frac{\partial u}{\partial n} = 0, \frac{\partial m}{\partial n} = 0, \ \text{on}\ [0, T]\times \Gamma_N. 
\end{cases}
\end{align}
\end{comment}
\begin{comment}
\begin{align}
	\label{CrowdMd0}
	\begin{cases}
		-\frac{\partial u}{\partial t} - \nu\Delta u + \frac{1}{2}|\nabla u|^2 = \lambda_1 m + \lambda_2 Q(x),\\
		\frac{\partial m}{\partial t} - \nu\Delta m - \div(m\nabla u)=0,\\
		u(T, x) = 0, m(0, x)=m_0(x), \\
		u(t, x) = 0, m(t, x)=0\ \text{on}\ [0, T]\times\Gamma_D,\\
		\frac{\partial u}{\partial n} = 0, \frac{\partial m}{\partial n} = 0, \ \text{on}\ [0, T]\times \Gamma_N. 
	\end{cases}
\end{align}
\begin{align}
	\label{CrowdMd2}
	\begin{cases}
		-\frac{\partial u}{\partial t} + \frac{1}{2}|\nabla u|^2 = \lambda_1 (\ln m(x)+1) + \lambda_2 Q(x),\\
		\frac{\partial m}{\partial t} - \div(m\nabla u)=0,\\
		u(T, x) = 1 + \ln \rho(T, x) - \ln \rho_1(x), m(0, x)=\rho_0(x).
	\end{cases}
\end{align}
\end{comment}
\begin{align}
	\label{CrowdMd3}
	\begin{cases}
		-\frac{\partial u}{\partial t} + \frac{1}{2}|u_x|^2 = 0.01 m, \text{ in } (0, 1)\times\mathbb{R}, \\
		\frac{\partial m}{\partial t} - (mu_x)_x=0, \text{ in } (0, 1)\times\mathbb{R},\\
		m(0, x)=\rho_0(x), m(1, x)=\rho_1(x), \text{ in } \mathbb{R},
	\end{cases}
\end{align}
where $\rho_0(x)=\rho_G(x, 0.5, 0.1)$ and  $\rho_1(x)=\rho_G(x, -0.5, 0.1)$. Since the time and the space domains have different variability, following \cite{chen2021solving}, we use the anisotropic kernel
\begin{align*}
	\kappa((t,s), (t',s'))=\operatorname{exp}(-(s-s')^2/\sigma_1^2-(t-t')^2/\sigma_2^2)
\end{align*}
for the GP method, where $(\sigma_1, \sigma_2)=(1/\sqrt{5}, 1/\sqrt{2})$. We choose the regularization parameters $\eta_1=\eta_2=\mu_1=\mu_2=10^{-5}$  in Remarks \ref{rmk1dnugget} and \ref{regFF}. For the FF method, we use orthogonal random Fourier features stated in Remark \ref{rmkrff}. In \eqref{mtxW}, we choose $\varsigma=0.2$ and select $200$ random Fourier features both for approximating $u$ and $m$. We write \eqref{CrowdMd3} in an analog form of \eqref{gsmfgslo} with $Q_b=D_b=0$, $Q=4$, $D=3$, and with linear operators $L_1(u)=u$, $L_2(u)=\partial_tu$, $L_3(u)=\partial_xu$, $L_4(u)=\partial^2_x u$, $J_1(m)=m$, $J_2(m)=\partial_tm$, and $J_3(m)=\partial_xm$. To solve the optimization problems, we apply the Gauss--Newton method with step size $0.1$. Both algorithms stop after {$320$}  iterations. In the experiments, we uniformly sample $M_\Omega=1200$ points in $(0, 1)\times (-2, 2)$, $200$ samples in $\{0\}\times (-2, 2)$, and $200$ points in $\{1\}\times (-2, 2)$.  To show the accuracy of our methods, we discretize the domain $[0, 1]\times [-2, -2]$ with $64\times 512$ uniformly distributed grid points, solve \eqref{CrowdMd3} via the Eulerian solver \footnote{https://github.com/EmoryMLIP/MFGnet.jl.git} used in \cite{ruthotto2020machine}, and save the results as a reference. Denote by $(u_{GP}, m_{GP})$ and  $(u_{FF}, m_{FF})$ the numerical solutions of the GP method and the FF algorithm, separately.  We represent the result of the Eulerian solver by $(u_{FV}, m_{FV})$.
Figure \ref{fig:PlanningMFG} plots numerical results for the planning MFG.  We plot the histories of Gauss--Newton iterations in Figures \ref{fig:Planning:GPloss}-\ref{fig:Planning:FFloss}. Both the GP and the FF methods use the same set of samples, which is shown in Figure \ref{fig:Planning:samples}. Since we care more about the evolution of the probability density in the planning problem, we show the pointwise error between $m_{GP}$ and $m_{FV}$ in Figure \ref{fig:Planning:mGPerrors}, and plot the error between $m_{FF}$ and $m_{FV}$ in Figure \ref{fig:Planning:mFFerrors}. We see that the errors are larger near the peaks of the probability density. This phenomenon suggests the need to study non-uniformly distributed sample points in the future. Figures \ref{fig:Planning:m14}-\ref{fig:Planning:ffm34} compare various time slices of $m_{GP}$, $m_{FF}$, and $m_{FV}$ at time $t=1/4, 1/2, 3/4$ to highlight the accuracy of our methods.

\section{Conclusions and Future Work}
\label{secConft}
This paper presents two meshless algorithms, the GP method and the FF algorithm, to solve MFGs. The GP method adapts the  algorithm in \cite{chen2021solving} to solve MFGs and finds numerical solutions in RKHSs. The convergence analysis of the GP method does not rely on the Lasry--Lions monotonicity condition. Hence, we plan to apply the GP method to solve MFGs with other monotonicity conditions or with non-monotone couplings in future work. To get better performance, we introduce the FF algorithm seeking approximations in Fourier features spaces. Compared to the GP method, the FF algorithm consumes less precomputation time without losing  accuracy by using a prescribed number of base functions. Since the Fourier features space lacks compactness, we cannot use the same convergence arguments of the GP method. Instead, the Lasry--Lions monotonicity ensures the boundedness of numerical errors and gives the convergence of the FF method.  We plan to study the convergence of the FF algorithm in displacement monotone or non-monotone settings in future work. We believe that one can also use the FF method to solve general PDEs. 
As we have observed in the numerical experiments, the choices of base functions and parameters in our methods profoundly influence the accuracy of numerical solutions. Hence, in future work, we also plan to investigate methods of hyperparameter learning and different sampling techniques.

\section*{Acknowledgment}
We want to thank Yifan Chen for sharing the code of \cite{chen2021solving} online and discussing the implementation. We also thank Professor Lars Ruthotto for sharing the Matlab code about the crowd motion model in \cite{ruthotto2020machine} {and thank Mathieu Lauri\`{e}re for  discussions about the finite difference methods of MFGs}. C. Mou gratefully acknowledges the support by CityU Start-up Grant 7200684 {and Hong Kong RGC Grant ECS 9048215}. C. Zhou is supported by Singapore MOE (Ministry of Education), AcRF Grants R-146-000-271-112 and R-146-000-284-114, and NSFC Grant No. 11871364.

\newpage
\bibliographystyle{plain}
\IfFileExists{"mfgpaper.bib"}{\bibliography{mfgpaper}}

\newpage
\appendix
\section{Proofs of the results}
\label{appendix}
\begin{proof}[Proof of Theorem \ref{convUpBdFF1}]
	Let $(u^*, m^*, \overline{H}^*)$ be the solution of \eqref{1dsmfg}. By the explicit formulas in \eqref{exp1du}, $u^*$ and $m^*$ are smooth. Thus, according to the approximation theorem of Fourier series \cite[Chapter 14]{Vasy2015PartialDE}, for any $\epsilon>0$, there exist $N\in\mathbb{N}$ and functions $u^N, m^N\in\mathcal{G}^N$ such that
	\begin{align}
	\label{ffbd1}
	\max_{i\leq 2}\sup_{x\in \mathbb{T}}|\dif^{(i)} u^N(x)-\dif^{(i)} u^*(x)|\leq \epsilon
	\end{align}
and 
\begin{align}
	\label{ffbd2}
	\max_{i\leq 1}\sup_{x\in \mathbb{T}}|\dif^{(i)} m^N(x)-\dif^{(i)} m^*(x)|\leq \epsilon, 
\end{align}
where $\dif^{(i)}$ represents the $i$th order derivative. 
Hence, by the definition of $\mathcal{Q}$ in  \eqref{defmq}, we have
\begin{align}
\label{Qbd}
\begin{split}
\mathcal{Q}(u^N, m^N, \overline{H}^*) =& \int_{\mathbb{T}}|e^{H(x, u_x^N)-\overline{H}^*}-m^N|^2\dif x+\biggl |\int_{\mathbb{T}}u^N\dif x\biggl|^2+\biggl|\int_{\mathbb{T}}m^N\dif x - 1\biggl|^2\\
&+\int_{\mathbb{T}}|m^N(u_{xx}^N+b_x(x))+m_x^N(u_x^N+b(x))|^2\dif x \\
\leq & \int_{\mathbb{T}}|e^{H(x, u_x^N)-\overline{H}^*}- e^{H(x, u^*_x)-\overline{H}^*}|^2\dif x+\int_{\mathbb{T}}|m^N- m^*|^2\dif x\\
&+\biggl|\int_{\mathbb{T}}u^N\dif x - \int_{\mathbb{T}}u^*\dif x\biggl|^2+\biggl|\int_{\mathbb{T}}m^N\dif x -  \int_{\mathbb{T}}m^*\dif x\biggl|^2 \\
&+\int_{\mathbb{T}}|m^N(u_{xx}^N+b_x(x))-m^*(u_{xx}^*+b_x(x))|^2\dif x\\
&+\int_{\mathbb{T}}|m_x^N(u_x^N+b(x))-m^*_x(u_x^*+b(x))|^2\dif x.
\end{split}
\end{align}
Next, we estimate the terms at the right-hand side (RHS) of \eqref{Qbd}. From \eqref{ffbd1} and \eqref{ffbd2}, we get
\begin{align}
	\label{convFFQbd1}
	\int_{\mathbb{T}}|m^N-m^*|^2\dif x+\biggl|\int_{\mathbb{T}}u^N\dif x - \int_{\mathbb{T}}u^*\dif x\biggl|^2+\biggl|\int_{\mathbb{T}}m^N\dif x -  \int_{\mathbb{T}}m^*\dif x\biggl|^2 \leq 3\epsilon^2. 
\end{align}
For the first term at the RHS of \eqref{Qbd}, we have
\begin{align}
	\label{Qbd1stt}
	\begin{split}
		&\int_{\mathbb{T}}|e^{H(x, u_x^N)-\overline{H}^*}- e^{H(x, u^*_x)-\overline{H}^*}|^2\dif x\\
		=&e^{-2\overline{H}^*}\int_{\mathbb{T}}\bigg|\int_0^1\dif e^{H(x, \hat{u}^N_x)}\bigg|^2\dif x\ (\text{where}\ \hat{u}_x:=su_x^N + (1-s)u_x^*, s\in [0, 1])\\
		=&e^{-2\overline{H}^*}\int_{\mathbb{T}}\bigg|\int_0^1e^{H(x,\hat{u}_x^N)}D_pH(x,\hat{u}^N_x)(u_x^N-u_x^*)\dif s\bigg|^2\dif x.
	\end{split}
\end{align}
Using \eqref{ffbd1}, we get
\begin{align}
	\label{bduhat}
	\sup_{x\in \mathbb{T}}|\hat{u}^N_x(x)|\leq \sup_{x\in \mathbb{T}}s|u_x^N-u_x^*|+\sup_{x\in\mathbb{T}}|u_x^*|\leq s\epsilon + \sup_{x\in \mathbb{T}}|u_x^*|.
\end{align}
Thus, by the definition of $H$ in \eqref{defH}, \eqref{ffbd1}, \eqref{Qbd1stt}, and \eqref{bduhat},  there exists a constant $C>0$ such that 
\begin{align}
	\label{Qbd1stt2}
	\int_{\mathbb{T}}|e^{H(x, u_x^N)-\overline{H}^*}- e^{H(x, u^*_x)-\overline{H}^*}|^2\dif x \leq C\epsilon^2. 
\end{align}
Next, we use \eqref{ffbd1} and \eqref{ffbd2}, and get
\begin{align}
	\label{Qbd5tt}
	\begin{split}
		&\int_{\mathbb{T}}|m^N(u_{xx}^N+b_x(x))-m^*(u_{xx}^*+b_x(x))|^2\dif x\\
		\leq & \int_{\mathbb{T}} |(m^N-m^*)(u_{xx}^N+b_x(x))|^2\dif x + \int_{\mathbb{T}}|m^*(u_{xx}^N-u_{xx}^*)|^2\dif x \leq C\epsilon^2. 
	\end{split}
\end{align}
Similarly, there exists a constant $C>0$ such that
\begin{align}
	\label{bd6tht}
	\int_{\mathbb{T}}|m_x^N(u_x+b(x))-m^*_x(u_x^*+b(x))|^2\dif x\leq C\epsilon^2. 
\end{align}
Combining \eqref{convFFQbd1}, \eqref{Qbd1stt2}, \eqref{Qbd5tt}, \eqref{bd6tht}, we get from \eqref{Qbd} that
\begin{align}
	\label{Qbdf}
	\mathcal{Q}(u^N, m^N, \overline{H}^*)\leq C\epsilon^2. 
\end{align}
Using \eqref{ffbd1} and \eqref{ffbd2} again, we get
\begin{align}
	\label{bd2n}
	\int_{\mathbb{T}}\|u^N\|^2\dif x\leq 2\|u^*\|_{L^2(\mathbb{T})}^2+2\epsilon^2\ \text{and} \ \int_{\mathbb{T}}\|m^N\|^2\dif x\leq 2\|m^*\|_{L^2(\mathbb{T})}^2+2\epsilon^2.
\end{align}
Therefore, we conclude \eqref{bdJ} by combining \eqref{Qbdf}, \eqref{bd2n}, and the definition of $J_{\gamma}$ in \eqref{defca}. 
\end{proof}

\begin{proof}[Proof of Corollary \ref{convFF1Q}] Let $(u^*, m^*, \overline{H}^*)$ be the solution of \eqref{1dsmfg}.
By Theorem \ref{convUpBdFF1}, for given small $\epsilon>0$ and $\gamma>0$, there exists a constant $C$ and functions $(u^N, m^N)\in \mathcal{G}^N\times\mathcal{G}^N$ such that
\begin{align}
	\label{bdJst}
	J_\gamma(u^N, m^N, \overline{H}^*) \leq 2\|u^*\|_{L^2}^2+2\|m^*\|_{L^2}^2+|\overline{H}^*|^2+C(1+\gamma)\epsilon^2.
\end{align}
Let $({u}^{N, \gamma}, {m}^{N, \gamma}, {\overline{H}}^{N, \gamma})$ be a minimizer in \eqref{1dcts} given $N$ and $\gamma$. Then, we have
\begin{align*}
J_\gamma({u}^{N, \gamma}, {m}^{N, \gamma}, {\overline{H}}^{N, \gamma})\leq J_\gamma(u^N, m^N, \overline{H}^*). 
\end{align*}
Thus, using \eqref{bdJst} and the definition of $J_\gamma$ in \eqref{defca}, we have
\begin{align*}
\gamma\mathcal{Q}({u}^{N, \gamma}, {m}^{N, \gamma}, {\overline{H}}^{N, \gamma})\leq 2\|u^*\|_{L^2}^2+2\|m^*\|_{L^2}^2+|\overline{H}^*|^2+C(1+\gamma)\epsilon^2,
\end{align*}
which yields
\begin{align*}
\mathcal{Q}({u}^{N, \gamma}, {m}^{N, \gamma}, {\overline{H}}^{N, \gamma})\leq \frac{2}{\gamma}(\|u^*\|_{L^2}^2+\|m^*\|_{L^2}^2+|\overline{H}^*|^2)+C(\frac{1}{\gamma}+1)\epsilon^2. 
\end{align*}
Therefore, we conclude \eqref{bdQ} by passing $\gamma$ to infinity. 
\end{proof}

\begin{proof}[Proof of Theorem \ref{conv1dFF}]
By Theorem \ref{convUpBdFF1} and Corollary \ref{convFF1Q}, there exists a sequence $\{(N_i, \gamma_i)\}_{i=1}^{\infty}$ such that the set $\{({u}^{i}, {m}^{i}, {\overline{H}}^{i})\}_{i=1}^\infty$, where $({u}^{i}, {m}^{i}, {\overline{H}}^{i})$ is a minimizer of \eqref{1dcts} given $N_i$ and $\gamma_i$, satisfies
\begin{align}
\label{Qconv}
\mathcal{{Q}}({u}^{i}, {m}^{i}, {\overline{H}}^{i})\rightarrow 0\ \text{as}\ i\rightarrow \infty. 
\end{align}
Meanwhile, there exits a constant $C$ such that for $i\geq 1$, 
\begin{align}
\label{bdJgamma}
J_{\gamma_i}(u^i, m^i, \overline{H}^i)\leq C.
\end{align}
Thus, by \eqref{bdJgamma} and the definition of $J_{\gamma_i}$ in \eqref{defca}, there exists a constant $C$ such that
\begin{align}
	\label{bdumld}
	\|{u}^i\|_{L^2(\mathbb{T})}\leq C,\ \|{m}^i\|_{L^2(\mathbb{T})}\leq C, \ \text{and}\ |{\overline{H}}^i|\leq C. 
\end{align}
Then, there exists $\overline{H}^\infty\in\mathbb{R}$ such that, up to a subsequence,  $\overline{H}^i$ converges to $\overline{H}^\infty$ in $\mathbb{R}$.
Let $h^i$ and $g^i$ be functions such that 
\begin{align}
\label{tpde}
\begin{cases}
	\frac{{(u^i)}_x^2}{2}+V(x)+b(x)u^i_x =  \ln({m}^i + h^i(x))+\overline{H}^i,\\
	-({m}^i({u}^i_x+b(x)))_x = g^i(x).
\end{cases}
\end{align}
By \eqref{Qconv}, we have
\begin{align}
\label{convhg}
\lim_{i\rightarrow \infty}\int_{\mathbb{T}}m^i\dif x=1, 
\lim_{i\rightarrow \infty}\|h^i\|_{L^2(\mathbb{T})}\rightarrow
 0 \text{ and } \lim_{i\rightarrow \infty}\|g^i\|_{L^2(\mathbb{T})}\rightarrow
 0. 
\end{align} 
Hence, integrating the first equation of \eqref{tpde},  we get from \eqref{convhg} and \eqref{bdumld} that
\begin{align*}
\|{u}^i_x\|_{L^2(\mathbb{T})}\leq C. 
\end{align*}
Thus, ${u}^i\in H^{1}(\mathbb{T})$. Since ${u}^i$ is the finite linear combination of trigonometric functions, there exists a constant $C$ such that
\begin{align}
\label{bdhatucc1}
\|{u}^i\|_{C^1(\mathbb{T})}\leq C\|{u}^i\|_{H^{1}(\mathbb{T})}\leq C. 
\end{align}

Next, we prove the convergence of $(u^i, m^i, \overline{H}^i)$ using the  Lasry-Lions monotonicity argument. 
Let $(u^*, m^*, \overline{H}^*)$ be the solution to \eqref{1dsmfg}. Then, from \eqref{tpde}, we get
\begin{align}
\label{diffpdes}
\begin{cases}
\frac{(u^i_x)^2}{2}-\frac{(u^*_x)^2}{2}+b(x)u^i_x-b(x)u^*_x =  \ln({m}^i + h^i(x)) -  \ln m^*+\overline{H}^i-\overline{H}^*,\\
-(m^i(u^i_x+b(x)))_x + (m^*(u^*_x+b(x)))_x = g^i(x).
\end{cases}
\end{align}
Multiplying the first equation of \eqref{diffpdes} by $m^i+h^i-m^*$ and the second equation in \eqref{diffpdes} by $u^i-u^*$, subtracting the resulting equations, and integrating by parts, we obtain 
\begin{align}
\label{1dmonoteq}
\begin{split}
&\frac{1}{2}\int_{\mathbb{T}}(m^i+h^i)|u_x^*-u_x^i|^2\dif x+\frac{1}{2}\int_{\mathbb{T}}m^*|u_x^*-u_x^i|^2\dif x\\
&+\int_{\mathbb{T}}(\ln(m^i+h^i)-\ln m^*)(m^i+h^i-m^*)\dif x\\
=& (\overline{H}^*-\overline{H}^i)\int_{\mathbb{T}}(m^i+h^i-m^*)\dif x+\int_{\mathbb{T}}g^i(u^i-u^*)\dif x \\
&+\int_{\mathbb{T}}h^i(u_x^i+b(x))(u_x^i-u_x^*)\dif x.
\end{split}
\end{align}
By passing $i\rightarrow 0$,  \eqref{convhg}, \eqref{bdhatucc1},  and  \eqref{1dmonoteq} yield
\begin{align}
\label{lim1dky}
\begin{split}
\lim\limits_{i\rightarrow 0}&\frac{1}{2}\int_{\mathbb{T}}(m^i+h^i)|u_x^*-u_x^i|^2\dif x+\frac{1}{2}\int_{\mathbb{T}}m^*|u_x^*-u_x^i|^2\dif x\\
&+\int_{\mathbb{T}}(\ln(m^i+h^i)-\ln m^*)(m^i+h^i-m^*)\dif x=0. 
\end{split}
\end{align}
From the first equation of \eqref{tpde}, we get
\begin{align*}
{m}^i+h^i(x) = e^{\frac{1}{2}({u}^i_x+b(x))^2+V(x)-{\overline{H}^i}-b^2(x)}\geq e^{-C}\min_{x\in \mathbb{T}}e^{V(x)-b^2(x)}.
\end{align*}
\begin{comment}
Next, we use the formula $u^*_x(x)=b(x)$ in \eqref{exp1du} and multiply the second equation in \eqref{tpde} by $\hat{u}_x-u_*$. Using integration by parts, we get
\begin{align*}
\int_{\mathbb{T}}\hat{m}|\hat{u}_x-u^*_x|^2\dif x = \int_{\mathbb{T}}g^i(x)(\hat{u}-u^*)\dif x,
\end{align*}
which is equivalent to 
\begin{align}
\label{bduhatmu}
\int_{\mathbb{T}}(\hat{m}+h^i(x))|\hat{u}_x-u^*_x|^2\dif x = \int_{\mathbb{T}}g^i(x)(\hat{u}_x-u^*_x)\dif x + \int_{\mathbb{T}}h^i(x)|\hat{u}_x-u^*_x)|^2\dif x. 
\end{align}
Thus, we have
\begin{align*}
\int_{\mathbb{T}}(\hat{m}+h^i(x))|\hat{u}_x-u^*_x|^2\dif x =& \int_{\mathbb{T}}g^i(x)(\hat{u}_x-u^*_x)\dif x + \int_{\mathbb{T}}h^i(x)|\hat{u}_x-u^*_x)|^2\dif x\\
=& \|g^i\|^2_{L^2}\|\hat{u}_x\|^2+\|g^i\|^2_{L^2}\|u^*_x\|_{L^2}\\
&+(\|\hat{u}\|^2_{C^1(\mathbb{T})}+\|u^*\|_{C^1(\mathbb{T})}^2)\|h^i(x)\|_{L^2}.
\end{align*}
By \eqref{bdmhi}, there exist constant $C>0$ such that
\begin{align}
\label{ubuxuh}
\begin{split}
\int_{\mathbb{T}}|\hat{u}_x-u^*_x|^2\dif x \leq & \int_{\mathbb{T}}g^i(x)(\hat{u}_x-u^*_x)\dif x + \int_{\mathbb{T}}h^i(x)|\hat{u}_x-u^*_x)|^2\dif x\\
=& \|g^i\|^2_{L^2}\|\hat{u}_x\|^2+\|g^i\|^2_{L^2}\|u^*_x\|_{L^2}\\
&+(\|\hat{u}\|^2_{C^1(\mathbb{T})}+\|u^*\|_{C^1(\mathbb{T})}^2)\|h^i(x)\|_{L^2}. 
\end{split}
\end{align}
\end{comment}
Then, the left-hand side of \eqref{lim1dky} is non-negative, and $m^i+h^i$ is uniformly bounded below. Hence,  \eqref{lim1dky} yields 
\begin{align*}
 u^i_x\rightarrow  u^*_x\ \text{in}\ L^2(\mathbb{T}). 
\end{align*}
Thus, up to a sub-sequence, $u^i_x$ converges to $u^*_x$ pointwisely. By the first equation of \eqref{tpde}, there exists a function $m^\infty$ such that, up to a sub-sequence, 
\begin{align*}
{m}^i\rightarrow m^\infty= e^{\frac{(u_x^*)^2}{2}+V(x)+b(x)u_x^*-{\overline{H}}^\infty} \text{ pointwisely}. 
\end{align*}
Since ${m}^i$ is uniformly bounded above, by the dominated convergence theorem, we conclude that ${m}^i$ converges to ${m}^\infty$ in $L^1(\mathbb{T})$. We note that $(u^*, {m}^\infty, {\overline{H}}^\infty)$ satisfies \eqref{1dsmfg}. By the uniqueness of the solution, we conclude that ${m}^\infty=m^*$ and ${\overline{H}}^\infty=\overline{H}^*$. 
\end{proof}
\begin{proof}[Proof of Theorem \ref{convUpBdFFnd}]\label{proof2dot7}
The arguments here are similar to the proof of Theorem \ref{convUpBdFF1}. Using the approximation theorem of Fourier series and \eqref{lpF}, we perform similar computations as in \eqref{Qbd}-\eqref{bd2n}. Then, we conclude \eqref{bdJnd} by the definition of $J_\gamma$.  
\end{proof}
\begin{proof}[Proof of Corollary \ref{convFFQnd}]
The proof is the same as in the arguments of Corollary \ref{convFF1Q}. 
\end{proof}
\begin{corollary}
\label{bdGN}
Let $\{f^i\}_{i=1}^{\infty}$ be a sequence such that $f^i\in \mathcal{G}^i$ defined in \eqref{defGnd}. If $\{f^i\}_{i=1}^{\infty}$ is uniformly bounded in $L^2(\mathbb{T}^d)$. Then, there exists a constant $C>0$ such that
\begin{align*}
	\|f^i\|_{L^\infty(\mathbb{T}^d)}\leq C. 
\end{align*}
Moreover, there exits a continuous function $f$ such that, up to a sub-sequence, $f^i$ converges to $f$ in $L^\infty(\mathbb{T}^d)$ as $i\rightarrow\infty$.  
\end{corollary}
\begin{proof}
Since $f^i\in \mathcal{G}^i$ for $i\geq 1$, there exist real numbers  $\{\alpha^i_j\}_{j\in \mathbb{Z}_i^d}$, $\{\beta^i_j\}_{j\in \mathbb{Z}_i^d}$, and $c^i$ such that
\begin{align*}
f^i(x) = c^i+\sum_{j\in \mathbb{Z}_i^d}\alpha_j^i\sin(2\pi j^T x)+\sum_{j\in \mathbb{Z}_i^d}\beta_{j}^i\cos(2\pi j^T x), x\in \mathbb{T}^d. 
\end{align*}
Thus, we have
\begin{align}
	\label{refi2}
\begin{split}
	(f^i(x))^2 =& (c^i)^2+\sum_{k, j\in \mathbb{Z}_i^d}\alpha_j^i\alpha_k^i\sin(2\pi j^T x)\sin(2\pi k^T x)\\
	&+\sum_{k, j\in \mathbb{Z}_i^d}\beta_{j}^i\beta_k^i\cos(2\pi j^T x)\cos(2\pi k^T x)\\
	&+2\sum_{k, j\in \mathbb{Z}_i^d}\beta_{j}^i\alpha_k^i\cos(2\pi j^T x)\sin(2\pi k^T x) + c^i\sum_{ j\in \mathbb{Z}_i^d}\alpha^i_j\sin(2\pi j^T x)\\
	&+c^i\sum_{ j\in \mathbb{Z}_i^d}\beta_j^i\cos(2\pi j^T x).
\end{split}
\end{align}
Since $\forall i, j\in \mathbb{Z}^d$, 
\begin{align*}
&\int_{\mathbb{T}^d}\cos(2\pi i^T x)\dif x = \int_{\mathbb{T}^d}\sin(2\pi i^T x) \dif x = 0,\\
&\int_{\mathbb{T}^d}\sin(2\pi j^T x)\sin(2\pi i^T x) \dif x = \begin{cases}
0, \ \text{if}\ i\not=j\ \text{or}\ i=j=0, \\
\frac{1}{2}, \ \text{if}\ i = j\not=0,
\end{cases}\\
&\int_{\mathbb{T}^d}\cos(2\pi j^T x)\cos(2\pi i^T x) \dif x = \begin{cases}
	0, \ \text{if}\ i\not=j, \\
	\frac{1}{2}, \ \text{if}\ i = j\not=0,\\
	1,\ \text{if}\ i = j = 0,
\end{cases}\\
\end{align*}
and 
\begin{align*}
\int_{\mathbb{T}^d}\cos(2\pi j^T x)\sin(2\pi i^T x) \dif x = 0, 
\end{align*}
we get from \eqref{refi2} that
\begin{align}
\label{refi3}
\begin{split}
	\int_{\mathbb{T}^d}(f^i)^2\dif x = (c^i)^2 + \frac{1}{2}\sum_{ j\in \mathbb{Z}_i^d}(\alpha_j^i)^2+\frac{1}{2}\sum_{j\in \mathbb{Z}_i^d}(\beta_{j}^i)^2, 
\end{split}
\end{align}
By \eqref{refi3} and the fact that  $\{f^i\}_{i=1}^\infty$ is uniformly bounded in $L^2(\mathbb{T}^d)$, there exits a constant $C>0$ such that
\begin{align}
\label{bdl2}
\begin{split}
\int_{\mathbb{T}^d}(f^i)^2\dif x = (c^i)^2 + \frac{1}{2}\sum_{ j\in \mathbb{Z}_i^d}(\alpha_j^i)^2+\frac{1}{2}\sum_{j\in \mathbb{Z}_i^d}(\beta_{j}^i)^2 \leq C. 
\end{split}
\end{align}
Thus, using \eqref{bdl2} and Young's inequality, we obtain
\begin{align*}
\|f^i\|_{L^\infty(\mathbb{T}^d)}^2\leq (|c^i| + \sum_{ j\in \mathbb{Z}_i^d}|\alpha^i_j| + \sum_{ j\in \mathbb{Z}_i^d}|\beta^i_j|)^2
\leq  C\int_{\mathbb{T}^d}(f^i)^2\dif x\leq C. 
\end{align*}
Hence, $\{f^i\}_{i=1}^\infty$ is uniformly bounded in $L^\infty(\mathbb{T}^d)$. 

Meanwhile, by \eqref{bdl2}, there exist  real numbers $c$, $\{\alpha_j\}_{j=1}^\infty$, and $\{\beta_j\}_{j=1}^\infty$ such that, up to subsequences, $c^i\rightarrow c$, $\alpha_j^i\rightarrow \alpha_j$, and $\beta_j^i\rightarrow \beta_j$, for all $j\in \mathbb{Z}^d_i$, as $i\rightarrow\infty$. Furthermore, for any $\epsilon>0$, there exists $N\in \mathbb{N}$ such that for any $|i|\geq N$, 
\begin{align}
\label{bdalphabeta}
\sum_{ j\in \mathbb{Z}^d_i, |j|\geq N}|\alpha_j^i| \leq \epsilon,  \sum_{ j\in \mathbb{Z}^d_i, |j|\geq N}|\beta_j^i| \leq \epsilon, 
\end{align}
\begin{align}
	\label{bdab}
	\sum_{ j\in \mathbb{Z}^d, |j|\geq N}|\alpha_j| \leq \epsilon, \ \text{and}\ \sum_{ j\in \mathbb{Z}^d, |j|\geq N}|\beta_j| \leq \epsilon. 
\end{align}
We define
\begin{align*}
f(x) = c + \sum_{ j\in \mathbb{Z}^d\backslash \{0\}}\alpha_j\sin(2\pi j^T x) + \sum_{ j\in \mathbb{Z}^d\backslash \{0\}}\beta_j\cos(2\pi j^T x) 
\end{align*}
and 
\begin{align*}
	f_N(x) = c + \sum_{ j\in \mathbb{Z}^d_N}\alpha_j\sin(2\pi j^T x) + \sum_{ j\in \mathbb{Z}^d_N}\beta_j\cos(2\pi j^T x). 
\end{align*}
Thus, for $|i|>N$, by \eqref{bdalphabeta} and \eqref{bdab}, there exists a constant $C>0$ such that 
\begin{align*}
\begin{split}
\|f^i-f\|_{L^\infty(\mathbb{T}^d)}\leq& \|f^i-f_N\|_{L^\infty(\mathbb{T}^d)} + \|f_N-f\|_{L^\infty}
\\ \leq  &|c^i-c|+\sum_{j\in \mathbb{Z}_N^d}|\alpha_j^i-\alpha_j| + \sum_{ j\in \mathbb{Z}_N^d}|\beta_j^i-\beta_j|\\
&+ \sum_{ j\in \mathbb{Z}^d_i, |j|>N}|\alpha_j^i| + \sum_{ j\in \mathbb{Z}^d_i, |j|>N}|\beta_j^i|\\
&+ \sum_{ j\in \mathbb{Z}^d, |j|>N}|\alpha_j| + \sum_{ j\in \mathbb{Z}^d, |j|>N}|\beta_j| \leq C\epsilon. 
\end{split}
\end{align*}
Therefore, we conclude that $f^i$ converges to $f$ in $L^\infty(\mathbb{T}^d)$. 
\end{proof}
\begin{proof}[Proof of Theorem \ref{convndFF}]
By Theorem \ref{convUpBdFFnd} and Corollary \ref{convFFQnd}, there exist a sequence $\{(N^i, \gamma^i)\}_{i=1}^\infty$ and a constant $C$ such that the sequence $\{(u^i, m^i, \overline{H}^i)\}_{i=1}^\infty$, where $(u^i, m^i, \overline{H}^i)$ is a minimizer of \eqref{ndcts} given $N^i$ and $\gamma^i$, satisfies 
\begin{align}
	\label{bdumnd}
	\|{u}^i\|_{L^2(\mathbb{T}^d)}\leq C,\ \|{m}^i\|_{L^2(\mathbb{T}^d)}\leq C, \ \text{and}\ |{\overline{H}}^i|\leq C. 
\end{align}
Meanwhile, let $h^i$ and $g^i$ be functions such that 
\begin{align}
	\label{tpdeddim}
	\begin{cases}
		-\Delta u^i + \frac{|\nabla u^i|^2}{2} + V(x) - F[m^i] - \overline{H}^i =  h^i(x),\\
		-\Delta m^i - \div(m^i\nabla u^i) = g^i(x).
	\end{cases}
\end{align}
Then, 
\begin{align}
	\label{convhgnd}
	\lim_{i\rightarrow \infty}\|h^i\|_{L^2(\mathbb{T}^d)}\rightarrow
	0 \text{ and } \lim_{i\rightarrow \infty}\|g^i\|_{L^2(\mathbb{T}^d)}\rightarrow
	0. 
\end{align}
Next, we study the regularity and the convergence of $\{u^i\}_{i=1}^\infty$. We split our arguments into three claims and prove  each claim. 
\begin{claim}
\label{cl1mc}
There exists a constant $C>0$ such that $\|u^i\|_{C^1(\mathbb{T}^d)}\leq C$ for all $i$  and exists a multi-valued function $\varsigma\in (L^\infty(\mathbb{T}^d))^d$ such that, up to a subsequence,  $\nabla u^i$ converges to $\varsigma$ in $(L^\infty(\mathbb{T}^d))^d$. 
\end{claim}
Integrating the first equation in  \eqref{tpdeddim}, we get
\begin{align}
\label{bdndDui}
\begin{split}
\int_{\mathbb{T}^d}\frac{|\nabla u^i|^2}{2}\dif x =-\int_{\mathbb{T}^d}V(x)\dif x + \int_{\mathbb{T}^d}F[m^i]\dif x + \overline{H}^i+ \int_{\mathbb{T}^d}h^i(x)\dif x. 	
\end{split}
\end{align}
Thus, for $i$ large enough, using Assumption \ref{hypF},  \eqref{bdumnd}, \eqref{convhgnd}, and the smoothness of $V$, we get that  $\nabla u^i$ is uniformly bounded in $(L^2(\mathbb{T}^d))^d$. Hence, by Corollary \ref{bdGN} and \eqref{bdumnd}, there exists a constant $C$ such that $\|u^i\|_{C^1(\mathbb{T}^d)}\leq C$ and exists a function $\varsigma\in (L^\infty(\mathbb{T}^d))^d$ such that, up to a sub-sequence, $\nabla u^i$ converges in $(L^\infty(\mathbb{T}^d))^d$ to  $\varsigma$.

The following claim gives the convergence of $\{m^i\}_{i=1}^\infty$ and the properties of the limit. 
\begin{claim} 
	\label{cl2mc}
	There exits a continuous function $m^\infty$ such that, up to a subsequence,  $m^i$ converges to $m^\infty$ in $H^{1}(\mathbb{T}^d)$ and in $L^\infty(\mathbb{T}^d)$ as $i\rightarrow \infty$. Moreover, $\int_{\mathbb{T}^d}m^\infty\dif x=1$ and  $m^\infty\geq 0$. 
\end{claim}
By \eqref{bdumnd} and  Corollary \ref{bdGN}, there exists a continuous function $m^\infty\in L^\infty(\mathbb{T}^d)$ such that $
m^i\rightarrow m^\infty \text{ in } L^\infty(\mathbb{T}^d)$.  
By Corollary \ref{convFFQnd}, $\lim_{i\rightarrow\infty}\int_{\mathbb{T}^d}m^i\dif x= 1$. Thus, $\int_{\mathbb{T}^d}m^\infty\dif x = 1$. Next, we show that $m^\infty\geq 0$. 

Multiplying the second equation of \eqref{tpdeddim} by $m^i$ and integrating by parts, we get
\begin{align}
\label{bdmi}
\begin{split}
\int_{\mathbb{T}^d}|\nabla m^i|^2\dif x =& -\int_{\mathbb{T}^d} m^i\langle \nabla u^i, \nabla m^i\rangle\dif x + \int_{\mathbb{T}^d}g^im^i\dif x\\
\leq & \frac{1}{2}\int_{\mathbb{T}^d}(m^i)^2|\nabla u^i|^2\dif x + \frac{1}{2}\int_{\mathbb{T}^d}|\nabla m^i|^2\dif x + \int_{\mathbb{T}^d}g^im^i\dif x, 
\end{split}
\end{align}
where we use Young's inequality in the above inequality. Hence, By \eqref{bdumnd}, \eqref{convhgnd},  Claim \ref{cl1mc}, and \eqref{bdmi}, there exists a constant $C$ such 
\begin{align*}
\|\nabla m^i\|_{L^2(\mathbb{T}^d)}\leq C. 
\end{align*}
Then, by Corollary \ref{bdGN}, up to a subsequence, $\nabla m^i$ converges to a vector field $\varXi\in (L^\infty(\mathbb{T}^d))^d$. Since $m^i\rightarrow m^\infty$ in $L^\infty(\mathbb{T}^d)$, $m^\infty$ is differentiable and $\varXi=\nabla m^\infty$. Hence, $m^i$ converges to $m^\infty$ in $H^{1}(\mathbb{T}^d)$. We multiply the second equation of \eqref{tpdeddim} by $\varphi\in C^\infty(\mathbb{T}^d)$,  integrate it by parts, and get
\begin{align*}
\int_{\mathbb{T}^d}\langle \nabla m^i, \nabla \varphi\rangle \dif x + \int_{\mathbb{T}^d} m^i\langle \nabla u^i, \nabla \varphi \rangle \dif x = \int_{\mathbb{T}^d} g^i\varphi \dif x.
\end{align*}
Passing $i$ to infinity, we obtain
\begin{align*}
\int_{\mathbb{T}^d}\langle \nabla m^\infty, \nabla \varphi\rangle \dif x + \int_{\mathbb{T}^d} m^\infty\langle \varsigma, \nabla \varphi \rangle \dif x = 0.
\end{align*}  
Hence, $m^\infty$ is a weak solution to 
\begin{align*}
-\Delta m^\infty - \div(\varsigma m^\infty) = 0 \text{ in } \mathbb{T}^d \text{ and } \int_{\mathbb{T}^d}m^\infty\dif x = 1. 
\end{align*}
Thus, by the ergodic theory, $m^\infty\geq 0$ (see Section 1.4 of \cite{bensoussan1987singular}).

Finally, we have the convergence of a subsequence of $(u^i, m^i, \overline{H}^i)_{i=1}^\infty$ to the solution of \eqref{ConGenMFG}. The proof is based on the  Lasry-Lions monotonicity argument. 
\begin{claim}
	Up to a subsequence, $u^i\rightarrow u^*$ in $H^{1}(\mathbb{T}^d)$, $m^i\rightarrow m^*$ in $H^{1}(\mathbb{T}^d)$, and $\overline{H}^i$ converges to $\overline{H}^*$ in $\mathbb{R}$ as $i$ goes to infinity. 
\end{claim}
Denote $w^i=u^i-u^*$, $\rho^i = m^i-m^*$, $\tau^i=\overline{H}^i-\overline{H}^*$. Then, we have
\begin{align}
\label{diffndeq}
\begin{cases}
-\Delta w^i + \langle \frac{\nabla u^i+\nabla u^*}{2}, \nabla w^i\rangle=F[m^i]-F[m^*]+\tau^i+h^i(x),\\
-\Delta \rho^i - \div(\nabla w^im^i)-\div(\rho^i\nabla u^*)=g^i(x).
\end{cases}
\end{align}
Multiplying the first equation in \eqref{diffndeq} by $\rho^i$ and integrating by parts, we obtain
\begin{align}
\label{wirhoi}
\begin{split}
&\int_{\mathbb{T}^d}\langle \nabla w^i, \nabla \rho^i\rangle\dif x + \int_{\mathbb{T}^d}\bigg\langle\frac{\nabla u^i+\nabla u^*}{2}, \nabla w^i\rho^i\bigg\rangle\dif x \\
=&\int_{\mathbb{T}^d}(F[m^i]-F[m^*])\rho^i\dif x + \tau^i\int_{\mathbb{T}^d}\rho^i\dif x + \int_{\mathbb{T}^d}h^i(x)\rho^i\dif x. 
\end{split}
\end{align}
Then, we integrate the second equation of \eqref{diffndeq} multiplied with $w^i$ over $\mathbb{T}^d$ and get
\begin{align}
\label{rhoiwi}
\begin{split}
\int_{\mathbb{T}^d}\langle \nabla \rho^i, \nabla w^i\rangle\dif x+ \int_{\mathbb{T}^d}|\nabla w^i|^2m^i\dif x + \int_{\mathbb{T}^d}\langle\rho^i\nabla u^*, \nabla w^i\rangle \dif x= \int_{\mathbb{T}^d}g^i(x)w^i\dif x. 
\end{split}
\end{align}
Subtracting \eqref{wirhoi} from \eqref{rhoiwi}, we have
\begin{align*}
\begin{split}
&\frac{1}{2}\int_{\mathbb{T}^d}|\nabla w^i|^2(m^i+m^*)\dif x+\int_{\mathbb{T}^d}(F[m^i]-F[m^*])\rho^i\dif x \\
=&-\tau^i\int_{\mathbb{T}^d}\rho^i\dif x - \int_{\mathbb{T}^d}h^i(x)\rho^i\dif x+\int_{\mathbb{T}^d}g^i(x)w^i\dif x. 
\end{split}
\end{align*}
Using \eqref{convhgnd}, $\int_{\mathbb{T}^d}m^*\dif x=\int_{\mathbb{T}^d}m^\infty\dif x=1$, and Claim \ref{cl2mc}, we obtain
\begin{align*}
\lim\limits_{i\rightarrow\infty}\frac{1}{2}\int_{\mathbb{T}^d}|\nabla w^i|^2(m^i+m^*)\dif x+\int_{\mathbb{T}^d}(F[m^i]-F[m^*])(m^i-m^*)\dif x = 0.  
\end{align*}
Hence, we have
\begin{align*}
\frac{1}{2}\int_{\mathbb{T}^d}|\nabla w^i|^2(m^\infty+m^*)\dif x+\int_{\mathbb{T}^d}(F[m^i]-F[m^*])(m^i-m^*)\dif x \rightarrow 0. 
\end{align*}
Since $m^\infty\geq 0$ and $m^*\geq 0$, we get
\begin{align*}
	\lim_{i\rightarrow 0}\int_{\mathbb{T}^d}(F[m^i]-F[m^*])(m^i-m^*)\dif x = 0. 
\end{align*}
By Assumption \ref{hypF} and Claim \ref{cl2mc}, 
\begin{align*}
	\int_{\mathbb{T}^d}(F[m^\infty]-F[m^*])(m^\infty-m^*)\dif x = 0. 
\end{align*}
Thus, we have $m^\infty=m^*$. 
By Claim \ref{cl2mc}, we conclude that, up to a sub-sequence, $m^i$ converges to $m^*$ in $H^{1}(\mathbb{T}^d)$. Next, we show the convergence of $u^i$. 

Multiplying the first equation in \eqref{diffndeq} by $w^i$ and integrating by parts, we get
\begin{align}
	\label{bdwi}
	\begin{split}
		\int_{\mathbb{T}^d}|\nabla w^i|^2=&-\int_{\mathbb{T}^d}\bigg\langle\frac{\nabla u^i+\nabla u^*}{2}, \nabla w^i w^i\bigg\rangle \dif x+\int_{\mathbb{T}^d}(F[m^i]-F[m^*])w^i\dif x\\
		&+\tau^i\int_{\mathbb{T}^d}w^i\dif x + \int_{\mathbb{T}^d}h^i(x)w^i\dif x.
	\end{split}
\end{align}
\begin{comment}
Then, we integrate the second equation in \eqref{diffndeq} multiplied with $\rho^i$ and get
\begin{align}
	\label{bdrhoi}
	\begin{split}
		\int_{\mathbb{T}^d}|\nabla \rho^i|^2\dif x = - \int_{\mathbb{T}^d}\langle m^i\nabla w^i, \nabla \rho^i\rangle \dif x - \int_{\mathbb{T}^d}\langle\rho^i\nabla u^*, \nabla \rho^i\rangle\dif x + \int_{\mathbb{T}^d}g^i(x)\rho^i\dif x.
	\end{split}
\end{align}
Hence, passing $i$ to infinity, we get from \eqref{bdwi} that $\nabla w^i$ converges to $0$ in $L^2(\mathbb{T}^d)$. Then, using \eqref{bdrhoi}, $\nabla \rho^i$ converges to $0$ in $L^2(\mathbb{T}^d)$. Hence,  $u^i$ converges to $u^*$ in $H^{1}(\mathbb{T}^d)$, $m^i$ converges to $m^*$ in $H^{1}(\mathbb{T}^d)$. Finally, using \eqref{bdndDui},  $\overline{H}^i$ converges to $\overline{H}^*$ in $\mathbb{R}$. 
\end{comment}
Hence, passing $i$ to infinity, we get from \eqref{bdwi} that $\nabla w^i$ converges to $0$ in $L^2(\mathbb{T}^d)$. Hence,  $u^i$ converges to $u^*$ in $H^{1}(\mathbb{T}^d)$. Finally, using \eqref{bdndDui},  $\overline{H}^i$ converges to $\overline{H}^*$ in $\mathbb{R}$. 
\end{proof}

\end{document}